\documentclass[a4paper,12pt]{article}
\usepackage{pdfsync}
\usepackage{amsfonts}
\usepackage{amssymb}
\usepackage{amsmath}
\usepackage{mathptmx}
\usepackage{cancel}
\usepackage{pstricks}
\usepackage{psfrag}
\usepackage{mathrsfs}
\usepackage{graphicx}
\usepackage{pgf,tikz}
\usepackage{algorithm}
\usepackage{algpseudocode}
\usepackage{soul}
\def\be{\begin{equation}}
\def\ee{\end{equation}}
\def\bea{\begin{eqnarray}}
\def\eea{\end{eqnarray}}
\def\beann{\begin{eqnarray*}}
\def\eeann{\end{eqnarray*}}

\newcommand{\rank}{{\rm rank}}

\def\ns{\hspace{-1mm}}
\newcommand{\real}{{\mathbb{R}}}

\def\spacingset#1{\def\baselinestretch{#1}\small\normalsize}
\spacingset{1.30}

\newtheorem{lemma}{Lemma}
\newtheorem{theorem}{Theorem}
\newtheorem{remark}{Remark}
\newtheorem{corollary}{Corollary}

\newtheorem{definition}{Definition}
\newtheorem{problem}{Problem}
\newtheorem{example}{Example}[section]
\newtheorem{assumption}{Assumption}[section]

\def\be{\begin{equation}}
\def\ee{\end{equation}}
\def\bea{\begin{eqnarray}}
\def\eea{\end{eqnarray}}
\def\beann{\begin{eqnarray*}}
\def\eeann{\end{eqnarray*}}

\def\ns{\hspace{-1mm}}
\def\im{{\rm im}\ }
\def\proof{\noindent{\bf{\em Proof:}\ \ }}
\def\QED{\mbox{\rule[0pt]{1.5ex}{1.5ex}}}
\def\endproof{\hspace*{\fill}~\QED\par\endtrivlist\unskip}

\newcommand{\ima}{\operatorname{im}}
\newcommand{\diag}{\operatorname{diag}}
\newcommand{\defi}{\stackrel{\text{\tiny def}}{=}}

\definecolor{Royalblue}{cmyk}{1,0.30,0.2,0.2}

\newcommand{\complex}{{\mathbb{C}}}

\def\gA{{\cal A}}
\def\gB{{\cal B}}

\def\gD{{\cal D}}
\def\gE{{\cal E}}

\def\gJ{{\cal J}}

\def\gL{{\cal L}}

\def\gR{{\cal R}}

\def\gT{{\cal T}}
\def\gU{{\cal U}}
\def\gV{{\cal V}}

\def\gX{{\cal X}}
\def\gY{{\cal Y}}
\def\gZ{{\cal Z}}

\def\bmat{\left[ \begin{array}}
\def\emat{\end{array} \right]}

\def\bmat{\left[ \begin{array}}
\def\emat{\end{array} \right]}
\def\bsmat{\left[ \begin{smallmatrix}}
\def\esmat{\end{smallmatrix} \right]}

\def\l{{\lambda}}
\def\gA{{\cal A}}

\def\gB{{\cal B}}
\def\gU{{\cal U}}
\def\gL{{\cal L}}
\def\gR{{\cal R}}

\def\gV{{\cal V}}

\def\gT{{\cal T}}

\def\gX{{\cal X}}
\def\i{{i}}

\newcommand{\spanR}{\operatorname{span}}
\def\field{\mathbb{K}}
\def\gsR{{\cal R}^{\star}}
\def\gsV{{\cal V}^{\star}}

\begin{document}
\begin{titlepage}
\title{\vspace{-5mm}
A structural approach to state-to-output decoupling\vspace{10mm}}
\author{Emanuele~Garone$^\dagger$, \and Lorenzo Ntogramatzidis$^\star$ \and Fabrizio Padula$^\star$}
\date{
     $^\dagger${\small Department SAAS, Universit{\`e} Libre de Bruxelles, Brussels, Belgium \\[-3pt]
     E-mail: {\tt egarone@ulb.ac.be} }\\[2pt]
$^\star${\small Department of Mathematics and Statistics,\\[-2pt]
         Curtin University, Perth (WA), Australia \\[-3pt]
        E-mail:  {\tt \{L.Ntogramatzidis,Fabrizio.Padula\}@curtin.edu.au}}\\[8pt]
        }%
\thispagestyle{empty} \maketitle \thispagestyle{empty}
\begin{abstract}%
In this paper, we address a general eigenstructure assignment problem where the objective is to distribute the closed-loop modes over the components of the system outputs in such a way that, if a certain mode appears in a given output, it is unobservable from any of the other output components. By linking classical geometric control results with the theory of combinatorics, we provide necessary and sufficient conditions for the solvability of this problem, herein referred to as {\em state-to-output decoupling}, under very mild assumptions. 
We propose solvability conditions expressed in terms of the dimensions of suitably defined controlled invariant subspaces of the system. In this way, the  solvability of the problem can be evaluated \emph{a priori}, in the sense that it is given in terms of the problem/system data. Finally, it is worth mentioning that the proposed approach is constructive, so that when a controller that solves the problem indeed exists, it can be readily computed by using the machinery developed in this paper.
\end{abstract}

\begin{center}
\begin{minipage}{14.2cm}
\vspace{2mm}
\end{minipage}
\end{center}
\thispagestyle{empty}
\end{titlepage}
%
%%%%%%%%%%%%%%%%%%%%%%%%%%%%%%%%%%%%%%%%%%%%%%%%%%%%%%%%%%%%%%%%%%%%%%%%%%%%%

\section{Introduction}
\label{secintro}
The problem of mode allocation/distribution in the outputs of multiple-input multiple-output (MIMO) systems is central in systems and control theory.
The pioneering paper \cite{Moore-76} was the first to highlight the fact that this problem is, in essence, a problem of eigenstructure assignment for the closed-loop. In other words, imposing a certain distribution of closed-loop modes on the output components of a MIMO system is equivalent to suitably assigning the closed-loop eigenvalues as well as the corresponding eigenvectors.
This idea has been exploited in a variety of contexts, raging from 
fault diagnosis and isolation \cite{Chen-P-99} to aircraft control \cite{Mudge-P-88}, and extending also to areas such as matrix interpolation
\cite{Antsaklis-G-93}, active suppression of vibrations \cite{Mottershead-TJR-08} and design of autopilots \cite{Doll-LFM-01}.

In recent years, the eigenstructure assignment of \cite{Moore-76} has found new applications in the area of tracking control for MIMO systems. In \cite{SN-AUT-10}, a new control methodology was presented to tackle the problem of tracking a vector of step functions with no overshoot; the main idea behind that strategy, which has been very recently developed in \cite{NTSF-TAC-15} for the case of monotonic tracking, is to ensure that every component of the tracking error comprises a single closed-loop mode independently from the initial condition. This property was proved in \cite{NTSF-TAC-15} to be necessary and sufficient to guarantee that the system response is monotonic %(and, consequently, both non-overshooting and non-undershooting)
 from any initial state of the system. %In \cite{NTSF-TAC-15} the converse implication result was also proved: the only way to achieve monotonic tracking from any initial state is to distribute one closed-loop mode per output component. %Importantly, a set of necessary and sufficient conditions, given in terms of the problem data, was presented in that paper for the solvability of this problem. 

In this paper, for the first time in the literature, we provide necessary and sufficient conditions for the solvability of the eigenstructure assignment problem of an arbitrary number of closed-loop modes per output component under virtually no assumptions. In particular, this paper addresses the problem of ensuring that each output component comprises a preassigned set of closed-loop modes, possibly including the invariant zeros of the system. 
In order to prove this result, a new framework is introduced which links classical results of geometric control theory  \cite{Wonham-85,Basile-M-92,Trentelman-SH-01,Chen-LS-04,Hamilton-B-12} with the theory of combinatorics 
\cite{Rado-42,Ore-55,Mirsky-67} that enables the solvability conditions to be expressed in terms of specific and easily computable controlled invariant subspaces which are completely defined in terms of the parameters of the problem.
It is also worth mentioning that the methodology developed in this paper is constructive in nature, because it allows to immediately compute the suitable feedback matrix that solves the problem whenever such matrix exists.

We also establish that the above mentioned eigenstructure assignment problem can be reformulated as the problem of rendering the autonomous system associated with the system at hand equivalent, in a system-theoretic sense, to a set of decoupled autonomous systems. 
Hence, the eigenstructure assignment problem considered here is equivalent to finding a controller that achieves a decoupling between the state and the output; for this reason, hereafter this property will be referred to as {\em state-to-output decoupling}.

This property appears to be a particularly important feature of the problem considered in this paper. For example, it links with some problems of security in large-scale complex systems, see \cite{Mo-S-15}, \cite{Ugrinovskii-L-11} and the references cited therein. Indeed, the idea behind the state-to-output decoupling is the fact that, from each output component, only a certain subset of the system modes is observable; this means that, in the context of secure control, an attacker needs to have access to the information originating from all the sensors in order to reconstruct the state of the system. In this way, if the information coming from a sensor is compromised, it is not possible to reconstruct the entire state of the system, but only a portion of it.

Furthermore, %This consideration highlights that, in addition to providing a constructive solution to an open problem, 
the machinery developed in this paper can be used as a building block to solve a variety of other important control problems. For instance it allows to drastically reduce the computational burden in the calculation of the matrix exponential of the closed-loop system. Other applications arise in the context of the fault detection and non-interacting control literature, see e.g. \cite{Wahrburg-A-12}. Indeed, a number of those problems, for which only {\em a posteriori} solvability conditions are currently available in the literature, can be %for which our results provide a closed-form solution
viewed as reformulations of the state-to-output decoupling problem.
Thus, the methodology provided in this paper provides a solution to the aforementioned problems in terms of the problem data, which is therefore {\em a priori}.

 %, for which a closed-form solution was not available.

Among the problems that can be dealt with as state-to-output decoupling, one that stands out is the monotonic tracking control for those systems for which the necessary and sufficient conditions of \cite{NTSF-TAC-15} do not hold.
Indeed, such systems may still exhibit a non-overshooting and non-undershooting response, and the shape and size of the set of initial conditions for which this is the case depends on the number of closed-loop modes appearing in each output component. %: fewer modes result in a larger region of initial states for which both overshooting and undershooting can be avoided. 
Moreover, in practice it is not always necessary to impose a monotonic response in each output component. %; indeed, in some situations it may be acceptable to require a monotonic response only from a proper subset of the output components. 
These two fundamental relaxations of the problem dealt with in \cite{NTSF-TAC-15} require a richer machinery, which is the one developed in this paper.

%These two fundamental relaxations of the problem dealt with in \cite{NTSF-TAC-15} require a richer machinery, which is developed in this paper. More specifically, in this paper we establish necessary and sufficient conditions for three classes of problems in which, essentially, we require only certain modes to appear on each output component independently of the initial condition. The distinction among these problems is due to the role that the invariant zeros play in the allocation of the closed-loop eigenvalues and eigenvectors. 
 %New problems arise when generalizing the approach of \cite{NTSF-TAC-15} to more than one closed-loop mode per output component; a new challenging aspect, among the others, is the fact that complex modes are now allowed in the closed-loop spectrum. 

 The concept of state-to-output decoupling introduced in this paper is also relevant in the context of constrained distributed control, involving a number of subsystems with shared constraints and dynamics. Generally speaking, the prediction obtained using e.g. a {model predictive control} (MPC) scheme \cite{Bemporad-M-94} or a distributed command governor architecture \cite{Garone-KD-16} cannot neglect the influence that each subsystem has on the other subsystems. Hence, 
  even though the decoupling of the dynamics of these subsystems does not completely overcome the issue (because of the presence of the constraints which remain in general coupled), the technique presented here leads to simpler and more efficient distributed control strategies (see e.g. \cite{Casavola-GT-14}).

Finally we want to mention that an important by-product of the results established in this paper is the identification of a self-bounded output-nulling subspace, herein denoted by ${\cal L}$, which has interesting %and, to the best of our knowledge, 
system-theoretic properties that, to the best of our knowledge, have never been investigated, and which
plays a key role in the solution of the state-to-output decoupling problem.

%Another 

%Interestingly, the systematic approach of \cite{NTSF-TAC-15} shows that 
% objectives of
%achieving a rapid settling time, while at the same time avoiding
%overshoot and undershoot, are not necessarily competing objectives in the
%controller design.

%An example is the constrained distributed control; consider the case of several subsystems both with constraints and dynamics in common: the control problem in this case is very complicated (MPC: a subsystem makes a prediction but it cannot neglect the influence it has on the dynamics of the other subsystem). If the dynamics are decoupled, even if the constraints are still in common, there are techniques that enable the problem to be solved in a way that the loss of performance due to the decoupling is negligible. Some examples are MPC and the distributed command governor.

{\bf Notation.} The image and the kernel of matrix $A$ are denoted by $\ima\,A$ and $\ker\,A$, respectively. The Moore-Penrose pseudo-inverse of $A$ is denoted by $A^\dagger$, and $A^{\scriptscriptstyle -R}$ denotes a right inverse of $A$ when $A$ is right invertible. When $A$ is square, we denote by $\sigma(A)$ the spectrum of $A$. If $A: \gX \longrightarrow \gY$ is a linear map and if $\gJ \subseteq \gX$, the restriction of the map $A$ to $\gJ$ is denoted by $A\,|\gJ$. If $\gX=\gY$ and $\gJ$ is $A$-invariant, the eigenstructure of $A$ restricted to $\gJ$ is denoted by $\sigma\,(A\,| \gJ)$. If $\gJ_1$ and $\gJ_2$ are $A$-invariant subspaces and $\gJ_1\,{\subseteq}\,\gJ_2$, the mapping induced by $A$ on the quotient space $\gJ_2 / \gJ_1$ is denoted by $A\,| {\gJ_2}/{\gJ_1}$, and its spectrum is denoted by $\sigma\,(A\,| {\gJ_2}/{\gJ_1})$.
%The symbol $\oplus$ stands for the direct sum of subspaces.  The symbol $\uplus$ denotes union with any common elements repeated.
Given a map $A: \gX \longrightarrow \gX$ and a subspace $\gB$ of $\gX$, we denote by $\langle A\,|\, \gB \rangle$ the smallest $A$-invariant subspace of $\gX$ containing $\gB$. %The symbol $\i$ stands for the imaginary unit, i.e., $\i=\sqrt{-1}$. 
Given a complex matrix $M$, the symbols $\overline{M}$ and ${M}^\ast$ denote the conjugate and the conjugate transpose of $M$, respectively. Moreover, we denote by $M_i$ its $i$-th row and by $M^j$ its $j$-th column, respectively.
 Given a finite set $S$, the symbol $2^S$ denotes the power set of $S$, while ${\rm {card}}(S)$ stands for the cardinality of $S$.
%We denote by $\complex_g$ the set of complex numbers with strictly negative real part, or an arbitrary self conjugate region of $\complex$ to the left of the imaginary axis.

\section{Problem Statements}
In what follows, whether the underlying system evolves in continuous or discrete time is irrelevant and, accordingly, the time index set of any signal is denoted by $\mathbb{T}$, on the understanding that this represents either $\real^+$ in the continuous time or $\mathbb{N}$ in the discrete time. The symbol $\complex_g$ denotes either the open left-half complex plane $\complex^-$ in the continuous time or the open unit disc $\complex^\circ$ in the discrete time.\footnote{The results developed in this paper continue to hold even when $\complex_g$ is an arbitrary self conjugate region of $\complex$ to the left of the imaginary axis in the continuous time or inside the open unit circle in the discrete time.} Likewise, $\real_g$ denotes the set of strictly negative real numbers in the continuous time or the real numbers in $(-1,1)$ in the discrete time.
Consider the LTI system $\Sigma$ governed by
 \bea
 \Sigma: \
 \left\{ \begin{array}{rcl}
\gD\,{x}(t) \ns&\ns \!\! = \!\!  \ns&\ns A\,x(t)+B\,u(t),\;\;\;\; x(0)=x_0, \label{syseq1}\hfill\cr
 y(t) \ns&\ns  \!\! = \!\!  \ns&\ns C\,x(t)+D\,u(t),\hfill \end{array} \right.
 \label{sys}
 \eea
where, for all $t \in \mathbb{T}$, $x(t) \in \gX=\real^n$ is the
state, $u(t) \in \gU=\real^m$ is the control input, $y(t) \in \gY=\real^p$ is the output, and $A$, $B$, $C$ and $D$ are appropriate
dimensional constant matrices. The operator $\gD$ denotes either the time derivative in the continuous time, i.e., $\gD x(t)=\dot{x}(t)$, or the unit time shift in the discrete time, i.e., $\gD x(t)=x(t+1)$. Let the system $\Sigma$ described by (\ref{syseq1}) be identified with the quadruple $(A,B,C,D)$.
The following standing assumptions ensures that any given constant reference target $r(t)=\bar{r} \in \real^p$ can be tracked from any initial condition $x_0 \in \gX$:\\[-3mm]
 \begin{assumption}
 \label{Ass1}
System $\Sigma$ is right invertible and stabilizable. Moreover, $\Sigma$ has no invariant zeros at the origin in the continuous time or at 1 in the discrete time.\\[-3mm]
\end{assumption}

%We denote by $\gZ$ the set of invariant zeros of $\Sigma$.
 Let us consider the state-feedback control law
 \be
 u(t) = F\, x(t) + G\,r(t),
 \label{ulaw}
 \ee
 where $F$ is a stabilizing feedback, i.e., $\sigma(A+B\,F)\subset \complex_g$, and $G$ is a right inverse of the static gain of the quadruple $(A+B\,F,B,C+D\,F,D)$, i.e.,%\footnote{The inverse of $A+B\,F$  if the system is right invertible and without zeros at $0$.}
 \[
 G=-\bigl((C+D\,F)\,(A+B\,F)^{-1}\,B+D\bigr)^{\scriptscriptstyle -R} \quad \text{and} \quad  G=\bigl((C+D\,F)\,\bigl(I-(A+B\,F)\bigr)^{-1}\,B+D\bigr)^{\scriptscriptstyle -R}
  \]
  in the continuous and discrete time, respectively. Notice that a right inverse always exists in view of Assumption \ref{Ass1}, and it can be computed for example as a Moore-Penrose pseudo-inverse.
Applying (\ref{ulaw}) to (\ref{sys}), we obtain the closed-loop
system

 \be
 \Sigma_{\sc F,G}: \
\left\{ \begin{array}{rcl}
\gD\,{x}(t) \ns&\ns  \!\! = \!\!  \ns&\ns (A+B\,F)\,x(t)+B\,G\,r(t),\;\;\quad x(0)=x_0,\hfill\cr
 y(t) \ns&\ns  \!\! = \!\!  \ns&\ns (C+D\,F)\,x(t)+ D\,G\,r(t). \hfill \end{array} \right.
 \label{syschom}
 \ee
Since $r(t)=\bar{r}$ is constant, with a change of coordinates (\ref{syschom}) can be written in terms of the error $\epsilon \defi y-r$ as
 \be
 \Sigma_{\sc F,G}: \
\left\{ \begin{array}{rcl}
\gD\,{\xi}(t) \ns&\ns  \!\! = \!\!  \ns&\ns (A+B\,F)\,\xi(t),\;\;\quad \xi(0)=\xi_0,\hfill\cr
\epsilon(t) \ns&\ns  \!\! = \!\!  \ns&\ns (C+D\,F)\,\xi(t). \hfill \end{array} \right.
 \label{syschom1}
 \ee

%In this paper we will deal with two classes of problems. The first class will be referred to as {\em output decoupling}, while the second will be called {\em input/state to output decoupling}.

This paper deals with the problem of determining the state feedback matrix $F$ 
for (\ref{syschom1}) such that each output component comprises a number of closed-loop modes that are unobservable from any other output component. This problem will be referred to as {\em state-to-output decoupling}.
%Roughly speaking,

\begin{definition}{\sc [State-to-Output Decoupling]}\\
\label{def1}
 We say that a feedback matrix $F$ in (\ref{syschom1}) achieves {\em state-to-output decoupling} if, when $r(t)$ is constant, the error in (\ref{syschom1}) can be written for any initial condition as
\[
\epsilon(t)=\bmat{c} 
\beta_{1,1}\,e^{\lambda_{1,1}\,t}+\ldots+\beta_{1,\nu_1}\,e^{\lambda_{1,\nu_1}\,t}\\
\vdots\\
\beta_{p,1}\,e^{\lambda_{p,1}\,t}+\ldots+\beta_{p,\nu_p}\,e^{\lambda_{p,\nu_p}\,t}
\emat \quad \text{and} \quad 
\epsilon(t)=\bmat{c} 
\beta_{1,1}\,\lambda_{1,1}^t+\ldots+\beta_{1,\nu_1}\,\lambda_{1,\nu_1}^t\\
\vdots\\
\beta_{p,1}\,\lambda_{p,1}^t+\ldots+\beta_{p,\nu_p}\,\lambda_{p,\nu_p}^t
\emat
\]
in the continuous and discrete time, respectively, 
where $\l_{i,j}$ are the observable closed-loop eigenvalues 
%where we recall that the mapping $x_0\mapsto \{\alpha_{i,j}\}$ is surjective.
and
\begin{itemize}
\item if $\l_{i,j}$ is real, the coefficient $\beta_{i,j}$ can be made arbitrary by choosing a suitable initial state $\xi_0$;
\item if $\l_{i,j}$ is complex, there exists $k$ such that $\l_{i,k}=\overline{\l}_{i,j}$, and $\beta_{i,k}=\overline{\beta}_{i,j}$ where $\beta_{i,j}$ can be made arbitrary by choosing a suitable initial state $\xi_0$.
%either the real or the imaginary part of 
%$\beta_{i,j}$ can be made arbitrary
\end{itemize}
\end{definition}
In Definition \ref{def1}, for clarity we have distinguished the case where $\l_{i,j}$
is real from the case where $\l_{i,j}$ is complex. The two cases can be captured together by saying that for every $\l_{i,j}$ either the real or the imaginary part of the corresponding $\beta_{i,j}$ can be made arbitrary.

Note that in Definition \ref{def1} it has been implicitly assumed that no Jordan chains appear in the observable closed-loop eigenstructure. Indeed, in this paper we make the standing assumption that no Jordan chains are allowed in the closed-loop eigenstructure. The reason for this choice, together with a discussion of the technicalities to overcome this apparent limitation, will be detailed in Remark \ref{pizzardonejordan}.

%\textcolor{red}{(Notice that in the complex case we must have $\l_{i,j}\in \complex$ $\Rightarrow$ there exists $\l_{i,k}$ such that $\l_{i,k}=\l_{i,j}^\ast$ for some $k \in \{1,\ldots,\nu_i\}$, and therefore $\beta_{i,k}=\beta_{i,j}^\ast$; non e' piu vero che i beta sono arbitrari: posso rendere arbitraria la parte real e e quella immaginaria di $\beta_{i,k}$ e quindi anche quella di $\beta_{i,j}$, ma non li posso rendere arbitrari tra loro; il numero di gradi di liberta' si conserva. Nel Remark che segue, ipotizzare che i $\l$ sono reali for the sake of simplicity.)}

\begin{remark}
{\em
The requirement that the coefficients $\beta_{i,j}$ can be made arbitrary guarantees that each $\l_{i,j}$ defines the closed-loop dynamics along a different direction of the state space. In other words, each $\l_{i,j}$ is associated with a different closed-loop eigenvector. This implies that if two closed-loop eigenvalues are identical, they describe the dynamics along different directions, and therefore they correspond to two different closed-loop modes.
This consideration can be formalized as follows. 
%The reason for requiring the coefficients $\beta_{i,j}$ to be arbitrary in Definition \ref{def1} is the following. 
The solution of $\Sigma_{\sc F,G}$, say in the continuous time, can be written as 
$\epsilon(t) = (C+D\,F)\,e^{(A+B\,F)\,t}\,\xi_0$. 
Assume for simplicity that $A+B\,F$ 
 has $n$ real eigenvalues $\l_{0,1},\ldots,\l_{0,\nu_0},
 \l_{1,1},\ldots,\l_{1,\nu_1},\ldots,\l_{p,1},\ldots,\l_{p,\nu_p}$ (with $\nu_0=n-\nu_1-\ldots-\nu_p$) associated with the 
 linearly independent real eigenvectors $v_{0,1},\ldots,v_{0,\nu_0},
 v_{1,1},\ldots,v_{1,\nu_1},\ldots,$ $v_{p,1},\ldots,v_{p,\nu_p}$, so that we can write 
$\xi_0=\sum_{i=0}^p \sum_{j=1}^{\nu_i} \alpha_{i,j}\,v_{i,j}$ for suitable $\alpha_{i,j} \in \real$. 
Recall that $(A+B\,F)\,v_{i,j}=\lambda_{i,j}\,v_{i,j}$ ($i\in \{0,\ldots,p\}$ and $j \in \{1,\ldots,\nu_i\}$) implies $e^{(A+B\,F)\,t}\,v_{i,j}=e^{\lambda_{i,j}\,t}\,v_{i,j}$,  
which means that
\beann
\xi(t) \ns&\ns = \ns&\ns  e^{(A+B\,F)\,t}\,\xi_0=\sum_{i=0}^p \sum_{j=1}^{\nu_i} \alpha_{i,j}\,e^{\l_{i,j}\,t}\,v_{i,j}
\\
\epsilon(t) \ns&\ns = \ns&\ns (C+D\,F)\,\xi(t)=
\sum_{i=0}^p \sum_{j=1}^{\nu_i} \alpha_{i,j}\,e^{\l_{i,j}\,t}\,(C+D\,F)\,v_{i,j} = \sum_{i=1}^p \sum_{j=1}^{\nu_i} \alpha_{i,j}\,e^{\l_{i,j}\,t}\,(C+D\,F)\,v_{i,j}.
\eeann
Now it is clear that 
$(C_i+D_i\,F)\,v_{i,j} \neq 0$ implies that $(C_{h}+D_{h}\,F)\,v_{i,j}=0$ for all $h\neq i$. 
 Indeed, if we have $(C_i+D_i\,F)\,v_{i,j}\neq 0$ and $(C_h+D_h\,F)\,v_{i,j}\neq 0$ for some $h \neq i$, then in $\epsilon_i$ we would have the component 
$\beta_{i,j}\,e^{\lambda_{i,j}\,t}=
\alpha_{i,j}\,(C_i+D_i\,F)\,v_{i,j}\,e^{\lambda_{i,j}\,t}$ and in $\epsilon_h$ we would have the component 
$\beta_{h,j}\,e^{\lambda_{i,j}\,t}=\alpha_{i,j}\,(C_h+D_h\,F)\,v_{i,j}\,e^{\lambda_{i,j}\,t}$, which are proportional. Thus, $\beta_{i,j}$ and $\beta_{h,j}$ cannot be made arbitrary by choosing a suitable initial condition (which would affect only $\alpha_{i,j}$).
}
\end{remark}

%\textcolor{red}{
%\begin{lemma}
%If the coefficients $\beta_{i,j}$ can be made arbitrary choosing an appropriate initial condition $\xi_0$, each closed-loop eigenvalue is observable from only one output, i.e., 
%\[
%(C_i+D_i\,F)\,v_{i,j} \neq 0 \quad \Rightarrow \quad (C_{k}+D_{k}\,F)\,v_{i,j}=0 \;\; \forall \,k\neq i.
%\] 
%\end{lemma}
%%
%\proof By contradiction suppose $(C_i+D_i\,F)\,v_{i,j} \neq 0$ and $(C_k+D_k\,F)\,v_{i,j} \neq 0$ with $i\neq k$.  
%%For some $h\in \{1,\ldots,\nu_i\}$ and $s\in \{1,\ldots,\nu_k\}$
%Thus,% we have
%\beann
%\beta_{i,j}=\alpha_{i,j}\,(C_i+D_i\,F)\,v_{i,j} \\
%\beta_{k,j}=\alpha_{i,j}\,(C_k+D_k\,F)\,v_{i,j}
%\eeann
%%(il modo $e^{\l_j\,t}$ corrispondera' ad un certo $e^{\l_{i,j}\,t}$: stiamo usando due notazioni diverse per indicare gli eig della catena chiusa).
%Thus, $\beta_{i,j}$ and $\beta_{k,j}$ are proportional, and cannot be made arbitrary by choosing a suitable initial condition (which would affect only $\alpha_{i,j}$).
%}

%Referring again to the continuous time case for the sake of argument, it is easy to see that 
The following result shows that the state-to-output decoupling problem can be reformulated as the problem of existence of $p$ single-output systems $\Sigma_1,\ldots,\Sigma_p$ such that $\Sigma$ is equivalent (in a system-theoretic sense) to the Cartesian product of $\Sigma_1,\ldots,\Sigma_p$.

\begin{theorem}
\label{thequiv}
The state-to-output decoupling problem is equivalent to the existence
of matrices $A_1,\ldots,A_p$ and row vectors $C_1,\ldots,C_p$ such that:
\begin{itemize}
\item to any state $x\in \gX$ it is possible to associate $x_1\in \gX_1$, $x_2\in \gX_2$, $\ldots$, $x_p\in \gX_p$ such that the response of $\Sigma_{\sc F,G}$ from the initial condition $x_0$
% and the initial time $t_0=0$ 
with the reference $r=0$ coincides with the vectors of the responses $\bsmat y_1(\cdot) \\[-2mm] \vdots \\ y_p(\cdot) \esmat$ obtained from $(A_i,C_i)$, the initial condition $x_{i,0}$, i.e.,
% and the initial time $t_0=0$
\bea
(C+D\,F) e^{(A+B\,F)\,t} x_0=
\bsmat
C_1 e^{A_1 \,t} x_{1,0} \\[-2mm]
\vdots \\
C_p e^{A_p \,t} x_{p,0} 
\esmat\quad \forall \,t\ge 0;\label{equivalence0}
\eea
\item conversely, for any choice of initial states $x_{1,0}\in \gX_1$, $x_{2,0}\in \gX_2$, $\ldots$, $x_{p,0}\in \gX_p$ there exists an initial state $x_0\in \gX$ of $\Sigma_{F,G}$ such that (\ref{equivalence0}) holds true for $r=0$.
\end{itemize}
\end{theorem}
\proof
Consider the continuous time for the sake of argument. If 
\[
(C+D\,F) \,e^{(A+B\,F)\,t} x_0=
\bsmat
C_1 e^{A_1 \,t} x_{1,0} \\[-2mm]
\vdots \\
C_p e^{A_p \,t} x_{p,0} 
\esmat
\]
and $\sigma(A_i)=\{\lambda_{i,1},\ldots,\lambda_{i,\nu_i},\lambda_{i,\nu_i+1},\ldots,\lambda_{i,n_i}\}$ (where in general $n_i\ge \nu_i$ since $\Sigma_i$ needs not be in minimal form), there exists an invertible matrix $T_i$ such that
$T_i^{-1}\,A_i\,T_i=\diag\{\lambda_{i,1},\ldots,\lambda_{i,n_i}\}$. 
We find
\beann
(C+D\,F)\, e^{(A+B\,F)\,t} x_0 \ns&\ns = \ns&\ns
\bsmat
C_1\,T_1\,\diag\{\lambda_{1,1},\ldots,\lambda_{1,n_1}\}\,T_1^{-1}\, x_{1,0} \\[-1mm]
\vdots \\
C_p \,T_p\,\diag\{\lambda_{p,1},\ldots,\lambda_{p,n_p}\}\,T_p^{-1}\, x_{p,0} 
\esmat=
\bsmat
[\begin{array}{cccc} c_{1,1} & \ldots & c_{1,n_1} \end{array}] \bsmat e^{\lambda_{1,1}\,t}\,z_{1,0,1} \\[-1mm]
\vdots \\ 
e^{\lambda_{1,n_1}\,t}\,z_{1,0,n_1} 
\esmat
\\
\vdots \\
[\begin{array}{cccc} c_{p,1} & \ldots & c_{p,n_p} \end{array}] \bsmat e^{\lambda_{p,1}\,t}\,z_{p,0,1} \\[-1mm]
\vdots \\ 
e^{\lambda_{p,n_p}\,t}\,z_{p,0,n_p}
\esmat
\esmat \\
 \ns&\ns = \ns&\ns 
 \bmat{c} 
\beta_{1,1}\,e^{\lambda_{1,1}\,t}+\ldots+\beta_{1,n_1}\,e^{\lambda_{1,n_1}\,t}\\
\vdots\\
\beta_{p,1}\,e^{\lambda_{p,1}\,t}+\ldots+\beta_{p,n_p}\,e^{\lambda_{p,n_p}\,t}
\emat= \bmat{c} 
\beta_{1,1}\,e^{\lambda_{1,1}\,t}+\ldots+\beta_{1,\nu_1}\,e^{\lambda_{1,\nu_1}\,t}\\
\vdots\\
\beta_{p,1}\,e^{\lambda_{p,1}\,t}+\ldots+\beta_{p,\nu_p}\,e^{\lambda_{p,\nu_p}\,t}
\emat
\eeann
where $C_i\,T_i=[\begin{array}{cccc} c_{i,1} & \ldots & c_{i,\nu_i} \end{array}]$ with $c_{i,\nu_i+1}=\ldots=c_{i,n_i}=0$ since $\l_{i,\nu_i+1},\ldots,\l_{i,n_i}$ are unobservable, $z_{i,0}=T_i^{-1}\,x_{i,0}=\bsmat z_{i,0,1}\\[-1mm] \vdots \\ z_{i,0,n_i}\esmat$, and
where $\beta_{i,j}=c_{i,j}\,z_{i,0,j}$. The same steps can be reversed to prove the opposite implication.
%The same argument can be used in the case of non-trivial Jordan structures.
\endproof

In this paper we deal with three specific problems of state-to-output decoupling. Before proceeding with their definition, we recall that the Rosenbrock matrix is defined as the matrix pencil 
\bea
\label{ros}
P_{\scriptscriptstyle \Sigma}(\lambda) \defi \bmat{cc} A-\lambda\,I & B \\ C & D \emat
\eea 
in the indeterminate $\lambda \in \complex$.
The invariant zeros of $\Sigma$ are the values of $\lambda \in \complex$ for which the rank of $P_{\scriptscriptstyle \Sigma}(\lambda)$ is strictly smaller than its normal rank, see \cite{Aling-S-84}.
%and the invariant zero structure is given by the zeros, multiplicity included, of the greatest common divisor of the minors of order $n+\min \{m,p\}$ of $P_{\scriptscriptstyle \Sigma}(\lambda)$, see \cite{MacFarlane-K-76,Aling-S-84}. Thus, 
Given an invariant zero $z \in \complex$, the rank deficiency of $P_{\scriptscriptstyle \Sigma}(\lambda)$ at the value $\lambda=z$ is the geometric multiplicity of the invariant zero $z$, and is equal to the number of elementary divisors (invariant polynomials) of $P_{\scriptscriptstyle \Sigma}(\lambda)$ associated with the complex frequency $\lambda=z$. The degree of the product of the elementary divisors of $P_{\scriptscriptstyle \Sigma}(\lambda)$ corresponding to the invariant zero $z$ is the algebraic multiplicity of $z$, see \cite{MacFarlane-K-76}.
Thus, the algebraic multiplicity of an invariant zero in not smaller than its geometric multiplicity. 

In line with our standing assumption on the absence of Jordan chains in the closed-loop eigenstructure, the algebraic and geometric multiplicities of every minimum-phase invariant zero coincide, i.e., the minimum-phase invariant zeros have trivial (i.e., diagonal) Jordan form, see Remark \ref{pizzardonejordan}.
%\end{assumption}

%This assumption does not lead to a significant loss of generality. In fact, the case of coincident zeros can be dealt with by using the procedure described in \cite{Ntogramatzidis-14}.

 The set of invariant zeros of $\Sigma$ is denoted by $\gZ$, and the set of the minimum-phase invariant zeros is denoted by $\gZ_g \defi \gZ \cap \complex_g$.

We now present the three main problems that we address in this paper: they all deal with the issue of achieving tracking with state-to-output decoupling. In the first problem, the number of observable modes that are visible from each output is fixed, and these modes do not coincide with the minimum-phase invariant zeros of the system. The second problem differs from the first only by the fact that minimum-phase invariant zeros are allowed to be observable eigenvalues for the closed loop. In the last problem,  minimum-phase invariant zeros are still allowed to become observable from the output, but only an upper bound for the number of modes observable from each output is assigned.

Each of these three problems will be in turn divided into three subproblems, labelled as {\bf (A)}, {\bf (B)} and {\bf (C)}: Problem $i$\,A (for $i\in \{1,2,3\}$) refers to the case where both the observable and the unobservable eigenvalues are assigned; Problem $i$\,B is the case where only the observable eigenvalues are assigned. Finally, Problem $i$\,C considers the situation where none of the observable/unobservable eigenvalues are assigned.

We now formulate each problem, along with its subproblems, precisely. We begin with the first problem, which considers the case where each $\epsilon_i$ displays exactly $\nu_i$ modes and the invariant zeros are not selected as observable eigenvalues.

\begin{problem}%{\bf (Output Decoupling + Tracking - 1)}
\label{pro01}
%Let $r$ be constant.
%Let $\{\lambda_{i,j}\}_{j=1,\ldots,\nu_i}$ for $i =1,\ldots,p$ be in $\real^-$.
Determine under which conditions $F$ and $G$ exist such that:
\begin{enumerate}
\item The output asymptotically tracks any constant reference $r$, i.e., if $r(t)=\overline{r}$ for all $t \ge 0$, then $\lim_{t \rightarrow \infty} y(t) = \overline{r}$;
\item State-to-output decoupling is achieved;
\item For each $\epsilon_i$ there are exactly $\nu_i$ observable eigenvalues and they are not invariant zeros of the system, \footnote{The eigenvalues are assumed to be counted with their multiplicities. This is equivalent to saying that the unobservable subspace relative to the output $i$ has dimension $n-\nu_i$. As already pointed out, for the sake of simplicity, we will only consider the case where the geometric and algebraic multiplicities coincide.}
\end{enumerate}
in the following three cases:
\begin{description}
\item{\bf (A)} the eigenvalues that are observable from $\epsilon_i$ are exactly $\{\lambda_{i,j}\}_{j=1,\ldots,\nu_i}$; the unobservable eigenvalues are equal to $\{\lambda_{0,j}\}_{j=1,\ldots,\nu_0}$, where $\nu_0=n-\sum_{i=1}^p \nu_i$;
\item{\bf (B)} the eigenvalues that are observable from $\epsilon_i$ are exactly $\{\lambda_{i,j}\}_{j=1,\ldots,\nu_i}$; the unobservable eigenvalues are not assigned {\em a priori};
\item{\bf (C)} neither the observable nor the unobservable eigenvalues are  assigned {\em a priori}.
\end{description}
\end{problem}

As aforementioned, the second problem deals with the case where each output displays exactly $\nu_i$ eigenvalues, and we allow the selection of invariant zeros in $\complex_g$ as observable eigenvalues.

\begin{problem}%{\bf (Output Decoupling + Tracking - 2)}
\label{pro02}
%Let $r$ be constant.
%Let $\{\lambda_{i,j}\}_{j=1,\ldots,\nu_i}$ for $i =1,\ldots,p$ be in $\real^-$.
Determine under which conditions $F$ and $G$ exist such that:
\begin{enumerate}
\item The output asymptotically tracks any constant reference $\overline{r}$, i.e. if $r(t)=\overline{r},$ $\lim_{t \rightarrow \infty} y(t) = \overline{r}$;
\item Output decoupling is achieved;
\item For each $\epsilon_i$ there are exactly $\nu_i$ observable eigenvalues,
\end{enumerate}
in the following three cases:
\begin{description}
\item{\bf (A)} the eigenvalues observable from $\epsilon_i$ are exactly $\{\lambda_{i,j}\}_{j=1,\ldots,\nu_i}$; the unobservable eigenvalues are equal to $\{\lambda_{0,j}\}_{j=1,\ldots,\nu_0}$, where $\nu_0=n-\sum_{i=1}^p \nu_i$;
\item{\bf (B)} the eigenvalues that are observable from $\epsilon_i$ are exactly $\{\lambda_{i,j}\}_{j=1,\ldots,\nu_i}$; the unobservable eigenvalues are not assigned {\em a priori};
\item{\bf (C)} neither the observable nor the unobservable eigenvalues are  assigned {\em a priori}.
\end{description}
\end{problem}

The last problem considers the case where each output displays at most $\nu_i$ eigenvalues and we allow the selection of the minimum-phase invariant zeros as observable eigenvalues.

\begin{problem}%{\bf (Output Decoupling + Tracking - 3)}
\label{pro03}
%Let $r$ be constant.
%Let $\{\lambda_{i,j}\}_{j=1,\ldots,\nu_i}$ for $i =1,\ldots,p$ be in $\real^-$.
Determine under which conditions $F$ and $G$ exist such that:
\begin{enumerate}
\item The output asymptotically tracks any constant reference $\overline{r}$, i.e. if $r(t)=\overline{r}$, then $\lim_{t \rightarrow \infty} y(t) = \overline{r}$;
\item State-to-output decoupling is achieved;
\item For each $\epsilon_i$ there are at most $\bar{\nu}_i$ observable eigenvalues (so that $\nu_i \le \bar{\nu}_i$),
\end{enumerate}
in the following three cases:
\begin{description}
\item{\bf (A)} the eigenvalues of the closed-loop are $\{\lambda_{i,j}\}_{i=0,\ldots,p,j=1,\ldots,\bar{\nu}_i}$, and the observable eigenvalues from $\epsilon_i$ are a subset $\{\lambda_{i,j}\}_{j=1,\ldots,\nu_i}$ of $\{\lambda_{i,j}\}_{j=1,...,\bar{\nu}_i}$; the unobservable eigenvalues contain $\{\lambda_{0,j}\}_{j=1,\ldots,\bar{\nu}_0}$, where $\bar{\nu}_0=n-\sum_{i=1}^p \bar{\nu}_i$;
\item{\bf (B)} the observable eigenvalues of the closed-loop from $\epsilon_i$ are the subset $\{\lambda_{i,j}\}_{j=1,\ldots,\nu_i}$ of $\{\lambda_{i,j}\}_{j=1,\ldots,\bar{\nu}_i}$;
%, and the observable eigenvalues  from $\epsilon_i$ are a subset of $\{\lambda_{i,j}\}_{j=1,\ldots,\nu_i}$; 
\item{\bf (C)} neither the observable nor the unobservable eigenvalues are  assigned {\em a priori}.
\end{description}
\end{problem}

Notice that in Problem \ref{pro03}{\bf (A)}, some eigenvalues may be hidden from the output, but they still result as eigenvalues of the closed loop.

Notice also that if the eigenvalues which are observable from $\epsilon_i$ are constrained to be at most $\nu_i$, we have the option of hiding as many modes as possible for each output component; hiding more modes than what is strictly necessary may compensate for values of $\lambda_{i,j}$ that we will not effectively observe. For this reason, in the case of Problem \ref{pro03}, $\lambda_{i,j}$ will not necessarily all be observable eigenvalues. For example, if we are able to hide $n$ modes, then we can obtain $\epsilon=0$, and none of  $\{\lambda_{i,j}\}_{i=1,\ldots,p,\;j=1,\ldots,\nu_i}$ will need to be part of the closed-loop eigenstructure.

Before proceeding with the solutions of the problems formulated in this section, in the next two sections we will discuss some geometric and combinatorial preliminaries that are needed for the main proofs of this paper.

\section{Geometric preliminaries}

We denote by $\gV^\star$ the largest output-nulling subspace of $\Sigma$, i.e., the largest subspace $\gV$ of $\gX$ for which a matrix $F\,{\in}\,\mathbb{R}^{m\,{\times}\,n}$ exists such that $(A+B\,F)\,\gV\subseteq \gV \subseteq \ker (C+D\,F)$. Any real matrix $F$ satisfying this inclusion is called a {\it friend \/} of $\gV$. %We denote by $\mathfrak{F}(\gV)$ the set of friends of $\gV$.
 The symbol $\gR^\star$ denotes the so-called {\em reachability subspace} on $\gV^\star$.
 The closed-loop spectrum can be partitioned as
$\sigma(A+B\,F)=\sigma(A+B\,F\,|\,\gV^\star)\uplus \sigma(A+B\,F\,|\,\gX/\gV^\star)$. Further, we have 
$\sigma(A+B\,F\,|\,\gV^\star)=\sigma(A+B\,F |\gR^\star) \uplus \sigma\,(A+B\,F | {\gV^\star}/{\gR^\star})$, where 
 $\sigma(A+B\,F |\gR^\star)$ is freely assignable
 with a suitable friend $F$ of $\gV^\star$, whereas  $\sigma\,(A+B\,F | {\gV^\star}/{\gR^\star})$ is fixed for every friend $F$ of $\gV^\star$ {and coincide with the invariant zero structure of $\Sigma$, \cite[Theorem 7.19]{Trentelman-SH-01}}.
 %Thus, $\gV^\star$ is inner stabilisable if and only if $\gZ \subset \complex_g$. 
 Finally, the symbol $\gV^\star_g$ denotes the largest stabilizability  subspace of $\Sigma$. \\
%Dually, we denote by $\gS_\star$ the smallest input-containing subspace of $\Sigma$, i.e., the smallest subspace $\gS$ of $\gX$ for which a matrix $L\,{\in}\,\mathbb{R}^{n\,{\times}\,p}$ exists such that $(A+L\,C)\,\gS\subseteq \gS$ and $\gS \supseteq  \ima (B+L\,D)$. The orthogonal complement of $\gS_\star$ is the largest output-nulling subspace of the dual system $\Sigma^\top=(A^\top,C^\top,B^\top,D^\top)$. A simple method to compute $\gR^\star$ is via the formula $\gR^\star=\gV^\star\cap\, \gS_{\star}$.
An important result for the computation of a basis for $\gR^\star$, which also offers a great deal of insight into the properties of this subspace, is based on the null-space of the Rosenbrock system matrix pencil, when $\lambda$ assumes arbitrary values that are distinct from the invariant zeros of the system.

Given the $h$ self-conjugate complex numbers $\gL=\{\lambda_1,\ldots,\lambda_h\}$ including exactly $s$ complex conjugate pairs, we say that $\gL$ is $s$-conformably indexed if $2\,s \le h$ and the first $2\,s$ values are complex, while the remaining are real, and for all odd $k \le 2\,s$ we have $\lambda_{k+1}=\bar{\lambda}_k$. 
The following important result holds, \cite[Proposition 4]{Moore-L-78}.

\begin{theorem}%{\sc [Moore-Laub, \cite[Prop. 4]{Moore-L-78}]}
\label{th1}
Let $r=\dim \gR^\star$. Let $\gL=\{\lambda_1,\lambda_2,\ldots,\lambda_r\}$ be an $s$-conformably indexed set of self-conjugate distinct complex numbers disjoint from the invariant zeros, i.e., $\gL \cap \gZ=\varnothing$. For all $k \in \{1,\ldots,r\}$, let us denote by $\bsmat X_k \\[1mm] Y_k \esmat$ a basis matrix for $\ker P_{\scriptscriptstyle \Sigma}(\lambda_k)$ partitioned conformably with $P _{\scriptscriptstyle \Sigma}(\lambda_k)$. Let this basis be chosen in such a way that $\bsmat X_{k+1} \\[1mm] Y_{k+1} \esmat= \bsmat \overline{X}_k \\[1mm] \overline{Y}_k \esmat$ when $k \le 2\,s$ is odd. Let
\bea
V_k \ns&\ns  \defi\ns&\ns \left\{ \begin{array}{ll} 
%\frac{1}{2} (X_k+X_{k+1})=
\mathfrak{Re}\{X_k\} & \text{if $k\le 2\,s$ is odd,} \\
%\frac{i}{2} (X_{k-1}-X_{k})=
\mathfrak{Im}\{X_k\} & \text{if $k\le 2\,s$ is even,} \\
{X}_{k} & \text{if $k>2\,s$,} \end{array} \right.  \qquad
W_k \defi\left\{ \begin{array}{ll} 
%\frac{1}{2} (Y_k+Y_{k+1})=
\mathfrak{Re}\{Y_k\} & \text{if $k\le 2\,s$ is odd,} \\
%\frac{i}{2} (Y_{k-1}-Y_{k})=
\mathfrak{Im}\{Y_k\} & \text{if $k\le 2\,s$ is even,} \\
{Y}_{k} & \text{if $k>2\,s$.}  \end{array} \right.\label{eqV}
\eea
Then, $\gR^\star= \ima [\begin{array}{cccc} V_1 & V_2 & \ldots & V_r\end{array}]$.
\end{theorem}

%As aforementioned, the assumption that the eigenvalues must be distinct can be dropped. In that case, however, it may not be sufficient to consider the null-spaces of the Rosenbrock matrix pencil to build basis vectors of $\gR^\star$. Indeed, in this case the vectors in the null-space of $P_{\scriptscriptstyle \Sigma}(\lambda)$ are the initial elements of a Jordan chain, which is built using the procedure outlined in \cite{Ntogramatzidis-14}.

The following corollary shows how the computation of a friend of $\gR^\star$ can be carried out. In particular, the values of $\lambda$ used to construct the basis of $\gR^\star$ will become, with such feedback $F$, eigenvalues of the closed-loop restricted to $\gR^\star$.

\begin{corollary}
\label{corcomplex}
Consider a basis for $\gR^\star$ as constructed in Theorem \ref{th1}. Let $\gR^\star=\ima [ \begin{array}{cccc} \! V_1 \!&\! V_2\!&\! \cdots \!&\!  V_r\! \end{array}]$. Let $\{v_1,\ldots,v_r\}$ be a set of columns extracted from the matrix $[\begin{array}{cccc}  \! V_1 \!&\! V_2\!&\! \cdots \!&\!  V_r\!\end{array}]$ to form a basis for $\gR^\star$, and let $\{w_1,\ldots,w_r\}$ denote the corresponding columns of $[\begin{array}{cccc}  \! W_1 \!&\! W_2\!&\! \cdots \!&\!  W_r\!\end{array}]$.
If $v_k$ is a column of $V_j$, let us denote by $\mu_k$ the eigenvalue $\lambda_j$. Let $\{v_1,\ldots,v_r\}$ be constructed in such a way that the multi-set $\{\mu_1,\ldots,\mu_r\}$ is self-conjugate. Then, the matrix
\bea
\label{ultima}
F=[\begin{array}{cccc} w_1 & w_2 & \ldots & w_r\end{array}]\,[\begin{array}{cccc} v_1 & v_2 & \ldots & v_r\end{array}]^\dagger
\eea
is a friend of $\gR^\star$, and $\sigma(A+B\,F\,|\,\gR^\star)=\{\mu_1,\ldots,\mu_r\}$.
\end{corollary}

Theorem \ref{th1} apparently requires the {\em a priori} knowledge of the dimension of $\gR^\star$ to determine a spanning set for $\gR^\star$. However, this knowledge is not necessary: in fact it is possible to compute a spanning set of $\gR^\star$ recursively, because when computing the sequence of subspaces $\{\ima V_k\}_{k \in \mathbb{N}}$, at each step $k$ the dimension of the subspace $\ima [\begin{array}{cccc}  \! V_1 \!&\! V_2\!&\! \cdots \!&\!  V_k\!\end{array}]$ increases with respect to the size of $\ima [\begin{array}{cccc}  \! V_1 \!&\! V_2\!&\! \cdots \!&\!  V_{k-1}\!\end{array}]$, until the dimension of $\gR^\star$ has been reached. In other words, considering the matrices $V_1,\ldots,V_r$ as obtained in Theorem \ref{th1}, 
for all $k \in \mathbb{N}$, 
we have $\rank [\begin{array}{cccc}  \! V_1 \!&\! V_2\!&\! \cdots \!&\!  V_{k-1}\!\end{array}]<r$ if and only if $\rank [\begin{array}{cccc}  \! V_1 \!&\! V_2\!&\! \cdots \!&\!  V_k\!\end{array}]> \rank [\begin{array}{cccc}  \! V_1 \!&\! V_2\!&\! \cdots \!&\!  V_{k-1}\!\end{array}]$.
 This follows from Theorem \ref{th1} and the Rosenbrock Theorem \cite{Rosenbrock-70}.

The second fundamental result is \cite[Proposition 5]{Moore-L-78}, and is about the construction of a basis matrix for $\gV^\star$ (resp. $\gV^\star_g$): the idea is essentially the same as the one for the construction of a basis for $\gR^\star$, but this time the invariant zeros (resp. minimum-phase invariant zeros) also have to be taken into account when choosing the $\lambda_i$ for which we compute the null-space of the Rosenbrock matrix.

\begin{theorem}%{\sc [Moore-Laub, \cite{Moore-L-78}]}\\
\label{th2}
Let $r=\dim \gR^\star$. Let $\gZ=\{z_1,z_2,\ldots,z_t\}$ be the $s$-conformably indexed set of self-conjugate invariant zeros (respectively, the minimum-phase invariant zeros). Let for all $k \in \{1,\ldots,t\}$ denote by $\bsmat X_k \\[1mm] Y_k \esmat$ a basis matrix for $\ker P _{\scriptscriptstyle \Sigma}(z_k)$ partitioned conformably with $P _{\scriptscriptstyle \Sigma}(z_k)$. Let this basis be chosen in such a way that $\bsmat X_{k+1} \\ Y_{k+1} \esmat= \bsmat \overline{X}_k \\ \overline{Y}_k \esmat$ when $k \le 2\,s$ is odd. Let $V_k$ and $W_k$ be constructed as in Theorem \ref{th1}.
Then, $\gV^\star= \gR^\star+\ima [\begin{array}{cccc}   \! V_1 \!&\! V_2\!&\! \cdots \!&\!  V_t\!\end{array}]$ (respectively, $\gV^\star_g= \gR^\star+\ima [\begin{array}{cccc}  \! V_1 \!&\! V_2\!&\! \cdots \!&\!  V_t\!\end{array}]$).
\end{theorem}

We finally recall that the following statements are equivalent:
\begin{itemize}
\item $\Sigma$ is right invertible;
\item $P_{\scriptscriptstyle \Sigma}(\l)$ is full row-rank for all but finitely many $\l\in \complex$;
\item the transfer function $G_{\scriptscriptstyle \Sigma}(\l)=C\,(\l\,I-A)^{-1}B+D$ is right invertible as a rational matrix. %;
%\item $\gV^\star+\gS_\star=\gX$.
\end{itemize}
%Dually, the following statements are equivalent:
%\begin{itemize}
%\item $\Sigma$ is left invertible;
%\item $P_{\scriptscriptstyle \Sigma}(\l)$ is full column-rank for all but finitely many $\l\in \complex$;
%\item the transfer function $G_{\scriptscriptstyle \Sigma}(\l)=C\,(\l\,I-A)^{-1}B+D$  is left invertible as a rational matrix;
%\item $\gR^\star=\gV^\star\cap \gS_\star=\{0\}$.
%\end{itemize}

\subsection{Preliminaries in combinatorial linear algebra and affine geometry}
Let $\field$ denote a field ($\real$ or $\complex$). 
We also recall that the dimension of a set $S$ of $\field^n$ is defined as the dimension of the smallest linear subspace that contains $S$ (i.e., the dimension of $\spanR_{\scriptscriptstyle \field} (S)$) or, equivalently, the maximum number of linearly independent vectors that it is possible to find in $S$. We recall that given two sets $S_1,S_2$ of the vector space ${\field}^n$, there holds
$\spanR_{\scriptscriptstyle \field} (S_1\cup S_2)=\spanR_{\scriptscriptstyle \field} S_1+\spanR_{\scriptscriptstyle \field} S_2$.

The following result is a cornerstone of Combinatorics, \cite[Theorem 3]{Rado-42}, and it will be the starting point of our investigation. %In what follows, we will use the notation $\field$ to denote a generic field.

\begin{theorem}{\sc [Rad{\'o}'s Theorem]}\\
\label{Rado}
Consider the sets $P_1,\ldots,P_q$ in the vector space $\field^n$. It is possible to find a linearly independent set $\{\xi_1,\ldots,\xi_q\}$ such that 
$\xi_1\in P_1$, $\xi_2\in P_2$, $\ldots$, $\xi_q\in P_q$ if and only if 
given $k$ numbers $\eta_1,\ldots,\eta_k\in \mathbb{N}$ such that 
$1\le \eta_1 < \eta_2 < \ldots < \eta_k \le q$ for all $k \in \{1,\ldots,q\}$,
the union $P_{\eta_1} \cup P_{\eta_2}\cup \ldots \cup P_{\eta_k}$ contains $k$ independent elements, i.e., if and only if
for any set $S\subseteq \{1,\ldots,q\}$ of cardinality $s=\operatorname{card}\,(S)$ there exist $s$ independent vectors $\zeta_1,\ldots,\zeta_s$ such that $\zeta_1,\ldots,\zeta_s \in \bigcup_{i\in S} P_i$.
\end{theorem}

The following corollary will be useful in the rest of the paper.

\begin{corollary}
\label{cor1}
Consider the sets $P_1,\ldots,P_q$ of vectors in the vector space $\field^n$. It is possible to find a set of linearly independent vectors $\{\xi_1,\ldots,\xi_q\}$ such that $\xi_1\in P_1$, $\xi_2\in P_2$, $\ldots$, $\xi_q\in P_q$ if and only if for any set $S\subseteq \{1,\ldots,q\}$  there holds
\[
\dim \Big(\sum_{i \in S} \spanR_{\scriptscriptstyle \field} P_i\Big)\ge \operatorname{card} S.
\]
\end{corollary}
\proof From Theorem \ref{Rado}, for any $S$ there exist $s=\operatorname{card}\,(S)$ vectors $\zeta_1,\ldots,\zeta_s \in \bigcup_{i\in S} P_i$ that are linearly independent if and only if 
$\dim \Big( \spanR_{\scriptscriptstyle \field} \bigl(\bigcup_{i \in S}P_i\bigr)\Big)\ge s$.
The statement follows noting that 
$\spanR_{\scriptscriptstyle \field}\bigl(\bigcup_{i \in S}P_i\bigr)=\sum_{i \in S} \spanR_{\scriptscriptstyle \field} P_i$.
%which concludes the proof.
\endproof
\ \\

The following corollary is a generalization of the latter.

\begin{corollary} 
\label{cor2}
Consider the sets $P_1,\ldots,P_q$ of vectors in $\field^n$ and $\nu_1,\ldots,\nu_q\in \mathbb{N}$. It is possible to find a set of linearly independent vectors $\{\xi_{1,1},\ldots,\xi_{1,\nu_1},\ldots,\xi_{q,1},\ldots,\xi_{q,\nu_q}\}$ such that 
$\xi_{i,1},\ldots,\xi_{i,\nu_i}\in P_i$ for $i \in \{1,\ldots,q\}$
 if and only if for any set $S\subseteq \{1,\ldots,q\}$ there holds 
\[
\dim \Big(\sum_{i \in S} \spanR_{\scriptscriptstyle \field} P_i\Big)\ge \sum_{i\in S} \nu_i.
\]
\end{corollary}
\proof The claim follows by considering the problem of finding a set of linearly independent vectors $\{\xi_{1,1},\ldots,\xi_{1,\nu_1},\ldots,\xi_{q,1},\ldots,\xi_{q,\nu_q}\}$ 
 such that $\xi_{1,1}\in P_{1,1},\ldots,\xi_{1,\nu_1}\in P_{1,\nu_1}$, $\ldots$, $\xi_{q,1}\in P_{q,1},\ldots,\xi_{q,\nu_q}\in P_{q,\nu_q}$, writing the condition of Corollary \ref{cor1} under the assumption
$P_{i,1}=P_{i,2}=\ldots=P_{i,\nu_i}=P_i$ for $i \in \{1,\ldots,q\}$.
\endproof

The following corollary %, although very simple, 
 highlights the fact that, when we are interested in selecting linearly independent vectors, what really matters is the span of the set $P_i$, rather than the set itself.

\begin{corollary}
\label{cor3}
Let $P_1,\ldots,P_q$ be sets of vectors in $\field^n$ and let $Q_1,\ldots,Q_q$ be sets of $\field^n$ such that $\spanR_{\scriptscriptstyle \field} P_i=\spanR_{\scriptscriptstyle \field} Q_i$ for $i \in \{1,\ldots,q\}$. It is possible to find a set of linearly independent vectors $\{\xi_1,\ldots,\xi_q\}$ such that $\xi_i\in P_i$ for all $i \in \{1,\ldots,q\}$ if and only if it is possible to find a set of linearly independent vectors $\{\zeta_1,\ldots,\zeta_q\}$ such that $\zeta_i\in Q_i$ for all $i \in \{1,\ldots,q\}$.
\end{corollary}

\proof Applying Corollary \ref{cor1} to $P_1,\ldots,P_q$ and $Q_1,\ldots,Q_q$, the two sets of conditions, for any set $S\subseteq\{1,\ldots,q\}$, are that
$\dim \Big(\sum_{i \in S} \spanR_{\scriptscriptstyle \field} P_i\Big)\ge \operatorname{card} S$ and $\dim \Big(\sum_{i \in S} \spanR_{\scriptscriptstyle \field} Q_i\Big)\ge \operatorname{card} S$. Since $\spanR_{\scriptscriptstyle \field}P_i=\spanR_{\scriptscriptstyle \field} Q_i$ for all $i \in \{1,\ldots,q\}$, the result readily follows.
\endproof

The previous result provides a guideline on the selection of the vectors in $P_1,\ldots,P_q$ by restricting the attention to the vectors of each $P_i$ that forms a basis for the subspace $\spanR_{\scriptscriptstyle \field} P_i$.

\begin{corollary}
\label{cor4} 
Let the vectors in $Q_i\subseteq P_i$ be basis vectors for  $\spanR_{\scriptscriptstyle \field} P_i$. 
It is possible to find a set of linearly independent vectors $\{\xi_1,\ldots,\xi_q\}$ such that $\xi_i\in P_i$ for all $i \in \{1,\ldots,q\}$ if and only if it is possible to find a set of linearly independent vectors $\{\zeta_1,\ldots,\zeta_q\}$ such that $\zeta_i\in Q_i$ for all $i \in \{1,\ldots,q\}$.
\end{corollary}

\proof The statement follows directly from Corollary \ref{cor3}, since a vector of $Q_i$ also belongs to $P_i$.
\endproof

%\tb{
We now consider another generalization of Rad{\'o}'s theorem, which considers the case where we want to extract at most $ k $ linearly independent vectors from $ q>k $ subspaces. The following theorem provides a solution to this problem.\footnote{This result is usually presented in the literature, see e.g. \cite[Theorem 1.3]{Mirsky-67} and \cite[Theorem 1.1]{Ore-55}, in terms of sets in an Euclidean space and expressed in terms of numbers of linearly independent vectors belonging to the unions of these sets. However, one can repeat verbatim the argument in the proof of 
Corollary \ref{cor1} to rewrite the same result in terms of the spans of these sets.}
     
\begin{theorem}
	\label{Ore}
	Consider the sets  ${P}_1,\ldots,{P}_q$ in the vector space $ \field^n $. It is possible to find a set of linearly independent vectors $\{\xi_1,\ldots,\xi_k\}$ such that $\xi_1\in {P}_{i_1}$, $\xi_2\in {P}_{i_2}$, $\ldots$, $\xi_q\in {P}_{i_k}$ for some 
	$ 1\leq i_1<i_2< \ldots<i_k\leq q $ if and only if  there holds
	\[
	\dim \Big(\sum_{i \in S} \spanR_{\scriptscriptstyle \field} P_i\Big)\ge \operatorname{card} S -(q-k)
	\]
	for all $ S\in\left\{\mathfrak{S}\in 2^{\{1,\ldots,q\}} |  \operatorname{card}\mathfrak{S}>q-k \right\}$. 
	\end{theorem}
%}

%\tb{
\begin{corollary}
\label{corollary_fabrizio}		
Consider the sets  $P_g, {P}_1,\ldots,{P}_q$ in the vector space $ \field^n $. Let $ h \geq n-q $ be the dimension of $ P_g $. There exists a linearly independent set of vectors $ \{\xi_{g_1},\ldots,\xi_{g_{n-k}},\xi_{i_1},\ldots,\xi_{i_k}\}$ such that $ \{\xi_{g_1},\ldots,\xi_{g_{n-k}}\} \in {P}_g $ and $ \xi_{i_j} \in {P}_{i_j} $ for some $ 1\leq i_1<i_2< \ldots<i_k\leq q $ and for some $ k \leq q $ if and only if
\be
\label{eq:fabrizio}
\dim \Big(\spanR_{\scriptscriptstyle \field}{P}_g + \sum_{i \in S} \spanR_{\scriptscriptstyle \field}{P}_i\Big)\ge n-q+\operatorname{card} S
\ee
holds for all $ S\in\left\{\mathfrak{S}\in 2^{\{1,\ldots,q\}} |  \operatorname{card}\mathfrak{S}>h-(n-q) \right\}$.
\end{corollary}
\proof
It is clear that if there exists the linearly independent set for some $ k $ such that $ n-k<h $ there always exists another linearly independent set for $ n- k=h $. Then, it is sufficient to prove the theorem when $ k=n-h $.\\
Let $ \field^n={\cal X}_1\oplus {\cal X}_2 $, where $ {\cal X}_1=\spanR_{\scriptscriptstyle \field}{P}_g $. In these coordinates a basis matrix of $ \spanR_{\scriptscriptstyle \field}{P}_g $ is given by $ \bsmat I_h \\[1mm] 0_{k\times h} \esmat $. Denote by $ \bsmat \Pi_{i,1} \\ \Pi_{i,2} \esmat $ a basis matrix for $ \spanR_{\scriptscriptstyle \field}{P}_i $, where $ \Pi_{i,1} $ and $ \Pi_{i,2} $ have $ h $ and $ k $ rows, respectively. We can find a linearly independent set $ \{\xi_{g_1},\ldots,\xi_{g_{h}},\xi_{i_1},\ldots,\xi_{i_k}\} $ such that $ \{\xi_{g_1},\ldots,\xi_{g_{h}}\} \in {P}_g $ and $ \xi_{i_j} \in {P}_{i_j} $ for some $ 1\leq i_1<i_2< \ldots<i_k\leq q $ and for $ k = q $ if and only if there exist $ \tilde \xi_{i_1} \in \im \Pi_{i_1,2},\ldots,\tilde \xi_{i_k} \in \im \Pi_{i_k,2} $ such that the set $ \{\tilde \xi_{i_1},\ldots,\tilde \xi_{i_k} \} $ is linearly independent. In view of Theorem \ref{Ore} this happens if and only if $ \dim \Big(\sum_{i \in S} \spanR_{\scriptscriptstyle \field}{P}_i\Big)\ge \operatorname{card} S -(q-k) $, for all $ S\in\left\{\mathfrak{S}\in 2^{\{1,\ldots,q\}} |  \operatorname{card}\mathfrak{S}>q-k \right\} $. Considering that $ k=n-h $ and that $ \spanR_{\scriptscriptstyle \field}{P}_g \cap \spanR_{\scriptscriptstyle \field}\{ \xi_{i_1},\ldots,\xi_{i_k} \} = \{0\}  $, \eqref{eq:fabrizio} is readily obtained.       
\endproof   
%}

We now specialize these results to the case where the field $\field$ is equal to $\complex$, see \cite[Lemma 1]{Kimura-75}.

\begin{theorem}{\sc [Kimura's Theorem]}
\label{kimura}
Consider the sets $P_1,\ldots,P_q\subseteq \complex^n$. It is possible to find a set of linearly independent vectors $\{\xi_1,\ldots,\xi_q\}$ such that 
$\xi_1\in P_1$, $\xi_2\in P_2$, $\ldots$, $\xi_q\in P_q$ if and only if for any set $S\subseteq \{1,\ldots,q\}$ of cardinality $s=\operatorname{card}\,(S)$ there holds
\[
\dim \left(\sum_{i \in S} \spanR_{\scriptscriptstyle \complex} P_i\right) \ge \operatorname{card}\,S.
\]
Moreover, for any pair $P_i, P_j$ that are linear subspaces such that $P_i=\overline{P}_j$ it is possible to guarantee that the further constraint $\xi_i=\overline{\xi}_j$ is satisfied.
\end{theorem}

The following result is an extension of Theorem \ref{kimura} to the case of affine sets of $\complex^n$.

\begin{theorem}
\label{kimura1}
Consider the sets $P_1,\ldots,P_q\subseteq \complex^n$. It is possible to find a set of linearly independent vectors $\{\xi_1,\ldots,\xi_q\}$ such that 
$\xi_1\in P_1$, $\xi_2\in P_2$, $\ldots$, $\xi_q\in P_q$ if and only if for any set $S\subseteq \{1,\ldots,q\}$ of cardinality $s=\operatorname{card}\,(S)$ there holds
\[
\dim \left(\sum_{i \in S} \spanR_{\scriptscriptstyle \complex} P_i\right) \ge \operatorname{card}\,S.
\]
Moreover:
\begin{itemize}
\item for every $P_i$ such that there exists a set of real vectors $Q_i \subseteq P_i$ for which $\spanR_{\scriptscriptstyle \complex} Q_i=\spanR_{\scriptscriptstyle \complex} P_i$, we can guarantee also that $\mathfrak{Im}\{\xi_i\}=0$; 
\item for any pair $P_i, P_j$ such that there exist two affine subspaces $Q_i \subseteq P_i$ and $Q_j \subseteq P_j$ and $Q_i=\overline{Q}_j$ and $\spanR_{\scriptscriptstyle \complex} Q_i=\spanR_{\scriptscriptstyle \complex} P_i$, we can guarantee also that the further constraint $\xi_i=\overline{\xi}_j$ is satisfied.
\end{itemize}
\end{theorem}
\proof
The proof of the first part follows from Corollary \ref{cor3}. Indeed, the existence of a vector in $P_i$ which is linearly independent from all the others is equivalent from the existence of a (real) vector from $Q_i$.

We prove the second point. Let us assume, with no loss of generality, that $P_1$ 
and $P_2$ are sets from which we want to extract two vectors $p_1\in P_1$ and $p_2\in P_2$ that are complex conjugate and linearly independent. Let $Q_1\subseteq P_1$ and $Q_2\subseteq P_2$ be such that $\spanR_{\scriptscriptstyle \complex} Q_1=\spanR_{\scriptscriptstyle \complex} P_1$ and $\spanR_{\scriptscriptstyle \complex} Q_2=\spanR_{\scriptscriptstyle \complex} P_2$, and $Q_i=\overline{Q}_j$; a linearly independent set $\{\xi_1,\ldots,\xi_q\}$  exists such that $\xi_i \in P_i$ for all $i \in \{1,\ldots,q\}$ if and only if a linearly independent set $\{\zeta_1,\ldots,\zeta_q\}$  exists such that $\zeta_1\in Q_1$, $\zeta_2\in Q_2$ and $\zeta_i \in P_i$ for all $i \in \{3,\ldots,q\}$. 

If $Q_i$ is an affine subspace, given two vectors $v_1,v_2\in Q_i$, for every $\l\in \complex$ their affine combination $\l\,v_1+(1-\l)\,v_2$ is in $Q_i$.

If $Q_1=\overline{Q}_2$ and we assume $\xi_1\neq \overline{\xi}_2$ such that $\xi_1\in Q_1$ and $\xi_2\in Q_2$, it is possible to construct the vectors
$w_1= \overline{\gamma}_1\,\xi_1+\gamma_2\,\overline{\xi}_2$ and $w_2 = {\gamma}_1\,\overline{\xi}_1+\overline{\gamma}_2\,{\xi}_2$,
where $\gamma_1,\gamma_2\in \complex$, such that by construction
\begin{enumerate}
\item $w_1=\overline{w}_2$;
\item since $\overline{\xi}_1\in \overline{P}_1=P_2$ and 
$\overline{\xi}_2\in \overline{P}_2=P_1$, then $w_1\in P_1$ and $w_2\in P_2$  if 
$\overline{\gamma}_1+\gamma_2=1$, i.e., if
\bea
\label{constgamma}
\mathfrak{Re}\{\gamma_1\}+\mathfrak{Re}\{\gamma_2\} = 1 \qquad \text{and} \qquad
\mathfrak{Im}\{\gamma_1\}=\mathfrak{Im}\{\gamma_2\}.
\eea
\end{enumerate}

We now have to prove that it is possible to find $\gamma_1,\gamma_2\in \complex$ such that (\ref{constgamma}) holds and such that the set of vectors $\{w_1,w_2,\xi_3,\ldots,\xi_q\}$ is linearly independent. 
The vectors $\overline{\xi}_1\in P_2$  and $\overline{\xi}_2\in P_1$ can be written as
\bea
\overline{\xi}_1 \ns&\ns = \ns&\ns \alpha_1\,\xi_1+\alpha_2\,\xi_2+\ldots+\alpha_n\,\xi_n+t_1\label{xi1} \\
\overline{\xi}_2 \ns&\ns = \ns&\ns \beta_1\,\xi_1+\beta_2\,\xi_2+\ldots+\beta_n\,\xi_n+t_2,\label{xi2} 
\eea
where $t_1$ and $t_2$ are suitable vectors orthogonal to $\sum_{i=1}^n \spanR_{\scriptscriptstyle \complex} \{\xi_i\}$.

%Our task is to replace $\xi_1$ and $\xi_2$ with two vectors $q_1\in \widehat{R}_i(\lambda)$ and $q_2\in \widehat{R}_i(\overline{\lambda})$ such that $\{q_1,q_2,\xi_3,\ldots,\xi_n\}$ is linearly independent and, moreover, $q_1=\overline{q}_2$.

To prove this point, we proceed similarly to what is done in \cite[Lemma 1]{Kimura-75} and we first show that, for all $\alpha_1,\alpha_2,\beta_1,\beta_2\in \complex$ there exist $\gamma_1,\gamma_2\in \complex$ such that 
\begin{enumerate}
\item $\gamma_1+\overline{\gamma}_2=1$;
\item the determinant of $\bsmat \overline{\gamma}_1+\gamma_2\,\beta_1 & \gamma_1\,\alpha_1 \\[1mm]
\gamma_2\,\beta_2 & \gamma_1\,\alpha_2+\overline{\gamma}_2\esmat$ is different from zero.
\end{enumerate}
Three cases must be considered:
\begin{enumerate}
\item if $\beta_1\neq 0$, choose $\gamma_1=0$ and $\gamma_2=1$, so that $\left|\begin{smallmatrix} \overline{\gamma}_1+\gamma_2\,\beta_1 & \gamma_1\,\alpha_1 \\[1mm]
\gamma_2\,\beta_2 & \gamma_1\,\alpha_2+\overline{\gamma}_2\end{smallmatrix} \right|=\left|\begin{smallmatrix} \beta_1 & 0 \\[1mm]
\beta_2 & 1\end{smallmatrix} \right|=\beta_1\neq 0$;
\item if $\beta_1=0$ and $\alpha_2\neq 0$, choose $\gamma_1=1$ and $\gamma_2=0$, so that 
$\left|\begin{smallmatrix} \overline{\gamma}_1+\gamma_2\,\beta_1 & \gamma_1\,\alpha_1 \\[1mm]
\gamma_2\,\beta_2 & \gamma_1\,\alpha_2+\overline{\gamma}_2\end{smallmatrix} \right|=\left|\begin{smallmatrix} 1 &\alpha_1 \\[1mm]
0 & \alpha_2\end{smallmatrix} \right|=\alpha_2\neq 0$;
\item if $\beta_1=0$ and $\alpha_2=0$, we have $\left|\begin{smallmatrix} \overline{\gamma}_1+\gamma_2\,\beta_1 & \gamma_1\,\alpha_1 \\[1mm]
\gamma_2\,\beta_2 & \gamma_1\,\alpha_2+\overline{\gamma}_2\end{smallmatrix} \right|=\left|\begin{smallmatrix}\overline{\gamma}_1 & \gamma_1\,\alpha_1 \\[1mm]
\gamma_2\,\beta_2 & \overline{\gamma}_2\end{smallmatrix} \right|= \overline{\gamma}_1\, \overline{\gamma}_2-\gamma_1\,\gamma_2\,\alpha_1\,\beta_2$.  Here we have to consider two subcases:
\begin{itemize}
\item if $\alpha_1\,\beta_2\neq 1$, by choosing $\gamma_1=\gamma_2=\frac{1}{2}$ the determinant becomes $\frac{1}{4}\,(1-\alpha_1\,\beta_2)\neq 0$;
\item if $\alpha_1\,\beta_2= 1$, by choosing $\gamma_1=1+\i$ and $\gamma_2=\i$ the determinant becomes $\overline{\gamma}_1\, \overline{\gamma}_2-\gamma_1\,\gamma_2=2\,\i$.
\end{itemize}
\end{enumerate}

We now show that $\{w_1,w_2,\xi_3,\ldots,\xi_q\}$ is linearly independent. Suppose by contradiction that there exist $\kappa_1,\ldots,\kappa_q\in \complex$ not all zero such that $\kappa_1\,w_1+\kappa_2\,w_2+\kappa_3\,\xi_3+\ldots+\kappa_q\,\xi_q=0$.
Using the definition of $w_1$ and $w_2$, and (\ref{xi1}-\ref{xi2}), we find
\beann
&&\hspace{-1cm} (\kappa_1\,\overline{\gamma_1}+\kappa_1\,\gamma_2\,\beta_1+\kappa_2\,\gamma_1\,\alpha_1)\,\xi_1+
(\kappa_1\,{\gamma_2}\,\beta_2+\kappa_2\,\gamma_1\,\alpha_2+\kappa_2\,\overline{\gamma}_2)\,\xi_2+\\
&& +\sum_{i=3}^q (\kappa_1\,\gamma_2\,\beta_i+\kappa_2\,\gamma_1\,\alpha_i+\kappa_i)\,\xi_i+\kappa_1\,\gamma_2\,p_2+\kappa_2\,\gamma_1\,p_1=0.
\eeann
Since $\{\xi_1,\ldots,\xi_q\}$ is a linearly independent set, all the coefficients in the latter are zero. Thus, in particular
$\bsmat \overline{\gamma}_1+\gamma_2\,\beta_1 & \gamma_1\,\alpha_1 \\[1mm]
\gamma_2\,\beta_2 & \gamma_1\,\alpha_2+\overline{\gamma}_2\esmat\bsmat \kappa_1\\[1mm] \kappa_2 \esmat=0$. Since the determinant of the matrix in the left hand side is non-zero, the only solution is $\kappa_1=\kappa_2=0$, and therefore also $\kappa_3=\ldots=\kappa_q=0$. This is a contradiction.
\endproof

\section{Solution of Problem 1}
\label{sec:problem1_real}
For the sake of simplicity, in this part of the paper we consider only the case where the eigenvalues to be assigned and the stable invariant zeros are real. The change that occurs where invariant zeros or eigenvalues to be assigned are in complex conjugate pairs will be discussed in Section \ref{sec:complex}. 
Nevertheless, whenever possible, the definitions of the subspaces used in the sequel will be given in the general case where the indeterminate is complex to avoid repetitions.
%\textcolor{red}{(si ma ora si usano i complessi)}
Let for all $\lambda \in \complex$
\[
\gR(\lambda) \defi \left\{v\in \complex^n\,\Big|\,\exists \,w\in \complex^m: \,\bmat{cc} A-\lambda\,I & B \\ C & D\emat\bmat{c} v \\ w\emat =0 \right\}.
\]
%Clearly, when $\lambda\in \real$, then $\gR(\lambda)=\left\{v\in \real^n\,|\,\exists \,w\in \real^m: \,P_{\scriptscriptstyle \Sigma}(\l) \bsmat v \\[1mm] w\esmat =0 \right\}$ is a subspace of $\real^n$; when $\l\in \complex$,  then $\gR(\lambda)$ is a subspace of $\complex^n$. 
Notice that 
$\bsmat A-\lambda\,I && B \\[1mm] C && D\esmat\bsmat v \\[1mm] w\esmat =0$ if and only if $\bsmat A-\overline{\lambda}\,I && B \\[1mm] C && D\esmat\bsmat \overline{v} \\[1mm] \overline{w}\esmat =0$, from which we find
$\gR(\l)=\overline{\gR(\overline{\l})}$.
Let us also define
\[
\gR_{i}(\lambda) \defi \left\{v\in \complex^n\,\Big|\,\exists \,w\in \complex^m: \,\bmat{cc} A-\lambda\,I & B \\ C_{(i)} & D_{(i)} \emat\bmat{c} v \\ w\emat=0 \right\},
\]
where $C_{(i)}$ and $D_{(i)}$ are matrices obtained from $C$ and $D$ by removing their $i$-th rows.

As a direct consequence of Theorem \ref{th1}, denoting by $r$ the dimension of $\gR^\star$, if $\lambda_1,\ldots,\lambda_r$ are real, distinct and different from the invariant zeros, there holds
\bea
\label{par}
\gR^\star=\gR(\lambda_1)+\gR(\lambda_2)+\ldots+\gR(\lambda_r),
\eea
and if the minimum-phase invariant zeros $z_1,\ldots,z_t$ are real, we also have
\[
\gV^\star_g=\gR(\lambda_1)+\ldots+\gR(\lambda_r)+\gR(z_1)+\ldots+\gR(z_t).
\]

It is worth observing that (\ref{par}) cannot be used to exhaustively parameterize the vectors of $\gR^\star$; in other words, given an arbitrary $v\in \gR^\star$, there might not exist $\lambda \in \complex$ such that $v \in \gR(\lambda)$. Consider for example the quadruple
$A=\bsmat -1 &&    0   && 0\\[1mm]
 3&&    2  &&   0 \\[1mm]
     -1   &&  2  &&   0 \esmat$, $B=\bsmat   0 \\[1mm]
 1\\[1mm]
 0\esmat$, $C=\bsmat 1    && 0  &&   0 \esmat$ and $D=0$.
     In this case $\gR^\star=\ima \bsmat 0 && 0 \\[1mm] 1 && 0 \\[1mm] 0 && 1 \esmat$,
but if we take $v=[\begin{array}{cccc}0 & 1 & 0 \end{array} ]^\top\in \gR^\star$, we cannot find $w\in \complex^m$ and $\lambda \in \complex$ such that $\bsmat A-\l\,I && B \\[1mm] C && D \esmat \bsmat v \\[1mm] w \esmat=0$; to see this, we can re-write this as the linear equation
$\bsmat -v && B \\[1mm] 0 && D \esmat\bsmat \l \\[1mm] w \esmat =-\bsmat A \\[1mm] C \esmat \,v$
and notice that $\bsmat A \\[1mm] C \esmat \,v=\bsmat 0\\[1mm] 2\\[1mm] 2\\[1mm] 0\esmat \notin \ima \bsmat -v && B \\[1mm] 0 && D \esmat=\spanR\left\{\bsmat 0\\[1mm] 1\\[1mm] 0\\[1mm] 0 \esmat\right\}$.

It is easily seen that for all $\lambda \in \real$, the sets $\gR(\lambda)$ and $\gR_i(\lambda)$ are subspaces of $\gX$ for all $i \in \{1,\ldots,p\}$.
For all $\lambda \in \complex$, we also define the set
\[
\widehat{R}_i(\lambda) \defi \left\{v\in \complex^n
\,\Big|\,\exists \,w\in \complex^m, \;\; \exists \,\delta \in \real\setminus \{0\}: \,\bmat{cc} A-\lambda\,I & B \\[0mm] C & D\emat\bmat{c} v \\[0mm] w\emat=\bmat{c} 0\\[0mm] \delta\,e_i \emat \right\}.
\]
Clearly, in general, the set $\widehat{R}_i(\lambda)$ is not a subspace of $\complex^n$.
%; rather, it can be seen geometrically as the set-theoretic difference between a linear subspace in $\complex^n$ and a linear subspace strictly contained in it. 
For reasons that will be clearer later, it is worth also to define the sets
\[
\widehat{W}_i(\lambda)=\left\{v\in \complex^n
\,\Big|\,\exists \,w\in \complex^m \;:\; \bmat{cc} A-\lambda\,I & B \\[0mm] C & D\emat\bmat{c} v \\[0mm] w\emat=\bmat{c} 0\\[0mm] e_i \emat \right\}.
\]
Notice that for every vector $v \in \widehat{R}_i(\l)$, there exist a vector $v^\prime \in \widehat{W}_i(\l)$ which is parallel to $v$ (so that, in particular, their spans coincide). Indeed, $\bsmat A-\lambda\,I & B \\[1mm] C & D\esmat\bsmat v \\[1mm]  w\esmat=\bsmat 0\\[1mm]  \delta\,e_i \esmat$, for $\delta \neq 0$, can be rewritten as $\bsmat A-\lambda\,I & B \\[1mm] C & D\esmat\bsmat \frac{1}{\delta} \, v \\[1mm]   \frac{1}{\delta} \,w\esmat=\bsmat 0\\[1mm]  e_i \esmat$. Notice also that $\widehat{W}_i(\lambda)$ is an affine set in $\complex^n$. Indeed, given $v_1,v_2 \in \widehat{W}_i(\lambda)$,
there exist $w_1,w_2\in \complex^m$ such that $\bsmat A-\l\,I && B \\[1mm] C && D \esmat \bsmat v_i \\[1mm] w_i \esmat=\bsmat 0 \\[1mm] e_i \esmat$ for $i\in \{1,2\}$; for any $\alpha \in \complex$, the vector $\alpha\,v_1+(1-\alpha)\,v_2$ is also in $\widehat{W}_i(\lambda)$. This can be seen by taking $v=\alpha\,v_1+(1-\alpha)\,v_2$ and $w=\alpha\,w_1+(1-\alpha)\,w_2$.

%Thus, we can also write
%\[
%\widehat{R}_i(\lambda)=\left\{v\in \complex^n \setminus \{0\} 
%\,\Big|\,\exists \,w\in \complex^m \;:\; \bmat{cc} A-\lambda\,I & B \\[0mm] C & D\emat\bmat{c} v \\[0mm] w\emat=\bmat{c} 0\\[0mm] e_i \emat \right\}.
%\]
%In this way, we can also view $\widehat{R}_i(\lambda)$ as an affine subspace parallel to $\ker P_{\scriptscriptstyle \Sigma}(\lambda)$. 

Finally, we notice that there holds
$\widehat{R}_i(\lambda)=\overline{\widehat{R}_i(\overline{\lambda})}$ and $\widehat{W}_i(\lambda)=\overline{\widehat{W}_i(\overline{\lambda})}$.

Given $\lambda \in \real$, the set $\widehat{R}_i(\lambda)$  contains the non-zero initial states for which a state-feedback control $u=F\,x$ exists for which every output except the $i$-th is zero, while the $i$-th is given by a single exponential. Indeed, consider $v\in \widehat{R}_i(\lambda)$, and let $w\in \gU$ and $\delta \in \real \setminus \{0\}$. Since $v \neq 0$, choosing $F\,v=w$ gives
\beann
(A+B\,F)\,v \ns&\ns = \ns&\ns \l\,v \\
(C+D\,F)\,v \ns&\ns = \ns&\ns \delta\,e_i.
\eeann
Let $x_0=v$. Then, recalling that $e^{(A+B\,F)\,t}\,v=e^{\l\,t}\,v$, we find that from
\beann
\dot{x}(t)\ns&\ns = \ns&\ns (A+B\,F) \,x(t), \qquad x_0=v, \\
{y}(t)\ns&\ns = \ns&\ns (C+D\,F) \,x(t), 
\eeann
we get
\beann
y(t)=(C+D\,F)\,e^{(A+B\,F)\,t}\,v=(C+D\,F)\,e^{\l\,t}\,v=\delta\,e_i\,e^{\l\,t}=\bmat{c} 0 \\ \vdots \\ \delta\,e^{\l\,t} \\ \vdots \\ 0 \emat. 
\eeann

The next result shows that the only invariant zeros that it is necessary to compute are those of the original system, because the invariant zeros of all the systems obtained by removing outputs are a subset of the former.

\begin{lemma}
For all $i\in \{1,\ldots,p\}$ 
\[
\gZ(A,B,C,D)\supseteq \gZ(A,B,C_{(i)},D_{(i)}).
\]
\end{lemma}
\proof 
The statement follows directly from the right invertibility of the system.
\endproof

%For the sake of simplicity, we assume that the observable eigenvalues are all real. The case where these are complex can be addressed as shown in the sequel for the unobservable eigenvalues. \footnote{We have no choice on the invariant zeros}

\begin{lemma}
\label{lem1}
Let $\mu \in \complex \setminus \gZ$. 
For all $j \in \{ 1,\cdots , p\}$, there holds
\[
\gR_j(\mu) \supset \gR(\mu).
\]
\end{lemma}
\proof
First, notice that $\gR_j(\mu) \supseteq \gR(\mu)$.
The row $[\begin{array}{cc} C_{j} & D_{j}\end{array}]$ is linearly independent from every row of $\left[ \begin{smallmatrix} A-\mu\,I_n & B \\[1mm] C_{(j)} & D_{(j)} \end{smallmatrix} \right]$. This implies that $\dim \gR(\mu) < \dim \gR_j (\mu)$.
\endproof

\begin{lemma}
\label{lem2}
Let $\mu \in  \complex  \setminus \gZ$. 
For all $j \in \{ 1,\cdots , p\}$, there holds
\bea
\label{emeq}
\gR_j(\mu)\supseteq \widehat{R}_j(\mu)\supseteq \gR_j(\mu)\setminus \gR(\mu).
\eea
\end{lemma}
\proof The fact that $\gR_j(\mu)\supseteq \widehat{R}_j(\mu)$ follows directly from the definition. We now show that $\widehat{R}_j(\mu)\supseteq \gR_j(\mu)\setminus \gR(\mu)$. Let $v\in \gR_j(\mu)\setminus \gR(\mu)$; since $v\in \gR_j(\mu)$ there exists $w\in \gU$ such that 
\[
\bmat{cc} A-\mu\,I & B \\ C_{(j)} & D_{(j)} \emat \bmat{c} v \\ w \emat=0.
\]
On the other hand, since $v\notin  \gR(\mu)$, there are no $\omega\in \gU$ for which 
\[
\bmat{cc} A-\mu\,I & B \\ C & D \emat \bmat{c} v \\ \omega \emat=0.
\]
Thus, there must hold
\[
\bmat{cc} A-\mu\,I & B \\ C & D \emat \bmat{c} v \\ w \emat=\bmat{c} 0 \\ \delta\,e_j \emat
\]
for some $\delta \neq 0$. Thus, $v\in \widehat{R}_j(\mu)$.
\endproof

\begin{lemma}
\label{lemincl}
Let $\mu \in {\complex}  \setminus \gZ$. 
For all $j \in \{1,\cdots , p\}$%, there holds
\[
\spanR_{\scriptscriptstyle {\complex}} \widehat{R}_j(\mu)=\gR_j(\mu).
\]
\end{lemma}
\proof Taking the span on each term of (\ref{emeq}) we get
\[
\spanR_{\scriptscriptstyle \complex}\gR_j(\mu) \supseteq \spanR_{\scriptscriptstyle \complex}\widehat{R}_j(\mu)\supseteq \spanR_{\scriptscriptstyle \complex} \bigl(\gR_j(\mu)\setminus \gR(\mu)\bigr).
\]
We have $\spanR_{\scriptscriptstyle \complex} \gR_j(\mu)=\gR_j(\mu)$, because $\gR_j(\mu)$ is a linear subspace. Recall that given two linear subspaces $\gA$ and $\gB$ such that $\gA\subset \gB$ (which means that $\gA\subseteq \gB$ and $\dim \gA < \dim \gB$) we have $\spanR_{\scriptscriptstyle \complex}(\gB\setminus \gA)=\spanR_{\scriptscriptstyle \complex} \gB=\gB$. Thus, $\spanR_{\scriptscriptstyle \complex}\bigl(\gR_j(\mu)\setminus \gR(\mu)\bigr)=\spanR_{\scriptscriptstyle \complex}\bigl(\gR_j(\mu)\bigr)=\gR_j(\mu)$. Thus, we find
\[
\gR_j(\mu) \supseteq \spanR_{\scriptscriptstyle \complex}\widehat{R}_j(\mu)\supseteq \gR_j(\mu),
\]
which immediately implies that $\spanR_{\scriptscriptstyle \complex}\bigl(\widehat{R}_j(\mu)\bigr)=\gR_j(\mu)$.
\endproof

\subsection{Problem \ref{pro01}A}
We begin by giving a necessary and sufficient condition for the solvability condition of Problem \ref{pro01}A that, even if not expressed in terms of the problem data, will turn out to be constructive for the calculation of the feedback matrix whenever the problem admits solutions. 

\begin{lemma}
\label{lemmac1}
Let $\nu_0=n-\nu_1-\ldots-\nu_p$.
Problem \ref{pro01}A is solvable if and only if there exist
\beann
&& v_{0,k} \in \gR(\lambda_{0,k})\qquad \,k\in\{1,\ldots,\nu_0\}\\
&& v_{i,j} \in \widehat{R}_i(\lambda_{i,j})\qquad  \,i\in\{1,\ldots,p\}, \;\; \,j\in \{1,\ldots,\nu_i\}
\eeann
such that the set $\bigl\{v_{0,1},\ldots,v_{0,\nu_0},v_{1,1},\ldots,v_{1,\nu_1},\ldots,v_{p,1},\ldots,v_{p,\nu_p}\bigr\}$ is linearly independent.
\end{lemma}
%
%\textcolor{blue}{Se un eig sta su due uscite $i_1,i_2$, devo definire
%\[
%\widehat{R}_{i_1,i_2}=\Big\{... \bmat{c} 0 \\ e_{i_1}+\alpha\,e_{i_2}\emat\Big\}
%\]
%e da questo estrarre un solo vettore.
%}

\proof
Let us prove sufficiency. Since $v_{0,k} \in \gR(\lambda_{0,k})$ for $k \in \{1,\ldots,\nu_0\}$, there exist $w_{0,1},\ldots,w_{0,\nu_0}\in \real^m$ such that 
$\bsmat A-\l_{0,k}\,I && B \\[1mm] C && D \esmat \bsmat v_{0,k} \\[1mm] w_{0,k} \esmat=0$ for $k \in \{1,\ldots,\nu_0\}$.
Moreover, from $v_{i,j} \in \widehat{R}_i(\lambda_{i,j})$ for $i\in \{1,\ldots,p\}$ and $j\in \{1,\ldots,\nu_i\}$, we have
$\bsmat A-\l_{i,j}\,I && B \\[1mm] C && D \esmat \bsmat v_{i,j} \\[1mm] w_{i,j} \esmat=\bsmat 0 \\[1mm] \delta_{i,j}\, e_i\esmat$, for $i\in \{1,\ldots,p\}$ and $j \in \{1,\ldots,\nu_i\}$
for some $\delta_{i,j}\neq 0$. Since $\nu_0=n-\nu_1-\ldots-\nu_p$ and $\{v_{0,1},\ldots,v_{0,\nu_0},v_{1,1},\ldots,v_{1,\nu_1},\ldots,$ $v_{p,1},\ldots,v_{p,\nu_p}\}$ are linearly independent, then $\{v_{0,1},\ldots,v_{0,\nu_0},v_{1,1},\ldots,v_{1,\nu_1},\ldots,v_{p,1},\ldots,v_{p,\nu_p}\}$ is a basis for $\gX$, and we can define
\bea
F \ns&\ns = \ns&\ns [\begin{array}{cccccccccccccccc}
w_{1,1} & \!\! \ldots \!\! & w_{1,\nu_1} & | & \ldots  & | &  w_{p,1} & \!\! \ldots \!\!  & w_{p,\nu_p}  & | &  w_{0,1} & \!\! \ldots \!\!  & w_{0,\nu_0}
\end{array}] \nonumber \\
\ns&\ns  \ns&\ns \hspace{2cm}
\times[\begin{array}{ccccccccccccccc}
v_{1,1} & \!\! \ldots \!\!  & v_{1,\nu_1} & | &  \ldots & | &  v_{p,1} &\!\! \ldots \!\!  & v_{p,\nu_p}  & | &  v_{0,1} & \!\! \ldots \!\!  & v_{0,\nu_0}
\end{array}]^{-1},\label{compF}
\eea
from which we find
\beann
&&\hspace{-.2cm} (A+B\,F)[\begin{array}{ccccccccccccccc}
\!\! v_{1,1} & \!\! \ldots \!\!  & v_{1,\nu_1} \!\! &\!\!  | \!\! &\!\!   \ldots\!\!  &\!\!  | \!\! &  \!\! v_{p,1}\!\!  &\!\! \ldots \!\!  & \!\! v_{p,\nu_p} \!\!  &\!\!  |\!\!  & \!\!  v_{0,1} \!\! & \!\! \ldots \!\!  & \!\! v_{0,\nu_0}\!\! 
\end{array}]  \\
&& \hspace{-.2cm}
=\diag\{\lambda_{1,1},\ldots,\lambda_{1,\nu_1},\ldots,\lambda_{p,1},\ldots,\lambda_{p,\nu_p},\lambda_{0,1},\ldots,\lambda_{0,\nu_0}\}\\
&& \qquad \times
%\bmat{cccccccccccc}
% \! \lambda_{1,1 \! } &\\[-2mm]
% \! & \! \!\!\ddots\!\! \\[-2mm]
% \! & \! & \!  \lambda_{1,\nu_1} \!  \\[-2mm]
% \! & \! & \! & \!  \!\!\ddots\!\! \!  \\[-2mm]
% \! & \! & \! & \! & \! \lambda_{p,1} \!  \\[-2mm]
% \! & \! & \! & \! & \! & \!  \!\!\ddots\!\! \! \\[-2mm]
% \! & \! & \! & \! & \! & \! &  \! \lambda_{p,\nu_p} \!  \\[-2mm]
% \! & \! & \! & \! & \! & \! & \! &  \! \lambda_{0,1}  \! \\[-2mm]
% \! & \! & \! & \! & \! & \! & \! & \! & \!  \!\!\ddots\!\! \!\\[-2mm]
% \! & \! & \! & \! & \! & \! & \! & \! & \! &  \! \lambda_{0,\nu_0}
%  \emat
 [\begin{array}{ccccccccccccccc}
\!\! v_{1,1} & \!\! \ldots \!\!  & v_{1,\nu_1} \!\! &\!\!  | \!\! &\!\!   \ldots\!\!  &\!\!  | \!\! &  \!\! v_{p,1}\!\!  &\!\! \ldots \!\!  & \!\! v_{p,\nu_p} \!\!  &\!\!  |\!\!  & \!\!  v_{0,1} \!\! & \!\! \ldots \!\!  & \!\! v_{0,\nu_0}\!\! 
\end{array}] \\
&& \hspace{-.2cm}(C+D\,F)[\begin{array}{ccccccccccccccc}
\!\! v_{1,1} & \!\! \ldots \!\!  & v_{1,\nu_1} \!\! &\!\!  | \!\! &\!\!   \ldots\!\!  &\!\!  | \!\! &  \!\! v_{p,1}\!\!  &\!\! \ldots \!\!  & \!\! v_{p,\nu_p} \!\!  &\!\!  |\!\!  & \!\!  v_{0,1} \!\!& \!\! \ldots \!\!  & \!\! v_{0,\nu_0}\!\! 
\end{array}] \\
&& \hspace{2cm}
= \left\{\begin{array}{ll} \delta_{i,j} \, e_i & i\in \{1,\ldots,p\}\\
0 & i \in \{p+1,\ldots,n\}\end{array} \right.
\eeann
The first says that
\[
e^{(A+B\,F)\,t}\,v_{i,j}=\exp(\lambda_{i,j}\,t)\,v_{i,j}
\]
 for $i\in \{0,\ldots,p\}$ and $j\in\{1,\ldots,\nu_i\}$.  Let $\xi_0=\xi(0)$ be the initial error state, and define 
 \[
 \alpha \defi [\begin{array}{ccccccccccccccc}
\! v_{1,1} & \! \ldots \!  & v_{1,\nu_1} \! &\!  | \! &\!   \ldots\!  &\!  | \! &  \! v_{p,1}\!  &\! \ldots \!  & \! v_{p,\nu_p} \!  &\!  |\!  & \!  v_{0,1} \!& \! \ldots \!  & \! v_{0,\nu_0}\! 
\end{array}]^{-1}\,\xi_0.
\]
 The second yields
 \beann
 \epsilon(t)   \ns&\ns = \ns&\ns (C+D\,F)\,e^{(A+B\,F)\,t}\,\xi_0\\
 \ns&\ns = \ns&\ns (C+D\,F)\,\sum_{i=0}^p \sum_{j=1}^{\nu_i} \exp(\lambda_{i,j}\,t)\,v_{i,j}\,\alpha_{i,j}\\
  \ns&\ns = \ns&\ns \sum_{i=1}^p \sum_{j=1}^{\nu_i} \delta_{i,j}\,e_i\,\exp(\lambda_{i,j}\,t)\,\alpha_{i,j}=
  \bmat{c} 
  \delta_{1,1}\,\alpha_{1,1}\,e^{\lambda_{1,1}\,t}+\ldots+\delta_{1,\nu_1}\,\alpha_{1,\nu_1}\,e^{\lambda_{1,\nu_1}\,t}\\
  \vdots \\
   \delta_{p,1}\,\alpha_{p,1}\,e^{\lambda_{p,1}\,t}+\ldots+\delta_{p,\nu_p}\,\alpha_{p,\nu_p}\,e^{\lambda_{p,\nu_p}\,t}\emat
 \eeann
as required. We now establish necessity. Suppose we have 
 \bea
 \label{confronto1}
 \epsilon(t)   \ns&\ns = \ns&\ns
  \bmat{c} 
  \gamma_{1,1}\,e^{\lambda_{1,1}\,t}+\ldots+\gamma_{1,\nu_1}\,e^{\lambda_{1,\nu_1}\,t}\\
  \vdots \\
   \gamma_{p,1}\,e^{\lambda_{p,1}\,t}+\ldots+\gamma_{p,\nu_p}\,e^{\lambda_{p,\nu_p}\,t}\emat,
 \eea
{where $\gamma_{i,j}$ can be made arbitrary by suitably choosing $\xi_0$}.
It follows that $n-\nu_1-\ldots-\nu_p=\nu_0$ closed-loop modes are unobservable. We denote these modes by $\lambda_{0,1},\ldots,\lambda_{0,\nu_0}$. Since $\lambda_{0,k}$ is not observable and is not an invariant zero, the corresponding closed-loop eigenvector $v_{0,k}$ is in $\gR^\star$ for all $k \in \{1,\ldots,\nu_0\}$. %Then, there must hold  $v_{0,j}\in\gR^\star(\lambda_{0,j})$, because defining $w_{0,j}=F\,v_{0,j}$ for all $j \in \{1,\ldots,\nu_0\}$ there must hold $\bsmat A-\l_{0,j}\,I && B \\[1mm] C && D \esmat\bsmat v_{0,j} \\[1mm] w_{0,j} \esmat=0$. 
Similarly, denoting by $v_{i,j}$ the closed-loop eigenvector associated with $\l_{i,j}$  and defining $\alpha=[\begin{array}{ccccccc}v_{0,1} &\ldots &v_{0,\nu_0} & \ldots & v_{p,1} & \ldots & v_{p,\nu_p} \end{array}]^{-1}\,\xi_0$, and $w_{i,j}=F\,v_{i,j}$ for all $i \in \{0,\ldots,p\}$ and $j \in \{1,\ldots,\nu_i\}$, we find
\bea
\epsilon(t) \ns&\ns = \ns&\ns (C+D\,F)\,e^{(A+B\,F)\,t}\,\xi_0 \nonumber  \\
\ns&\ns = \ns&\ns(C+D\,F)\,e^{(A+B\,F)\,t}\,[\begin{array}{ccccccccccc}v_{0,1} & \ldots & v_{0,\nu_0} & v_{1,1} & \ldots & v_{1,\nu_1} & \ldots & v_{p,1} & \ldots & v_{p,\nu_p} \end{array}]\,\alpha \nonumber  \\
%\ns&\ns = \ns&\ns(C+D\,F)\,[\begin{array}{ccccccccccc}v_{0,1} & \ldots & v_{0,\nu_0} & v_{1,1} & \ldots & v_{1,\nu_1} & \ldots & v_{p,1} & \ldots & v_{p,\nu_p} \end{array}]\nonumber  \\
%\ns&\ns\ns&\ns \;\; \qquad \times
%\diag\{\lambda_{0,1},\ldots,\lambda_{0,\nu_0},\lambda_{1,1},\ldots,\lambda_{1,\nu_1},\ldots,\lambda_{p,1},\ldots,\lambda_{p,\nu_p}\}\alpha  \nonumber  \\
\ns&\ns = \ns&\ns \sum_{i=0}^p \sum_{j=0}^{\nu_i}
(C+D\,F)\,v_{i,j}\,e^{\l_{i,j}\,t}\alpha_{i,j},\label{confronto2}
\eea
where $\alpha=[\begin{array}{cccccccccc} \alpha_{0,1} &\ldots&\alpha_{0,\nu_0}&\alpha_{1,1}&\ldots&\alpha_{1,\nu_1}&\ldots&\alpha_{p,1}&\ldots&\alpha_{p,\nu_p} \end{array}]^\top$ is partitioned conformably. Comparing (\ref{confronto1}) with (\ref{confronto2}), there must hold:
\begin{itemize}
\item
 $(C+D\,F)\,v_{0,k}=0$ for all $k \in \{1,\ldots,\nu_0\}$. It follows that $\bsmat A-\l_{0,k}\,I && B \\[1mm] C && D \esmat \bsmat v_{0,k} \\[1mm] w_{0,k} \esmat=0$ for all  $k \in \{1,\ldots,\nu_0\}$, which proves that $v_{0,k}\in \gR(\lambda_{0,k})$ or all  $k \in \{1,\ldots,\nu_0\}$;
 \item defining $(C+D\,F)\,v_{i,j}\,\alpha_{i,j}=\sum_{\ell =1}^p e_{\ell}\,\phi_{\ell}$ for some coefficients $\phi_1,\ldots,\phi_p$; we must have $\phi_{\ell}=0$ for all $\ell\neq i$, or else the coefficients $\gamma_{i,j}$ would not be arbitrary. Thus, $\phi_i=\gamma_{i,j}$ so that 
 $(C+D\,F)\,v_{i,j}\,\alpha_{i,j}=e_i\,\gamma_{i,j}$.
 Hence, for an initial state such that $\alpha_{i,j} \neq 0$ we have 
  %$\alpha_{i,j} \neq 0$ for all $i \in \{1,\ldots,p\}$ and $j \in \{1,\ldots,\nu_i\}$, 
  $[\begin{array}{cc} C &D \end{array} ]\bsmat v_{i,j} \\[1mm] w_{i,j} \esmat=e_i\,\frac{\gamma_{i,j}}{\alpha_{i,j}}$, which together with $[\begin{array}{cc} A-\l_{i,j}\,I &B \end{array} ]\bsmat v_{i,j} \\[1mm] w_{i,j} \esmat=0$ implies $v_{i,j}\in \widehat{R}_{i}(\l_{i,j})$.
 \end{itemize}
\endproof

%
%
%Now we apply Rado.
%
%We define
%\beann
%\gX_{i,j}=\left\{\begin{array}{ll}
%\gR^\star_i(\lambda_{i,j})\setminus \gR^\star(\lambda_{i,j}) & i=1,\ldots,p\;j=1,\ldots,\nu_i \\
%\gR^\star(\lambda_{0,j})& \forall \,j=1,\ldots,\nu_0\end{array} \right.
%\eeann

%
%\begin{lemma}
%Let $\nu_0=n-\nu_1-\ldots-\nu_p$.
%Problem \ref{pro01}A\footnote{We recall that we are considering the case where the visible eigenvalues are not invariant zeros, but the invisible can be.} is solvable if and only if
%\beann
%\left. \begin{array}{lll} \dim \gX_{i,j} \ge 1 & i=0,\ldots,p\; j=1,\ldots,\nu_i \\
%\dim (\gX_{i,j}\cup \gX_{k,h}) \ge 2 & i,k=0,\ldots,p\; j=1,\ldots,\nu_i,\;\; h=1,\ldots,\nu_k \; (i,j)\neq (k,h)  \\
%\vdots
%\end{array} \right.
%\eeann
%or, concisely, if and only if
%\[
%\dim \bigcup_{(i,j)\in P} \gX_{i,j} \ge \operatorname{card}\,P
%\]
%for all $P$ in the power set $2^I$ where $I=\{(0,1),\ldots,(0,\nu_0),\ldots,(p,1),\ldots,(p,\nu_p)\}$.
%\end{lemma}
%
%la dimostrazione e' essenzialmente Rado.
%Non c'e' niente da fare.
For conciseness of notation, we define
$\gR_0(\lambda_{0,k}) \defi \gR(\lambda_{0,k})$ for $k \in \{1,\ldots,\nu_0\}$.

The following result provides a necessary and sufficient condition for the solvability of Problem \ref{pro01}A written in terms of the parameters of the problem.

\begin{theorem}
\label{pro1Asol}
Let $\nu_0=n-\nu_1-\ldots-\nu_p$.
Problem \ref{pro01}A is solvable if and only if
\bea
\label{cond1A}
\dim \left(\sum_{(i,j)\in P} \gR_{i}(\lambda_{i,j})\right) \ge \operatorname{card}\,P
\eea
for all $P$ in the power set $2^I$ where $I=\{(0,1),\ldots,(0,\nu_0),\ldots,(p,1),\ldots,(p,\nu_p)\}$.
\end{theorem}

\proof
From Lemma \ref{lemincl}, 
 there holds $\spanR_{\scriptscriptstyle \real}\widehat{R}_j(\mu)=\gR_j(\mu)$
for $\mu \in \real \setminus \gZ$ and $j \in \{ 1,\cdots , p\}$. Applying Corollary \ref{cor1}, the statement immediately follows.
\endproof

With Theorem \ref{pro1Asol}, we have obtained a set of necessary and sufficient conditions for the solution to Problem \ref{pro01}A. These conditions are very easy to check, because they are expressed in terms of the subspaces $\gR_i(\lambda_{i,j})$. In order to construct the feedback matrix, we can use the result in Lemma \ref{lemincl}. Indeed, if the conditions of Theorem \ref{pro1Asol} are satisfied, almost all choices of vectors 
$v_{0,k} \in \gR(\lambda_{0,k})$ for $k\in\{1,\ldots,\nu_0\}$ and 
$v_{i,j} \in \widehat{R}_i(\lambda_{i,j})$ for $i\in\{1,\ldots,p\}$ and $j\in \{1,\ldots,\nu_i\}$
will be such that the set $\bigl\{v_{0,1},\ldots,v_{0,\nu_0},\ldots,v_{p,1},\ldots,v_{p,\nu_p}\bigr\}$ is linearly independent, as the following result establishes.

%If we denote by $w_{0,j}$ and $w_{i,j}$ the corresponding lower part, so that
%\[
%\bsmat A-\lambda_{0,j}\,I && B \\[1mm] C && D\esmat\bsmat v_{0,j} \\[1mm] w_{0,j} \esmat=0 \quad \forall j\in\{1,\ldots,\nu_0\} 
%\]
%and 
%\[
%\bsmat A-\lambda_{i,j}\,I && B \\[1mm] C && D\esmat\bsmat v_{i,j} \\[1mm] w_{i,j} \esmat=\bsmat 0 \\[1mm] e_i \esmat\quad \forall \,i\in\{1,\ldots,p\}\;\; \forall \,j\in \{1,\ldots,\nu_i\}
%\]
%then we can compute $F$ as in (\ref{compF}).

\begin{theorem}
\label{theF1A}
Let the conditions of Theorem \ref{pro1Asol} hold true. 
Let  $V_{0,k}$ and $W_{0,k}$ be such that $\bsmat V_{0,k} \\[1mm] W_{0,k}\esmat$ be a basis matrix for $\ker \bsmat A-\l_{0,k}\,I && B \\[1mm] C && D \esmat$ for all $k \in \{1,\ldots,\nu_0\}$ and let $\bsmat V_{i,j} \\[1mm] W_{i,j}\esmat$ be a basis matrix for $\ker \bsmat A-\l_{i,j}\,I && B \\[1mm] C_{(i)} && D_{(i)} \esmat$ for all $i \in \{1,\ldots,p\}$ and $j \in \{1,\ldots,\nu_i\}$. Let $k_{i,j}$ be parameter vectors of suitable size, for $i \in \{0,\ldots,p\}$ and $j \in \{1,\ldots,\nu_i\}$, such that we can define
\beann
V_{k_{i,j}} \ns&\ns = \ns&\ns  [\begin{array}{ccc|c|ccc} 
V_{0,1} k_{0,1} & \ldots & V_{0,\nu_0} k_{0,\nu_0} & \ldots &
V_{p,1} k_{p,1} & \ldots & V_{p,\nu_p} k_{p,\nu_p} \end{array}],\\
W_{k_{i,j}} \ns&\ns = \ns&\ns  [\begin{array}{ccc|c|ccc} 
W_{0,1} k_{0,1} & \ldots & W_{0,\nu_0} k_{0,\nu_0} & \ldots &
W_{p,1} k_{p,1} & \ldots & W_{p,\nu_p} k_{p,\nu_p} \end{array}].
\eeann
Then:
\begin{enumerate}
\item the rank of $V_{k_{i,j}}$ is equal to $n$ for almost all parameters $k_{i,j}$, $i \in \{0,\ldots,p\}$ and $j \in \{1,\ldots,\nu_i\}$;
\item For almost all $k_{i,j}$, $i \in \{0,\ldots,p\}$ and $j \in \{1,\ldots,\nu_i\}$ such that $\rank V_{k_{i,j}}=n$, the feedback matrix
\bea
\label{paramF}
F=W_{k_{i,j}}\,V_{k_{i,j}}^{-1},
\eea
solves Problem \ref{pro01}A.
\end{enumerate}
\end{theorem}
\proof  First, we observe that there exist $k_{i,j}$, $i \in \{0,\ldots,p\}$ and $j \in \{1,\ldots,\nu_i\}$ such that the matrix
\beann
\Omega \ns&\ns = \ns&\ns [\begin{array}{ccc|c|ccc} 
V_{0,1} k_{0,1} & \ldots & V_{0,\nu_0} k_{0,\nu_0} & \ldots &
V_{p,1} k_{p,1} & \ldots & V_{p,\nu_p} k_{p,\nu_p} \end{array}]\\
 \ns&\ns = \ns&\ns
 [\begin{array}{ccc|c|ccc} 
V_{0,1}& \ldots & V_{0,\nu_0} & \ldots &
V_{p,1} & \ldots & V_{p,\nu_p} \end{array}]\diag\{ k_{0,1},\ldots,k_{0,\nu_0},\ldots,k_{p,1}, \ldots,k_{p,\nu_p}\}
\eeann
has rank equal to $n$. The rank of matrix $ [\begin{array}{ccc|c|ccc} 
V_{0,1}& \ldots & V_{0,\nu_0} & \ldots &
V_{p,1} & \ldots & V_{p,\nu_p} \end{array}]$ is equal to $n$ from the condition (\ref{cond1A}). Thus, $\Omega$ loses rank only for values of $k_{i,j}$, $i \in \{0,\ldots,p\}$ and $j \in \{1,\ldots,\nu_i\}$ for which a set of linear equations are satisfied. This proves the first point.
%To prove the second point, we need to show that the parameterization given by (\ref{paramF}) in $k_{i,j}$, $i \in \{0,\ldots,p\}$ and $j \in \{1,\ldots,\nu_i\}$  is exhaustive. 
We now prove the second point. We first show that every feedback matrix $F$ that solves Problem \ref{pro01}A can be written as in (\ref{paramF}).
Let $F$ be a feedback matrix that solves Problem \ref{pro01}A. Let $v_{0,k}\in \gR(\l_{0,k})$ for $k \in \{1,\ldots,\nu_0\}$. Then, $F$ satisfies $\bsmat A+B\,F \\[1mm] C+D\,F\esmat v_{0,k}=\bsmat v_{0,k} \\[1mm] 0 \esmat\,\l_{0,k}$ for $k \in \{1,\ldots,\nu_0\}$. Likewise, let $v_{i,j}\in \widehat{R}_i(\l_{i,j})$ for $i \in \{1,\ldots,p\}$ and $j \in \{1,\ldots,\nu_i\}$. 
Then, since $\widehat{R}_i(\l_{i,j}) \subseteq \gR_i(\l_{i,j})$, matrix $F$ satisfies $\bsmat A+B\,F \\[1mm] C_{(i)}+D_{(i)}\,F\esmat v_{i,j}=\bsmat v_{i,j} \\[1mm] 0 \esmat\,\l_{i,j}$ for $i \in \{1,\ldots,p\}$ and $j \in \{1,\ldots,\nu_i\}$. 
Defining $w_{i,j}=F\,v_{i,j}$ for $i \in \{0,\ldots,p\}$ and $j \in \{1,\ldots,\nu_i\}$, we obtain $\bsmat A-\l_{0,k}\,I && B \\[1mm] C && D \esmat \bsmat v_{0,k} \\[1mm] w_{0,k} \esmat=0$ for $k \in \{1,\ldots,\nu_0\}$ and 
 $\bsmat A-\l_{i,j}\,I && B \\[1mm] C_{(i)} && D_{(i)} \esmat \bsmat v_{i,j} \\[1mm] w_{i,j} \esmat=0$ for $i \in \{1,\ldots,p\}$ and $j \in \{1,\ldots,\nu_i\}$. Thus, $F$ satisfies $F\,[\begin{array}{ccc|c|ccc} 
v_{0,1}& \ldots & v_{0,\nu_0} & \ldots &
v_{p,1} & \ldots & v_{p,\nu_p} \end{array}]=[\begin{array}{ccc|c|ccc} 
w_{0,1}& \ldots & w_{0,\nu_0} & \ldots &
w_{p,1} & \ldots & w_{p,\nu_p} \end{array}]$. Moreover, $[\begin{array}{ccc|c|ccc} 
w_{0,1}& \ldots & w_{0,\nu_0} & \ldots &
w_{p,1} & \ldots & w_{p,\nu_p} \end{array}]$ and $[\begin{array}{ccc|c|ccc} 
v_{0,1}& \ldots & v_{0,\nu_0} & \ldots &
v_{p,1} & \ldots & v_{p,\nu_p} \end{array}]$ can be written as $W_{k_{i,j}}\,\diag\{ k_{0,1},\ldots,k_{0,\nu_0},\ldots,k_{p,1}, \ldots,k_{p,\nu_p}\}$ and $V_{k_{i,j}}\,\diag\{ k_{0,1},\ldots,k_{0,\nu_0},\ldots,$ $k_{p,1}, \ldots,k_{p,\nu_p}\}$ for a suitable choice of the parameters $k_{i,j}$, $i \in \{0,\ldots,p\}$ and $j \in \{1,\ldots,\nu_i\}$. We conclude the proof by noting that the set of parameters $k_{i,j}$ for which $v_{i,j} \in \gR_i(\l_{i,j})\setminus \widehat{R}_i(\l_{i,j})$ has zero Lebesgue measure.
\endproof

%First, recall that the sum of two linear subspaces can be defined as the smallest subspace containing their union. 
%Recall that the dimension of the union of two sets is equal to the dimension of the sum of the subspaces spanned by those sets from Lemma \ref{nac}. To conclude the proof it is enough to show that since 
%\[
%\dim \gR^\star_i(\lambda_{i,j})> \dim \gR^\star(\lambda_{i,j})
%\]
%(from the right invertibility, since $\lambda_{i,j}$ is not an invariant zero, when removing the $i$-th row of $[C\,\;\;D]$ the Rosenbrock matrix thus obtained is still full row-rank) then 
%\[
%\spanR \bigl(\gR^\star(\lambda_{i,j})\setminus \gR^\star(\lambda_{i,j})\bigr)=
%\gR^\star(\lambda_{i,j}).
%\]

\subsection{Problem \ref{pro01}B}
%Recall that we are excluding the case of complex invariant zeros.
We now consider the problem in which the unobservable closed-loop eigenvalues are not assigned but stable. To this end, we
 define the set
\beann
E_g \ns&\ns \defi \ns&\ns \bigcup_{\lambda \in \real_g} \gR(\lambda) \\
\ns&\ns = \ns&\ns
 \left\{v\in \real^n\,\Big|\,\exists \,\lambda\in \real_g,\;\exists \,w\in \real^m\,:
\bmat{cc} A-\lambda\,I & B \\ C & D \emat \bmat{c} v \\ w \emat=0\right\}.
\eeann

\begin{lemma}
\label{lemagg}
There holds $\spanR_{\scriptscriptstyle \real} E_g=\gV^\star_g$.
\end{lemma}

\proof This is a simple consequence of Theorem \ref{th2}. Indeed, a spanning set for the subspace $\gV^\star_g$ is therein constructed exactly by taking vectors of $E_g$.
\endproof

\begin{lemma}
\label{lemmac2}
Let $\nu_0=n-\nu_1-\ldots-\nu_p$.
Problem \ref{pro01}B is solvable if and only if there exist
\beann
&& v_{0,k} \in E_g \qquad \forall \,k\in\{1,\ldots,\nu_0\}\\
&& v_{i,j} \in \widehat{R}_i(\lambda_{i,j}) \qquad \forall \,i\in\{1,\ldots,p\} \;\; \forall \,j\in \{1,\ldots,\nu_i\}
\eeann
 such that $\{v_{0,1},\ldots,v_{0,\nu_0}, \ldots,v_{p,1},\ldots,v_{p,\nu_p}\}$ is linearly independent.
\end{lemma}

\proof 
The proof can be carried along the same lines of that of Lemma \ref{lemmac1}. 
Indeed, in the part of sufficiency the only difference is that $v_{0,k} \in E_g$ implies that there exist $\l_{0,k}\in \real_g$ and $w_{0,k}\in \real^m$ such that $\bsmat A-\l_{0,k}\,I && B \\[1mm] C && D \esmat\bsmat v_{0,k} \\[1mm] w_{0,k} \esmat=0$ for all $k \in \{1,\ldots,\nu_0\}$. Necessity is the same as in the proof of Lemma \ref{lemmac1}, since $\gR(\l_{0,k})\subset E_g$ for all $\l_{0,k}$.
\endproof

%
%\begin{lemma}
%Let $\nu_0=n-\nu_1-\ldots-\nu_p$.
%Problem \ref{pro01}B is solvable if and only if
%\beann
%\left. \begin{array}{lll} 
%\dim \gE_g \ge \nu_0 & \\
%\dim \gX_{i,j} \ge 1 & i=1,\ldots,p\; j=1,\ldots,\nu_i \\
%\dim \gE_g\cup \gX_{i,j} \ge 1+\nu_0 & i=1,\ldots,p\; j=1,\ldots,\nu_i \\
%\dim (\gX_{i,j}\cup \gX_{k,h}) \ge 2 & i,k=1,\ldots,p\; j=1,\ldots,\nu_i,\;\; h=1,\ldots,\nu_k \; (i,j)\neq (k,h)  \\
%\dim (\gE_g\cup \gX_{i,j}\cup \gX_{k,h}) \ge 2+\nu_0 & i,k=1,\ldots,p\; j=1,\ldots,\nu_i,\;\; h=1,\ldots,\nu_k \; (i,j)\neq (k,h)  \\
%\vdots
%\end{array} \right.
%\eeann
%or, concisely, if and only if
%\[
%\dim (\bigcup_{(i,j)\in P} \gX_{i,j})\cup \gE_g \ge \operatorname{card}\,P+\nu_0
%\]
%and
%\[
%\dim \bigcup_{(i,j)\in P} \gX_{i,j} \ge \operatorname{card}\,P
%\]
%for all $P$ in the power set $2^I$ where $I=\{(1,1),\ldots,(1,\nu_1),\ldots,(p,1),\ldots,(p,\nu_p)\}$.
%\end{lemma}
%%
%\proof We apply Rado's theorem by considering the set $\gE^\star_g$ $\nu_0$ times. Given $(i,j)\in P$, we have the conditions
%\[
%\dim (\bigcup_{(i,j)\in P} \gX_{i,j})\cup \gE_g \ge \operatorname{card}\,P+\ell 
%\]
%with $\ell \in \{1,\ldots,\nu_0\}$, so that the conditions for  $\ell \in \{1,\ldots,\nu_0-1\}$ are redundant and can eliminated.
%\endproof

\begin{theorem}
\label{pro1Bsol}
Let $\nu_0=n-\nu_1-\ldots-\nu_p$.
Problem \ref{pro01}B is solvable if and only if
\bea
\label{cond1B}
\dim \left(\gV^\star_g+ \sum_{(i,j)\in P} \gR_{i}(\lambda_{i,j})\right) \ge \operatorname{card}\,P+\nu_0
\eea
and
\[
\dim \left(\sum_{(i,j)\in P} \gR_{i}(\lambda_{i,j})\right) \ge \operatorname{card}\,P
\]
for all $P$ in the power set $2^I$ where $I=\{(1,1),\ldots,(1,\nu_1),\ldots,(p,1),\ldots,(p,\nu_p)\}$.\end{theorem}

\proof
%Using the same argument we obtain the N\&S conditions 
%\[
%\dim \bigl((\gE_g\cup \bigcup_{(i,j)\in P} \gR^\star_{i}(\lambda_{i,j})\setminus\gR^\star(\lambda_{i,j}) \bigr) \ge \operatorname{card}\,P+\nu_0
%\]
%and
%\[
%\dim \sum_{(i,j)\in P} \gR^\star_{i}(\lambda_{i,j}) \ge \operatorname{card}\,P
%\]
Since $\gR_{i}(\lambda_{i,j})=\spanR_{\scriptscriptstyle \real}\widehat{R}_{i}(\lambda_{i,j})$ 
and $\gV^\star_g$ is the smallest subspace containing $E_g$ because $\gV^\star_g=\spanR_{\scriptscriptstyle \real}E_g$ in view of Lemma \ref{lemagg}, then we can apply Corollary \ref{cor2} and the statement follows.\endproof

The next result shows how to construct the feedback matrix that solves Problem \ref{pro01}B.

\begin{theorem}
\label{thepar1B}
Let the conditions of Theorem \ref{pro1Bsol} hold true. 
Let  $V_{0,k}$ and $W_{0,k}$ be such that $\bsmat V_{0,k} \\[1mm] W_{0,k}\esmat$ is a basis matrix for $\ker \bsmat A-\l_{0,k}\,I && B \\[1mm] C && D \esmat$ for some $\l_{0,k} \in \real_g$, possibly including minimum-phase invariant zeros, for all $k \in \{1,\ldots,\nu_0\}$ and let $\bsmat V_{i,j} \\[1mm] W_{i,j}\esmat$ be a basis matrix for $\ker \bsmat A-\l_{i,j}\,I && B \\[1mm] C_{(i)} && D_{(i)} \esmat$ for all $i \in \{1,\ldots,p\}$ and $j \in \{1,\ldots,\nu_i\}$. Let $k_{i,j}$ be parameter vectors of suitable size, for $i \in \{0,\ldots,p\}$ and $j \in \{1,\ldots,\nu_i\}$, such that we can define
\beann
V_{k_{i,j}} \ns&\ns = \ns&\ns  [\begin{array}{ccc|c|ccc} 
V_{0,1} k_{0,1} & \ldots & V_{0,\nu_0} k_{0,\nu_0} & \ldots &
V_{p,1} k_{p,1} & \ldots & V_{p,\nu_p} k_{p,\nu_p} \end{array}],\\
W_{k_{i,j}} \ns&\ns = \ns&\ns  [\begin{array}{ccc|c|ccc} 
W_{0,1} k_{0,1} & \ldots & W_{0,\nu_0} k_{0,\nu_0} & \ldots &
W_{p,1} k_{p,1} & \ldots & W_{p,\nu_p} k_{p,\nu_p} \end{array}].
\eeann
Then:
\begin{enumerate}
\item the rank of $V_{k_{i,j}}$ is equal to $n$ for almost all parameters $k_{i,j}$, $i \in \{0,\ldots,p\}$ and $j \in \{1,\ldots,\nu_i\}$;
\item For almost all $k_{i,j}$, $i \in \{0,\ldots,p\}$ and $j \in \{1,\ldots,\nu_i\}$ such that $\rank V_{k_{i,j}}=n$, the feedback matrix
\bea
\label{paramF1}
F=W_{k_{i,j}}\,V_{k_{i,j}}^{-1},
\eea
solves Problem \ref{pro01}B.
%the set of all feedback matrices solving Problem \ref{pro01}B is given by
%\bea
%\label{paramF1}
%F=W_{k_{i,j}}\,V_{k_{i,j}}^{-1},
%\eea
%where $k_{i,j}$, $i \in \{0,\ldots,p\}$ and $j \in \{1,\ldots,\nu_i\}$ are such that $\rank V_{k_{i,j}}=n$.
\end{enumerate}
\end{theorem}
The proof can be carried out along the same lines of the proof of Theorem \ref{theF1A}, and it is therefore omitted.

\subsection{Problem \ref{pro01}C}
We finally consider the case where none of the closed-loop eigenvalues is assigned.
Define
\[
E_i \defi \left\{v\in \real^n \,\Big|\,\,\exists\, \lambda \in \real_g\setminus \gZ, \;\;\exists \,w \in \real^m, \; \exists \,\delta\in \real\setminus \{0\}: \bmat{cc} A-\lambda\,I & B \\ C & D \emat\bmat{c} v \\ w \emat=\bmat{c} 0 \\ \delta\, e_i \emat\right\}.
\]

\begin{lemma}
\label{lemagg2}
For all $i \in \{1,\ldots,p\}$ there holds 
\[
 \spanR_{\scriptscriptstyle \real}E_i=\gR^\star_i.
 \]
 \end{lemma}
 \proof
By definition we have
$E_i=\bigcup_{\lambda\in \real_g\setminus \gZ} \widehat{R}_i(\lambda)$. % \setminus  \gR^\star(\lambda)
 Thus,
 \beann 
 \spanR_{\scriptscriptstyle \real} E_i \ns&\ns = \ns&\ns \spanR_{\scriptscriptstyle \real}\Big(\bigcup_{\lambda\in \real_g\setminus \gZ} \widehat{R}_i(\lambda)\Big) \\
  \ns&\ns = \ns&\ns \sum_{\lambda\in \real_g\setminus \gZ} \spanR_{\scriptscriptstyle \real} \widehat{R}_i(\lambda) \\
    \ns&\ns = \ns&\ns \sum_{\lambda\in \real_g\setminus \gZ} \gR_i(\lambda)=\gR^\star_i,
    \eeann
    where the last equality follows from Theorem \ref{th1}. 
    \endproof
 
% and also 
% \[
% \spanR\{\gE_i\}=\spanR \{\gR^\star_i(\lambda) \setminus  \gR^\star(\lambda)\,|\,\lambda \in \real_g\setminus \gZ_g\}
% \]

\begin{lemma}
Let $\nu_0=n-\nu_1-\ldots-\nu_p$.
Problem \ref{pro01}C is solvable if and only if there exist
\beann
&& v_{0,k} \in E_g \qquad \forall \,k\in\{1,\ldots,\nu_0\}\\
&& v_{i,j} \in E_i \qquad \forall \,i\in\{1,\ldots,p\} \;\; \forall \,j\in\{1,\ldots,\nu_i\}
\eeann
 such that $\{v_{0,1},\ldots,v_{0,\nu_0}, \ldots,v_{p,1},\ldots,v_{p,\nu_p}\}$  is linearly independent.
\end{lemma}

\proof
This result follows by adapting the proof of Lemma \ref{lemmac2} considering this time that the sets $E_i$ represent the sets from which the closed-loop eigenvalues can be effectively extracted using an arbitrary closed-loop eigenvalue. Thus, 
in the sufficiency the only difference is that $v_{i,j} \in E_i$ implies that there exist $\l_{i,j}\in \real_g$ and $w_{i,j}\in \real^m$ such that $\bsmat A-\l_{i,j}\,I && B \\[1mm] C && D \esmat\bsmat v_{i,j} \\[1mm] w_{i,j} \esmat=\bsmat 0 \\[1mm] \delta_{i,j}\,e_i\esmat$ for all $i \in \{1,\ldots,p\}$ and $j\in \{1,\ldots,\nu_i\}$. Necessity is the same as in the proof of Lemma \ref{lemmac1}, since $\gR(\l_{0,k})\subseteq E_g$ for all $\l_{0,k}$ and 
$\widehat{R}_i(\l_{i,j})\subseteq \gR_i(\l_{i,j}) \subseteq E_i$.
\endproof

%
%\begin{lemma}
%Let $\nu_0=n-\nu_1-\ldots-\nu_p$.
%Problem \ref{pro01}C is solvable if and only if
%\beann
%\left. \begin{array}{lll} 
%\dim \gE_g \ge \nu_0 & \\
%\dim \gE_i \ge \nu_i & i=1,\ldots,p\\
%\dim \gE_g\cup \gE_i \ge \nu_i+\nu_0 & i=1,\ldots,p \\
%\dim (\gL_i\cup \gE_j) \ge \nu_i+\nu_j & i,j=1,\ldots,p\; i\neq j  \\
%\dim (\gE_g\cup \gE_i\cup \gE_j) \ge \nu_0+\nu_i+\nu_j & i,j=1,\ldots,p\; \;\; i\neq j  \\
%\vdots
%\end{array} \right.
%\eeann
%or, concisely, if and only if
%\[
%\dim (\bigcup_{i\in P} \gE_{i})\cup \gE_g \ge \nu_0+\sum_{i\in P} \nu_i
%\]
%and
%\[
%\dim(\bigcup_{i\in P} \gE_{i}) \ge \sum_{i\in P} \nu_i
%\]
%for all $P$ in the power set $2^I$ where $I=\{1,2,\ldots,p\}$.
%\end{lemma}
%%
%\proof
%Choosing $\nu_i$ vectors from $\gL_i$ is equivalent, in Rado, to consider $\gE_i$ $\nu_i$ times. For example, we have
%\beann
%\dim \left(\gE_g\cup \bigl(\bigcup_{i\in P,\;i\neq j} \gE_{i}\bigr) \cup \gL_j\right) \ns&\ns \ge \ns&\ns  \nu_0+\sum_{i\in P,\;i\neq j} \nu_i+1 \\
%\dim \left(\gE_g\cup \bigl(\bigcup_{i\in P,\;i\neq j} \gE_{i}\bigr) \cup \gE_j\right)  \ns&\ns \ge  \ns&\ns \nu_0+\sum_{i\in P,\;i\neq j} \nu_i+2\\
% \ns&\ns  \vdots \ns&\ns \\
%\dim \left(\gE_g\cup \bigl(\bigcup_{i\in P,\;i\neq j} \gE_{i}\bigr) \cup \gE_j\right)  \ns&\ns \ge \ns&\ns  \nu_0+\sum_{i\in P,\;i\neq j} \nu_i+\nu_j 
%\eeann
%so that the first $\nu_j-1$ conditions are redundant. In a similar way we can prove
%\[
%\dim(\bigcup_{i\in P} \gE_{i}) \ge \sum_{i\in P} \nu_i
%\]
%\endproof

\begin{theorem}
\label{th:1C}
Let $\nu_0=n-\nu_1-\ldots-\nu_p$.
Problem \ref{pro01}C is solvable if and only if
\[
\dim \left(\gV^\star_g+ \sum_{i\in P} \gR^\star_i\right) \ge \sum_{i\in P} \nu_i+\nu_0
\]
and
\[
\dim \left(\sum_{i\in P} \gR^\star_i\right) \ge \sum_{i\in P} \nu_i
\]
for all $P$ in the power set $2^I$ where $I=\{1,2,\ldots,p\}$.\end{theorem}

\proof We recall that $\gV^\star_g=\spanR_{\scriptscriptstyle \real} E_g$ (see Lemma \ref{lemagg}) and that $\spanR_{\scriptscriptstyle \real} E_i=\gR^\star_i$ (see Lemma \ref{lemagg2}), then 
%There holds $\gL_i\subseteq \gR^\star_i \setminus \gR^\star$ and
$\dim E_i=\dim \gR^\star_i$.
Therefore, we can apply Corollary \ref{cor2} and we obtain the result. 
%If $\gR^\star_i \setminus \gR^\star$ is not zero, then $\gR^\star_i$ is the smallest subspace that contains $\gL_i$. In fact, given $\lambda$, $\gL_i(\lambda)$ is spanned by a single vector. If $\lambda$ is unspecified, $\gL_i(\lambda)$ is a one dimensional vector space parameterized in $\lambda$, and for each $\lambda$ the vector $v$ obtained by solving
%$\bmat{cc} A-\lambda\,I & B \\ C& D \emat \bmat{c} v \\ w \emat= \bmat{c} 0 \\ e_i \emat
%$ is a vector of $\gR^\star_i \setminus \gR^\star$. Choosing at most $n$ different values of $\lambda$ and computing the corresponding vectors $v$ in this fashion, we construct a basis for $\gR^\star_i$ by Moore and Laub. Thus, $\spanR \gL_i=\gR^\star_i$.
\endproof

The construction of the feedback matrix $F$ that solves Problem \ref{pro01}C is carried out exactly in the same way as described in Theorem \ref{thepar1B}.

\section{Solution of Problem 2}
\label{sec:problem2_real}
Let us now consider Problem \ref{pro02}. We recall that this problem requires that in output $i$ we can observe exactly $\nu_i$ modes, which, differently from Problem \ref{pro01}, this time can be chosen also among the minimum-phase invariant zeros. 
For all $i\in \{1,\ldots,p\}$, let us define the sets
\[
L_i \defi \left\{v\in \real^n\,\Big|\,\exists \,\lambda\in \real_g,\;\,\exists \,w\in \real^m,\;\;\exists \;\delta \in \real\setminus \{0\}:
\bmat{cc} A-\lambda\,I & B \\ C& D \emat \bmat{c} v \\ w \emat= \bmat{c} 0 \\ \delta\,e_i \emat\right\}.
\]
What distinguishes the set $L_i$ from the set $E_i$ defined earlier is the fact that in $L_i$ now we are allowing $\lambda$ to be a minimum-phase invariant zero.
We also define
\[
T_i \defi \left\{v\in \real^n\,\Big|\,\exists \,\lambda\in \real_g,\;\;\exists \,w\in \real^m\,:
\bmat{cc} A-\lambda\,I & B \\ C_{(i)} & D_{(i)}  \emat \bmat{c} v \\ w \emat= 0\right\}.
\]
We allow again $\lambda$ to be an minimum-phase invariant zero.
Notice that the span of $T_i$ is the supremal stabilizability subspace of the system $(A,B,C_{(i)},D_{(i)})$, that we also denote by  $\gV^\star_{g,i}$, so that $\spanR_{\scriptscriptstyle \real}T_i=\gV^\star_{g,i}$ (remember that right now we are assuming that the minimum-phase invariant zeros are real).\footnote{If the system is right invertible, it is possible to prove that there holds
$\gV^\star_{g,i}=\gV_g^\star+\gR^\star_i$.}

%The set $L_i$ may be considered as the extension of $E_i$ to the case where $\lambda$ is allowed to be an invariant zero. 

We have proved that $\spanR_{\scriptscriptstyle \real} E_i=\gR^\star_i$; in the same way, one would expect the identity $\spanR_{\scriptscriptstyle \real}L_i=\gV^\star_{g,i}$ to hold. However, it can be proved that this is not the case. In other words, $\gL_i \defi \spanR_{\scriptscriptstyle \real} L_i$ is not equal to $\spanR_{\scriptscriptstyle \real} T_i$ in general. In fact, when $\lambda$ is equal to an invariant zero, the system
\[
\bmat{cc} A-\lambda\,I & B \\ C& D \emat \bmat{c} v \\ w \emat= \bmat{c} 0 \\ e_i \emat
\]
may not admit solutions because in this case the Rosenbrock matrix might lose rank (and therefore its rows are no longer linearly independent). 
This happens when the row $[\begin{array}{cc} C_i & D_i\end{array}]$ becomes linearly dependent with the other rows. 

%
%\begin{example}
%{\em
%Let
%\[
%A=\bmat{ccc}
%0 & -2 & -2 \\ 0 & 0 & 8 \\ 0 & 1 & -5\emat,\quad B=\bmat{cc} 0 & 0 \\ 2 & 0 \\ 0 & -3 \emat,\quad C=\bmat{ccc}
%8 & -2 & 0 \\ 0 & 0 & -2 \emat, D=0_{2 \times 2}.
%\]
%This system is left and right invertible with an invariant zero at $-8$.
%The quadruple has the no invariant zeros, and
%\[
%\gR_1=\bsmat 1 && 0\\[1mm] 0 && 1 \\[1mm]  0 && 0 \esmat
%\]
%and therefore $\gV^\star_{g,1}=\gR_1$ is two dimensional. Let us calculate $\gL_1$. First, for $\mu\neq -8$, \textcolor{red}{verify that the matrix is right invertible so that the system with $[0;e_j]$ can be solved} we have
%\[
%\bsmat A-\mu\,I && B \\[1mm] C & D \esmat^\dagger=
%\bsmat A-\mu\,I && B \\[1mm] C & D \esmat^{-1}=
%\frac{1}{\mu+8} \bsmat
%-1 && 0 && 0 && 1 && 1 \\[1mm] 
%-4 && 0 && 0 && -\frac{1}{2} \mu && 4 \\[1mm] 
%0 && 0 && 0 && 0 && -\frac{1}{2}\,\mu-4\\[1mm] 
%-2\,\mu && \frac{1}{2}\,\mu+4 && 0 && -\frac{1}{4}\mu^2 && 4\,\mu+16 \\[1mm] 
%-\frac{4}{3} && 0 && -\frac{1}{3}\,\mu-\frac{8}{3} && 
%-\frac{1}{6}\,\mu && \frac{1}{6}\mu^2+\frac{13}{6}\mu+8 \esmat
%\]
%so that for all $\mu \neq -8$ the vectors $v$ satisfying
%\[
%\bsmat v \\[1mm] w \esmat=\bmat{cc} A-\mu\,I & B \\ C & D\emat^\dagger\bsmat 0 \\ \beta\,e_1 \esmat= \beta\,\bsmat 1 \\[1mm] -\frac{1}{2}\mu \\[1mm] 0 \\[1mm] -\frac{1}{4}\mu^2 \\[1mm] -\frac{1}{6}\mu \esmat
%\]
%for some $w$ are all vectors in the form $\bsmat \alpha \\[1mm] \beta  \\[1mm] 0 \esmat$. When $\mu=-8$, we have 
%\[
%\bsmat 0 \\ \beta\,e_1 \esmat \notin \ima \bsmat A-(-8)\,I && B \\[1mm] C & D \esmat
%\]
%so that $\gL_i=\gV^\star_{g,i}$.
%}
%\end{example}
%

\begin{example}
{\em
Consider the right invertible quadruple $(A,B,C,D)$ given by the matrices
\[
A=\bsmat
0 && 0 && -2 \\[1mm] 0 && -3 && 0 \\[1mm] 0 && 0 && 0\esmat,\quad B=\bsmat 3 && 2 \\[1mm] 0 && 0 \\[1mm] -1 && 2 \esmat,\quad C=\bsmat
0 && 0 && -2 \\[1mm] 2 && 0 && 0 \esmat,\quad D=0_{2 \times 2}.
\]
This quadruple has one invariant zero at $-3$. One can easily verify that  $(A,B,C_{(1)},D_{(1)})$ has the same invariant zero, and
$\gR_1=\bsmat 0\\[1mm] 0 \\[1mm] 1 \esmat$.
Moreover,
\beann
\ker \bsmat A-(-3)\,I && B \\[1mm] C_{(1)} && D_{(1)}\esmat=\ima \bsmat 0 && 0 \\[1mm] 1 && 0 \\[1mm] 0 && -8 \\[1mm] 
\hline \\[0mm] 0 && -10 \\[1mm] 0 &&7 \esmat
\eeann
 gives
$\gV^\star_{g,1}=\ima\bsmat 0 && 0 \\[1mm] 1 && 0 \\[1mm] 0 && 1 \esmat$.
The subspace $\gL_i$ is spanned by the vectors $v$ satisfying
$\bsmat A-\mu\,I && B \\[1mm] C && D\esmat\bsmat v \\[1mm] w \esmat=\bsmat 0 \\ \beta\,e_1 \esmat$
for some $w$, for arbitrary $\mu$ including the invariant zero. A calculation shows that when $\mu\neq -3$ we have
\[
\bmat{cc} A-\mu\,I & B \\ C & D\emat^{-1}\bmat{c} 0 \\ 0 \\ 0 \\ \delta \\ 0\emat=
\bmat{cccccc}
0 & 0 & 0 & 0 & \frac{1}{2} \\
0 & -\frac{1}{\mu+3} & 0 & 0 & 0 \\
0 & 0 & 0 & -\frac{1}{2} & 0 \\
\frac{1}{4} & 0 & -\frac{1}{4} & \frac{1}{8}\,\mu-\frac{1}{4} & \frac{1}{8}\mu \\
\frac{1}{8} & 0 & \frac{3}{8} & -\frac{3}{16}\mu-\frac{1}{8} & \frac{1}{16}\mu\emat\bmat{c} 0 \\ 0 \\ 0 \\ \delta \\ 0\emat=\delta\,\bmat{c} 0 \\ 0 \\-\frac{1}{2}\\[1mm]
\hline\\[-3mm]
 \frac{1}{8}\,\mu-\frac{1}{4}\\ -\frac{3}{16}\mu-\frac{1}{8}\emat.
\]
%which implies 
%\[
%\bmat{cc} A-\mu\,I & B \\ C & D\emat^\dagger\bmat{c} 0 \\ 0 \\ 0 \\ \delta \\ 0\emat=
%\delta\,\bmat{c} 0 \\ 0 \\-\frac{1}{2}\\[1mm]
%\hline\\[-3mm]
% \frac{1}{8}\,\mu-\frac{1}{4}\\ -\frac{3}{16}\mu-\frac{1}{8}\emat
%\]
When $\mu=-3$ we have
\[
\bmat{cc} A-(-3)\,I & B \\ C & D\emat^\dagger\bmat{c} 0 \\ 0 \\ 0 \\ \delta \\ 0\emat=\bmat{cccccc} 
0 & 0 & 0 &  0 & \frac{1}{2} \\[0mm]
0 & 0 & 0 & 0 & 0 \\[0mm] 
0 & 0 & 0 & -\frac{1}{2} & 0 \\[0mm] 
\frac{1}{4} & 0 & -\frac{1}{4} & -\frac{5}{8} & -\frac{3}{8} \\[0mm] 
\frac{1}{8} & 0 & \frac{3}{8} & \frac{7}{16} & -\frac{3}{16} \emat\bmat{c} 0 \\ 0 \\ 0 \\ \delta \\ 0\emat=\delta\,\bmat{c} 0 \\ 0 \\-\frac{1}{2}\\[1mm]
\hline\\[-3mm]
-\frac{5}{8} \\ \frac{7}{16}\emat,
\]

so that no other new vectors are added from the invariant zeros 
and $\gL_i=\ima \bsmat 0 \\[1mm] 0 \\[1mm] 1 \esmat$. Hence, 
in this case $\gL_i$ is strictly contained in $\gV^\star_{g,i}$.
}
\end{example}

This example shows the necessity to introduce the new subspace $\gL_i$. The following result is instrumental in proving that, for all $i \in \{1,\ldots,p\}$, the subspace $\gL_i$ is ``between'' $\gR^\star_i$ and $\gV^\star_{g,i}$, i.e., $\gR^\star_i \subseteq \gL_i \subseteq \gV^\star_{g,i}$ for all $i \in \{1,\ldots,p\}$.

\begin{lemma}
	\label{lem:Li_real}
For all $i \in \{1,\ldots,p\}$ we have
\[
\gL_i=\gR^\star_i+\sum_{\lambda \in {\real_g\cap}\gZ} \spanR_{\scriptscriptstyle \real}\widehat{R}_i(\lambda).
\]
\end{lemma}
\proof
We have the following chain of identities:
\beann
\gL_i\ns&\ns = \ns&\ns\spanR_{\scriptscriptstyle \real}L_i  =
\spanR_{\scriptscriptstyle \real} \Big(\bigcup_{\lambda\in \real_g} \widehat{R}_i(\lambda)\Big) \\
\ns&\ns = \ns&\ns \sum_{\lambda\in \real_g} \spanR_{\scriptscriptstyle \real}\widehat{R}_i(\lambda) = \sum_{\lambda\in \real_g\setminus \gZ} \spanR_{\scriptscriptstyle \real}\widehat{R}_i(\lambda)+ \sum_{\lambda \in{\real_g\cap}\gZ} \spanR_{\scriptscriptstyle \real}\widehat{R}_i(\lambda)\\
\ns&\ns = \ns&\ns \sum_{\lambda\in \real_g\setminus \gZ} {\gR}_i(\lambda)+ \sum_{\lambda \in{\real_g\cap}\gZ} \spanR_{\scriptscriptstyle \real}\widehat{R}_i(\lambda)= \gR^\star_i+\sum_{\lambda \in {\real_g\cap}\gZ} \spanR_{\scriptscriptstyle \real}\widehat{R}_i(\lambda).
\eeann
\endproof

\begin{theorem}
There holds
\[
\gR^\star_i \subseteq \gL_i \subseteq \gV^\star_{g,i}.
\]
\end{theorem}
\proof
From the previous result it is obvious that $\gR^\star_i \subseteq \gL_i$. Moreover, as already observed we have $\sum_{\lambda \in \gZ} {\gR}_i(\lambda)\subseteq \gV^\star_{g,i}$. Since  we have shown that ${\gR}_i(\lambda)\supseteq \widehat{R}_i(\lambda)$, we can conclude that $\sum_{\lambda \in \gZ} \spanR_{\scriptscriptstyle \real}\widehat{R}_i(\lambda)\subseteq \gV^\star_{g,i}$.
\endproof

\subsection{Problem \ref{pro02}A}
The counterpart of Lemma \ref{lemmac1} appears  to be written exactly as Lemma \ref{lemmac1} itself. However, recall that in Problem \ref{pro02}A the closed-loop eigenvalues are allowed to coincide with minimum-phase invariant zeros.

\begin{lemma}
\label{lemmac4}
Let $\nu_0=n-\nu_1-\ldots-\nu_p$.
Problem \ref{pro02}A is solvable if and only if there exist
\beann
&& v_{0,k} \in \gR(\lambda_{0,j})\qquad \forall \,k\in\{1,\ldots,\nu_0\}\\
&& v_{i,j} \in \widehat{R}_i(\lambda_{i,j})\qquad \forall \,i\in\{1,\ldots,p\} \;\; \forall \,j\in\{1,\ldots,\nu_i\}
\eeann
such that $\{v_{0,1},\ldots,v_{0,\nu_0},v_{1,1},\ldots,v_{1,\nu_1},\ldots,v_{p,1},\ldots,v_{p,\nu_p}\}$ is linearly independent.
\end{lemma}

\proof 
The proof follows directly from the one of Lemma \ref{lemmac1}. 
\endproof

We denote $\widehat{R}_{0}(\lambda)=\gR(\lambda)$ for notational conciseness.

\begin{theorem}
\label{cond2A}
Let $\nu_0=n-\nu_1-\ldots-\nu_p$.
Problem \ref{pro02}A is solvable if and only if
\bea
\label{pippo}
\dim \left(\sum_{(i,j)\in P} \spanR_{\scriptscriptstyle \real}\widehat{R}_{i}(\lambda_{i,j})\right) \ge \operatorname{card}\,P
\eea
for all $P$ in the power set $2^I$ where $I=\{(0,1),\ldots,(0,\nu_0),\ldots,(p,1),\ldots,(p,\nu_p)\}$.
\end{theorem}

\proof
In both statements of Lemma \ref{lemmac1} and \ref{lemmac4}, the sets $\widehat{R}_i(\lambda_{i,j})$ are involved. However, while in the case where $\lambda_{i,j}$ are not invariant zeros the span of $\widehat{R}_i(\lambda_{i,j})$ is equal to $\gR_i(\lambda_{i,j})$, when $\lambda_{i,j}$ coincide with invariant zeros this may not necessarily be the case.
\endproof

We notice that condition (\ref{pippo}) is very easy to check since, whenever $\lambda_{i,j}\notin \gZ_{g}$, we have $\spanR_{\scriptscriptstyle \real}\widehat{R}_{i}(\lambda_{i,j})=\gR_i(\lambda_{i,j})$. The parameterization of all the feedback matrices $F$ that solves Problem \ref{pro02}A, when the necessary and sufficient conditions in Theorem \ref{cond2A} are satisfied, 
can be carried out exactly as in Theorem \ref{theF1A}, recalling that this time the observable eigenvalues $\{\l_{i,j}\}_{i=1,\ldots,p,\;j =1,\ldots,\nu_i}$ may contain invariant zeros.

\subsection{Problem \ref{pro02}B}

Using the same argument, for Problem \ref{pro02}B the following results hold.

\begin{lemma}
\label{lemmac5}
Let $\nu_0=n-\nu_1-\ldots-\nu_p$.
Problem \ref{pro02}B is solvable if and only if there exist
\beann
&& v_{0,k} \in E_g\qquad \forall \,k\in\{1,\ldots,\nu_0\}\\
&& v_{i,j} \in \widehat{R}_i(\lambda_{i,j})\qquad \forall \,i\in\{1,\ldots,p\} \;\; \forall \,j\in\{1,\ldots,\nu_i\}
\eeann
such that $\{v_{0,1},\ldots,v_{0,\nu_0},v_{1,1},\ldots,v_{1,\nu_1},\ldots,v_{p,1},\ldots,v_{p,\nu_p}\}$ is linearly independent.
\end{lemma}

\proof The proof follows from that of Lemma \ref{lemmac2}.
\endproof

\begin{theorem}
\label{pro2Bsol}
Let $\nu_0=n-\nu_1-\ldots-\nu_p$.
Problem \ref{pro02}B is solvable if and only if
\beann
\dim \left(\gV^\star_g+ \sum_{(i,j)\in P} \spanR_{\scriptscriptstyle \real}\widehat{R}_{i}(\lambda_{i,j})\right) \ge \operatorname{card}\,P+\nu_0
\eeann
and 
\beann
\dim \left(\sum_{(i,j)\in P} \spanR_{\scriptscriptstyle \real}\widehat{R}_{i}(\lambda_{i,j})\right) \ge \operatorname{card}\,P
\eeann
for all $P$ in the power set $2^I$ where $I=\{(1,1),\ldots,(1,\nu_1),\ldots,(p,1),\ldots,(p,\nu_p)\}$.
\end{theorem}

The proof follows immediately from the one of Theorem \ref{pro1Bsol}. Likewise, the parameterization of all the feedback matrices that solve Problem \ref{pro02}B are given exactly as that in Theorem \ref{thepar1B}, with the only difference that the set $\{\l_{i,j}\}_{i=1,\ldots,p,\;j =1,\ldots,\nu_i}$ is allowed to contain invariant zeros.

\subsection{Problem \ref{pro02}C}
Let us now consider Problem \ref{pro02}C.

\begin{lemma}
Let $\nu_0=n-\nu_1-\ldots-\nu_p$.
Problem \ref{pro02}C is solvable if and only if there exist
\beann
&& v_{0,k} \in E_g \qquad \forall \,k\in\{1,\ldots,\nu_0\}\\
&& v_{i,j} \in L_i \qquad \forall \,i\in\{1,\ldots,p\} \;\; \forall \,j\in\{1,\ldots,\nu_i \}
\eeann
 such that $\{v_{0,1},\ldots,v_{0,\nu_0},v_{1,1},\ldots,v_{1,\nu_1},\ldots,v_{p,1},\ldots,v_{p,\nu_p}\}$ is linearly independent.
\end{lemma}

\begin{theorem}
\label{th:2C}
Let $\nu_0=n-\nu_1-\ldots-\nu_p$.
Problem \ref{pro02}C is solvable if and only if
\[
\dim \left(\gV^\star_g+ \sum_{i\in P} \gL_i\right) \ge \sum_{i\in P} \nu_i+\nu_0
\]
and
\[
\dim \left(\sum_{i\in P} \gL_i\right) \ge \sum_{i\in P} \nu_i
\]
for all $P$ in the power set $2^I$ where $I=\{0,1,\ldots,p\}$.
\end{theorem}

\proof The statement follows on recalling that $\spanR_{\scriptscriptstyle \real}L_i=\gL_i$ and using Corollary \ref{cor2}.
\endproof

\section{Solution of Problem 3}
\label{sec:problem3_real}
Recall that in Problem 3 we need to observe at most $\nu_i$ modes on the $i$-th output. 

\subsection{Problem \ref{pro03}A}
Finally, in this section, we solve the third problem, in which only the maximum number of eigenvalues is assigned. 
\begin{lemma}
\label{lemmac7}
Let $\bar{\nu}_0=n-\bar{\nu}_1-\ldots-\bar{\nu}_p$.
Problem \ref{pro03}A is solvable if and only if there exist
\beann
&& v_{0,k} \in \gR(\lambda_{0,k})\qquad \forall \,k\in\{1,\ldots,\bar{\nu}_0\}\\
&& v_{i,j} \in {\gR}_i(\lambda_{i,j})\qquad \forall \,i \in\{1,\ldots,p\} \;\; \forall \,j \in\{1,\ldots,\bar{\nu}_i\}
\eeann
such that $\{v_{0,1},\ldots,v_{0,\bar{\nu}_0},v_{1,1},\ldots,v_{1,\bar{\nu}_1},\ldots,v_{p,1},\ldots,v_{p,\bar{\nu}_p}\}$ is linearly independent.
\end{lemma}

\proof Differently from the other two cases, in Problem  \ref{pro03}A some modes associated to a particular output component may not appear. Therefore, it is easy to see that in this case the result in Lemma \ref{lemmac1} holds true for sets defined as $\widehat{R}_i(\lambda_{i,j})$ but where $\delta_{i,j}$ is allowed to be zero. It is obvious that such set coincides with the linear space $\gR_i(\lambda_{i,j})$.
Moreover, $\gR(\l_{i,j})$ is contained in $\gR_i(\l_{i,j})$ for every $i \in\{1,\ldots,p\}$, so that if the eigenvector $v_{i,j}$ associated with a certain eigenvalue $\l_{i,j}$ is in $\gR(\l_{i,j})$, it is also in 
$\gR_i(\l_{i,j})$, which implies that the condition can be expressed in terms of the problem data $\bar{\nu}_i$ instead of $\nu_i$.
% \textcolor{red}{(spiegare meglio che certi $v_{i,j}$ possono appartenere a $\gR(\l_{i,j})$ che e' contenuto in $\gR_i(\l_{i,j})$; in questo modo il modo $\l_{i,j}$ non si vede ma la condizione di ris. puo essere espressa in termini di $\bar{\nu}$.)}
\endproof

Notice that in view of the analogy between Lemma \ref{lemmac1} and Lemma \ref{lemmac7}, the necessary and sufficient solvability conditions for Problem \ref{pro03}A are exactly the same as those of Problem \ref{pro01}A.

\begin{theorem}
\label{th:3A}
Let $\bar \nu_0=n-\bar \nu_1-\ldots-\bar \nu_p$.
Problem \ref{pro03}A is solvable if and only if
\[
\dim \left(\sum_{(i,j)\in P} \gR_{i}(\lambda_{i,j})\right) \ge \operatorname{card}\,P
\]
for all $P$ in the power set $2^I$ where $I=\{(0,1),\ldots,(0,\bar \nu_0),\ldots,(p,1),\ldots,(p,\bar \nu_p)\}$.
\end{theorem}
\subsection{Problem \ref{pro03}B}
Let us now consider Problem \ref{pro03}B. The same argument given before justify the following.  

\begin{lemma}
\label{lemmac8}
Let $\bar \nu_0=n-\bar \nu_1-\ldots-\bar \nu_p$. Problem \ref{pro03}B is solvable if and only if there exist $ {\nu}_i \leq \bar{\nu}_i,\;i\in\{1,
\ldots,p\}$ and $\nu_0 = n- \nu_1-\ldots- \nu_p \geq \bar \nu_0$ and
\beann
&& v_{0,k} \in E_g \qquad \forall \,k\in\{1,\ldots, \nu_0\}\\
&& v_{i,j} \in {\gR}_i(\lambda_{i,j}) \qquad \forall \,i\in\{1,\ldots,p\} \;\; \forall \,j\in\{1,\ldots, \nu_i\}
\eeann
 such that $\{v_{0,1},\ldots,v_{0,\nu_0},v_{1,1},\ldots,v_{1,\nu_1},\ldots,v_{p,1},\ldots,v_{p,\nu_p}\}$ is linearly independent.
\end{lemma}

In the case of Problem \ref{pro03}B, we cannot express the statement of Lemma \ref{lemmac8} only in terms of the parameters of the problem, because in this case the $\bar{\nu}_i-\nu_i$ closed-loop modes that are not effectively visible on the $i$-th output are not necessarily closed-loop unobservable modes.
In other words, the conditions in Lemma \ref{lemmac8} are expressed in terms of the numbers of closed-loop eigenvalues effectively observable from each output component. Nevertheless, it is desirable to express the solvability conditions in terms of the problem data. The following theorem addresses this point.

\begin{theorem}
\label{th:3B}
Let $\bar{\nu}_0=n- \bar\nu_1-\ldots-\bar \nu_p$.
Problem \ref{pro03}B is solvable if and only if
\[
\dim \left(\gV^\star_g+ \sum_{(i,j)\in P} \gR_{i}(\lambda_{i,j})\right) \ge \operatorname{card}\,P+\bar{\nu}_0
\]

for all $P$ in the power set $2^I$ where $I=\{(1,1),\ldots,(1,\bar{\nu}_1),\ldots,(p,1),\ldots,(p,\bar{\nu}_p)\}$.
\end{theorem}

%\tb{
\proof 
The statement follows from Corollary \ref{corollary_fabrizio} by considering that $q=\sum_{i=1}^p\bar {\nu}_i$, $ h=
\dim \gV^\star_{g} $, $ k=n-\dim \gV^\star_{g}  $ and recalling that $\spanR_{\scriptscriptstyle \real}E_g=\gV^\star_{g}$.
\endproof
%}

\subsection{Problem \ref{pro03}C}
Finally we consider Problem \ref{pro03}C.
\begin{lemma}
Let $\bar \nu_0=n-\bar \nu_1-\ldots-\bar \nu_p$.
Problem \ref{pro03}C is solvable if and only if there exist
\beann
&& v_{0,k} \in E_g \qquad \forall \,k\in \{1,\ldots,\bar \nu_0\}\\
&& v_{i,j} \in T_i \qquad \forall \,i\in \{1,\ldots,p\} \;\; \forall \,j\in \{1,\ldots,\bar \nu_i\}
\eeann
 such that $\{v_{0,1},\ldots,v_{0,\bar \nu_0},v_{1,1},\ldots,v_{1,\bar \nu_1},\ldots,v_{p,1},\ldots,v_{p,\bar \nu_p}\}$ is  linearly independent.
\end{lemma}

\proof
This result follows from the definition of $T_i$, by noting  that
since $E_g \subset T_i$, if a vector $v_{i,j}$ belongs to $E_g$, it also belongs to $T_i$, so that the condition can be expressed in terms of $\bar{\nu}_i$.
\endproof

\begin{theorem}
\label{th:3C}
Let $\bar{\nu}_0=n-\bar{\nu}_1-\ldots-\bar{\nu}_p$.
Problem \ref{pro03}C is solvable if and only if
\bea
\label{cond3c}
\dim \left(\gV^\star_g+ \sum_{i\in P} \gV^\star_{g,i}\right) \ge \sum_{i\in P} \bar{\nu}_i+\bar{\nu}_0
\eea
for all $P$ in the power set $2^I$ where $I=\{1,2,\ldots,p\}$.
\end{theorem}

\proof The statement follows from Corollary \ref{cor2} on recalling that $\spanR_{\scriptscriptstyle \real}T_i=\gV^\star_{g,i}$ and considering that $\gV^\star_{g,i} \supseteq \gV^\star_g$ for all $i \in \{1,\ldots,p\}$. %\tr{(Stesso discorso di prima)}
\endproof

 From the conditions obtained above we can see that whenever the closed-loop eigenvalues must be chosen to be different from the minimum-phase invariant zeros, requiring that a certain exact number will be observable from a certain output is entirely equivalent to requiring that at most the same number will be observable from that output. This fact seems rather counterintuitive, because at first sight the second problem appears to be a relaxation of the first. Nevertheless we have shown that no extra degrees of freedom arise when we only specify an upper bound on the number of modes we can observe, unless the closed-loop eigenvalues are chosen from within the minimum-phase invariant zeros. Indeed, in such case, it is no longer true that  requiring that a certain number of modes will be observable from a certain output is  equivalent to requiring that at most the same number will be observable from that output.

\begin{corollary}
	\label{cor:3C}
Let $\bar{\nu}_0=n-\bar{\nu}_1-\ldots-\bar{\nu}_p$.
Problem \ref{pro03}C is solvable if and only if
\bea
\label{cons1}
\dim \left(\sum_{i\in P} \gV^\star_{g,i}\right) \ge \sum_{i\in P} \bar{\nu}_i+\bar{\nu}_0
\eea
and
\bea
\label{cons2}
\dim \gV^\star_{g} \ge \nu_0
\eea
for all $P$ in the power set $2^I\setminus \varnothing$ where $I=\{1,2,\ldots,p\}$.
\end{corollary}

\proof 
Since $\gV^\star_{g,i} \supseteq \gV^\star_g$ for all $i \in \{1,\ldots,p\}$,  (\ref{cond3c})
is equivalent to (\ref{cons1}) for all $P\in 2^{\{1,\ldots,p\}}\setminus \varnothing$ and, when $P=\varnothing$, (\ref{cond3c}) reduces to (\ref{cons2}). 
\endproof

\begin{remark}
\label{pizzardonejordan}
{\em
As repeatedly mentioned, in this paper we have restricted our attention to the case where no Jordan structures occur, both for the assignable and unassignable eigenvalues (invariant zeros). The case of non-trivial Jordan structures requires a slightly different machinery, which involves the computation of Jordan chains of generalized closed-loop eigenspaces. For example, in Theorem \ref{th1}, a spanning set for $\gR^\star$ in the case of possibly coincident eigenvalues $\l_1,\ldots,\l_r$ involves the null-space of the Rosenbrock pencil complemented with a suitable chain of subspaces obtained in a recursive way starting from those null-spaces, see \cite{Ntogramatzidis-14}. The other subspaces $\gR(\cdot)$, $\gR_i(\cdot)$, $\gR^\star_i$, $\gL_i$, $\gV^\star_g$, $\gV^\star_{g,i}$ defined in the previous sections have to be generalized accordingly. While this extension does not pose conceptual difficulties, it does not lead to further insight and it considerably increases the notational burden; for this reason it has not been considered in this paper. It is also worth noting that allowing the case of non-trivial Jordan chains for the assignable eigenstructure does not enlarge the set of solvable problems. 
Finally, we observe that the most general definition of state-to-output decoupling, which takes into account the case of possibly non-trivial Jordan forms, is the one given in Theorem \ref{thequiv}; the adaptation of its proof to the case of Jordan chains is trivial.
}
\end{remark}

\section{The complex case}
\label{sec:complex}
The case of complex conjugate closed-loop eigenvalues and invariant zeros is significantly more difficult than the real case. The reason for this is immediately clear when one thinks that, in a case where $\gR^\star=\{0\}$ and the system has a single complex conjugate pair of invariant zeros in $\complex_g$ with single multiplicity, we cannot render a single closed-loop mode unobservable, because the complex conjugate vectors that we extract to build the feedback must be in pairs. This fact alone suggests that Rado's theorem may not be applied directly, because an additional constraint has to be added in some situations.\\ 
Consider, for example, the minimum-phase system
\beann
A=\bsmat
-6 && 0 && 0 && 3 \\[1mm]
0 && 0 && 4 && 0 \\[1mm]
0 && -4 && 0 && 0 \\[1mm]
0 && 0 && 3 && 0
\esmat
\qquad B=\bsmat
0 && 0 \\[1mm]
0 && 0 \\[1mm]
4 && 0 \\[1mm]
-4 && -1
\esmat\qquad
C=\bsmat
-5 && 0 && 0 && -1 \\[1mm]
3 && 0 && 7 && -2
\esmat\qquad D=\bsmat
0 && 0 \\[1mm]
0 && -1
\esmat
\eeann
which has the following zeros $ \gZ_g=\{-21,\,-2+\i\sqrt{7},\,-2-\i\sqrt(7)\} $. We aim to solve Problem \ref{pro01}B with $ \nu_1 = 1 $ and $ \l_{1,1}=-3 $, $ \nu_2=1  $ and $ \l_{2,2}=-5 $, $ \nu_0=2 $. For this systems we have  
\beann
\gR_1(-3)=
\spanR_{\scriptscriptstyle\real} \left\{\bsmat
0 \\[1mm]
-\frac{4}{3}\\[1mm]
1 \\[1mm]
0
\esmat \right\}
\qquad
\gR_2(-5)=\spanR_{\scriptscriptstyle\real}\left\{\frac{1}{10}\bsmat
63 \\[1mm]
-8 \\[1mm]
10 \\[1mm]
21
\esmat\right\} \qquad
\gsV_g=\ima\bsmat
-1 && 0 && 0 \\[1mm]
0 && 1 && 0 \\[1mm]
0 && 0 && 1 \\[1mm]
5 && 0 && 0
\esmat.
\eeann
It is immediate to check that the conditions of Theorem \ref{pro1Bsol} hold. Nevertheless, the problem is not solvable by using a real feedback matrix $ F $. Indeed, denoting by $\bsmat V_1 \\[1mm] W_1 \esmat$ a basis matrix of $\ker\bsmat A-(-2+\i\sqrt{7}) I && B\\[1mm] C && D\esmat$ and by $\bsmat V_2 \\[1mm] W_2 \esmat$ a basis matrix of $\ker\bsmat A-(-2-\i\sqrt{7}) I && B\\[1mm] C && D\esmat$, both partitioned conformably with the Rosenbrock matrix,
it can be noted that $ \gR_1(-3)\subseteq\ima V_1+\ima V_2\subseteq\gsV_g $. Hence, in order to have the mode $ \l_{1,1}=-3 $ appearing on the first output, we should only consider a subspace of dimension 1 of 
$\ima V_1+\ima V_2$
%\[
%\ker\bsmat A-(-2+\i\sqrt{7}) I && B\\[1mm] C &&D\esmat+\ker\bsmat A-(-2-\i\sqrt{7}) I && B\\[1mm] C && D\esmat 
%\]
 which, evidently, implies that such a subspace cannot contain complex conjugate elements. Thus, it is impossible to extract from that subspace pairs of complex conjugate linearly independent vectors, which is a necessary condition to obtain a real feedback matrix $ F $.\\
In other words, Rado's theorem provides necessary and sufficient conditions for the extraction of a set of linearly independent vectors, but is does not ensure that such a basis contains complex vectors that are not in complex conjugate pairs.\\

In the rest of the paper, for the sake of simplicity and with no loss of generality (see Remark \ref{pizzardonecomplessi}), we assume that the arbitrary modes that we select are real. The invariant zeros are allowed to be in complex conjugate pairs. With this simplifying assumption in mind, the solvability conditions for Problems \ref{pro01}A, \ref{pro02}A  and \ref{pro03}A do not change, {provided that $ \spanR_{\scriptscriptstyle\complex}(\cdot) $ is used in place of $ \spanR_{\scriptscriptstyle\real}(\cdot) $  }.\\ 
The situation is different for Problem \ref{pro01}B. The following corollary of Theorem \ref{th1} shows that the use of complex conjugate closed-loop eigenvalues that are not invariant zeros has no influence in the span of all the possible $\gR(\lambda)$.

\begin{corollary}
\label{MLcor}
There holds
\[
\spanR_{\scriptscriptstyle \complex} \Big(\bigcup_{\lambda \in \complex\setminus \gZ} \gR(\lambda)\Big)=
\spanR_{\scriptscriptstyle \complex} \Big(\bigcup_{\lambda \in \real\setminus \gZ_{\real}} \gR(\lambda)\Big).
\]
\end{corollary}
%
%\textcolor{red}{(da qui vanno tolti gli zeri ovviamente)}
\proof 
We only need to prove that $\spanR_{\scriptscriptstyle \complex} \bigl(\bigcup_{\lambda \in \complex\setminus \gZ} \gR(\lambda)\bigr)\subseteq
\spanR_{\scriptscriptstyle \complex} \bigl(\bigcup_{\lambda \in \real\setminus \gZ_{\real}} \gR(\lambda)\bigr)$, the opposite inclusion being obvious.
We recall that in view of Theorem \ref{th1} we have
\[
\gR^\star=\spanR_{\scriptscriptstyle \real} \Big(\bigcup_{\lambda \in \real\setminus \gZ_{\real}} \gR(\lambda)\Big)
\]
and for all $\lambda\in \complex \setminus \gZ$ we have $\mathfrak{Re}\{\gR(\lambda)\}\subseteq \gR^\star$ and 
$\mathfrak{Im}\{\gR(\lambda)\}\subseteq \gR^\star$, because for the construction of a basis for $\gR^\star$ the values of the closed-loop eigenvalues are arbitrary (provided they form a self conjugate set of distinct values that are different from the invariant zeros). 
Let $\{v_1,\ldots,v_r\}$ be a basis for $\gR^\star$. Let $v\in \gR(\lambda)$, where $\lambda \in \complex$. Since $\mathfrak{Re}\{v\},\mathfrak{Im}\{v\}\in \gR^\star$, we can write
$\mathfrak{Re}\{v\}=\alpha_1\,v_1+\ldots+\alpha_r\,v_r$ and $\mathfrak{Im}\{v\}=\beta_1\,v_1+\ldots+\beta_r\,v_r$, where $\alpha_i,\beta_i \in \real$ for $i\in \{1,\ldots,r\}$. Thus, $v=(\alpha_1+i\,\beta_1)\,v_1+\ldots+(\alpha_r+i\,\beta_r)\,v_r$.
\endproof

%The statement follows recalling that a basis for $\real^n$ is also a basis for the vector space $\complex^n$ over $\complex$.

Let $\gZ_{g,\complex}$ denote the set of invariant zeros in $\complex_g \setminus \real$.
%\begin{theorem}{\sc [Kimura's Theorem Extended]}\\
%\label{kimura2}
%Let $\lambda\in \complex$. Let $i \in \{1,\ldots,p\}$.
%Consider sets $\widehat{R}_i(\lambda),\widehat{R}_i(\overline{\lambda}),A_3,\ldots,A_n\subseteq \complex^n$ and let $I$ denote an independence relation. It is possible to find a set of independent vectors $\{\xi_1,\ldots,\xi_n\}$ such that 
%$\xi_1\in \widehat{R}_i(\lambda)$, $\xi_2\in \widehat{R}_i(\overline{\lambda})$, $\xi_3\in A_3$, $\ldots$, $\xi_n\in A_n$ and $\xi_1=\overline{\xi}_2$
%if and only if for any set $S\subseteq \{1,\ldots,n\}$ of cardinality $s=\operatorname{card}\,(S)$ 
%\[
%\dim \sum_{i \in S} \spanR\{A_i\} \ge \operatorname{card}\,(S).
%\]
%%Moreover, for any pair $A_i, A_j$ that are linear subspaces such that $A_i=\overline{A}_j$ it is possible to guarantee that the further constraint $\xi_i=\overline{\xi}_j$ is satisfied.
%\end{theorem}
%%
%\proof
%The proof is constructive and follows the same argument of \cite[Lemma 1]{Kimura-75}. 
%Let $\{\xi_1,\ldots,\xi_n\}$ be a linearly independent set of vectors such that 
%$\xi_1\in \widehat{R}_i(\lambda)$, $\xi_2\in \widehat{R}_i(\overline{\lambda})$, $\xi_3\in A_3$, $\ldots$, $\xi_n\in A_n$, whose existence is ensured by Rado's theorem. 
%
%
%\endproof
The following result is a counterpart of Corollary \ref{MLcor}, and can be proved using the same argument.

\begin{lemma}
\label{lemss}
There holds
\bea
\label{eqf1}
\spanR_{\scriptscriptstyle \complex} \Big(\bigcup_{\lambda \in \complex\setminus \gZ_{\complex}} \gR(\lambda)\Big)=
\spanR_{\scriptscriptstyle \complex}\Big( \bigcup_{\lambda \in \real}\gR(\lambda)\Big).
\eea
\bea
\label{eqf2}
\spanR_{\scriptscriptstyle \complex} \Big(\bigcup_{\lambda \in \complex_{{g}}\setminus \gZ_{\complex}} \gR(\lambda)\Big)=
\spanR_{\scriptscriptstyle \complex} \Big(\bigcup_{\lambda \in \real_{{g}}} \gR(\lambda)\Big).
\eea
\end{lemma}

We begin defining the set
\beann
E_g \ns&\ns = \ns&\ns \bigcup_{\lambda \in \complex_{{g}}} \gR(\lambda) \\
\ns&\ns = \ns&\ns
 \left\{v\in \complex^n\,\Big|\,\exists \,\lambda\in \complex_{{g}},\; \;\exists \,w\in \complex^m\,:
\bmat{cc} A-\lambda\,I & B \\ C & D \emat \bmat{c} v \\ w \emat=0\right\}\\
\ns&\ns = \ns&\ns E_{g,0} \cup \bigcup_{\lambda \in \gZ_{g,\complex}} \gR(\lambda)
\eeann
in $\complex^n$, 
where $E_{g,0}=\bigcup_{\lambda \in \complex_{{g}} \setminus \gZ_{g,\complex}} \gR(\lambda)$. If there are $c$ pairs of complex conjugate invariant zeros in  $\gZ_{g,\complex}$, we may write 
\[
\bigcup_{\lambda \in \gZ_{g,\complex}} \gR(\lambda)=E_{g,1} \cup E_{g,2}\cup \ldots\cup E_{g,2\,c},
\] 
where the $E_{g,i}$ are conformably indexed, i.e., 
where for all odd $i \in \{1,\ldots,2\,c-1\}$ we have $E_{g,i}=\overline{E}_{g,i+1}$. {By considering the definition of $ \gR(\l) $, it is immediate to note that $ E_{g,i}=\gR(\l_i) $ and $ E_{g,i+1}=\gR(\l_{i+1}) $, for all odd $i \in \{1,\ldots,2\,c-1\}$ such that $ \l_i=\bar\l_{i+1} $ and $ \l_i,\l_{i+1}\in \gZ_{g,\complex} $.}

\subsection{Problem \ref{pro01}}
We address in this section the solution of Problems \ref{pro01}B-\ref{pro01}C. Following the same structure used in Section \ref{sec:problem1_real}, we first propose the solution in terms of existence of linearly independent vectors and then in terms of dimension of suitable subspaces.    
\subsubsection{Problem \ref{pro01}B} 
Since, in Lemma \ref{lemmac2}, $\nu_0$ closed-loop eigenvectors are chosen from $E_g$, in the complex case there must exist $\nu_{0,0},\nu_{0,1},\ldots,\nu_{0,2\,c}$, with $\nu_{0,i}=\nu_{0,i+1}$ for each odd $i$, such that $\nu_0=\sum_{i=0}^{2\,c} \nu_{0,i}$.
Lemma \ref{lemmac2} is modified in the complex case as follows.
\begin{lemma}
\label{lemmac21}
Let $\nu_0=n-\nu_1-\nu_2-\ldots-\nu_p$.
Problem \ref{pro01}B is solvable if and only if there exist
\beann
&& v_{0,0,1},\ldots,v_{0,0,\nu_{0,0}} \in E_{g,0} \\
&& v_{0,i,j} \in E_{g,i} \\%\qquad \, \\
&& v_{0,i+1,j} =\overline{v}_{0,i,j} \in E_{g,i+1}=\overline{E}_{g,i}  \qquad i\in \{1,3,\ldots,2\,c-1\},\,j=\{1,\ldots,\nu_{0,i}\} \\
&& v_{i,j} \in \widehat{R}_i(\lambda_{i,j}) \qquad i\in\{1,\ldots,p\}, \;\; \,j\in\{1,\ldots,\nu_i\} 
\eeann
which are all linearly independent and such that $v_{0,0,1},\ldots,v_{0,0,\nu_{0,0}}$  are real.
\end{lemma}
\proof The only point that needs to be proved is the requirement that $v_{0,0,1},\ldots,v_{0,0,\nu_{0,0}}$ are real. Since $E_{g,0}$ is in $\complex^n$, but we want to obtain a real feedback, we can choose the vectors $v_{0,0,1},\ldots,v_{0,0,\nu_{0,0}}$ to be either real or in complex conjugate pairs. 
However, in view of Lemma \ref{lemss} and Corollary \ref{cor4}, selecting these vectors to be real or in complex conjugate pairs is irrelevant.
\endproof

In the previous lemma, the vectors were complex, because $E_{g,0},E_{g,1} ,\ldots,E_{g,2\,c}$ are sets in $\complex^n$.
This does not constitute an issue for the vectors in $E_{g,1} ,\ldots,E_{g,2\,c}$, because they will result in complex conjugate pairs. The problem lies in the vectors that we are free to choose from within the set $E_{g,0}$. In other words, when
using Rado's theorem, we learn that our ability to choose linearly independent vectors $\{v_{0,1},\ldots,v_{0,\nu_{0,0}}\}$ depends on the span, with complex coefficients, of $E_{g,0}$. On the other hand, Corollary \ref{MLcor} ensures that $\spanR_{\scriptscriptstyle \complex} E_{g,0}$ coincides with the span that is obtained by restricting ourselves to real values of $\lambda$. Thus, we have the following intermediate result.

\begin{lemma}
Let $\nu_0=n-\nu_1-\ldots-\nu_p$.
Problem \ref{pro01}B is solvable if and only if there exist
\beann
&& v_{0,0,1},\ldots,v_{0,0,\nu_{0,0}} \in \gE_{g,0}=\spanR_{\scriptscriptstyle \complex} E_{g,0}  \\
&& \vdots \\
&& v_{0,i,1},\ldots,v_{0,i,\nu_{0,i}} \in \gE_{g,1}=\spanR_{\scriptscriptstyle \complex} E_{g,i}   \\%\qquad \forall \,j=1,\ldots,\nu_{0,0} \\
&& \overline{v}_{0,i,1},\ldots,\overline{v}_{0,i,\nu_{0,i}} \in \overline{\gE}_{g,i}=\spanR_{\scriptscriptstyle \complex} \overline{E}_{g,i}   \\
&& \vdots \\
&& v_{i,j} \in \widehat{R}_i(\lambda_{i,j}) \qquad\,i\in\{1,\ldots,p\} \;\; j\in\{1,\ldots,\nu_i\}
\eeann
with $v_{0,j,k}=\overline{v}_{0,j+1,k}$ for all $1\leq k \leq \nu_{0,j}=\nu_{0,j+1}  $ and for each odd $j\in \{1,\ldots,2\,c\}$,
 such that $\{v_{0,1},\ldots,v_{p,\nu_p}\}$ are linearly independent and such that $v_{0,0,1},\ldots,v_{0,0,\nu_{0,0}}$ are real.
\end{lemma}

\begin{theorem}
\label{leFE}
Problem \ref{pro01}B is solvable if and only if there exist $\nu_{0,0},\nu_{0,1},\nu_{0,3},\ldots,\nu_{0,2c-1}\in \mathbb{N}$ such that $\nu_{0,0}+2\nu_{0,1}+2\nu_{0,3}+\ldots+2\nu_{0,2c-1}=\nu_{0}$ and
\bea
\label{cond1B_comp}
 \dim \left(\sum_{i\in Q_0}\gE_{g,0}+\sum_{i\in Q}\gE_{g,i}+\sum_{i\in Q^\prime}\overline{\gE}_{g,i}+\sum_{(i,j)\in P} \gR_i^c(\lambda_{i,j})\right) \ge 
\sum_{i\in Q_0} \nu_{0,i}+\sum_{i\in Q} \nu_{0,i}+\sum_{i\in Q^\prime} \nu_{0,i}
+\operatorname{card} P\qquad 
\eea
for all 
\begin{itemize}
\item $Q_0\in 2^{J_0}$ with $J_0=\{0\}$, 
\item $Q,Q^\prime\in 2^J$ with $J= \{1,3,\ldots,2c-1\}$, 
\item $P\in 2^I$, where $I=\{(1,1),\ldots,(1,\nu_1),\ldots,(p,1),\ldots,(p,\nu_p)\}$,
\end{itemize}
 and where $\gR_i^c(\lambda_{i,j})=\spanR_{\scriptscriptstyle \complex} \widehat{R}_i(\lambda_{i,j})$.
\end{theorem}
\proof {The result follows directly from Theorems \ref{kimura} and \ref{kimura1} by considering that i) $ E_{g,0} $ is a real set, ii) $ E_{g,i}=\gR(\l_i) $ and $ \bar{E}_{g,i}=E_{g,i+1}=\gR(\l_{i+1}) $ are subspaces, thus also affine subspaces and iii) that, since $ \lambda_{i,j}\in\real$, then $ \widehat R_i(\l_{i,j})$ always contains a real set $\real \supseteq Q \subseteq \widehat R_i(\l_{i,j}) $ such that $\spanR_{\scriptscriptstyle \complex} Q = \spanR_{\scriptscriptstyle \complex} \widehat R_i(\l_{i,j}) $. The first two points are obvious; the third one follows immediately by noting that, from its definition, the set $ \widehat R_i(\l_{i,j}) $ always comprises pairs of complex conjugate elements. For every set containing complex conjugate pairs there exists a real subset such that their complex spans coincide.} 
\endproof

Notice that if the conditions of Theorem \ref{leFE} are satisfied, $\dim \gE_{g,i}\ge \nu_{0,i}$ for all $i\in\{0,\ldots,c\}$.
\begin{remark}{\em 
The construction of the feedback in this case can be carried out by following the same procedure given in Theorem \ref{thepar1B}, where now the values $\l_{0,k}$ are allowed to also be in $\gZ_g$, with the constraint that if a complex value is chosen, its complex conjugate, say $\l_{0,\ell}$ is also chosen. Moreover, if $V_{0,k}$ and $W_{0,k}$ are such that $\bsmat V_{0,k} \\[1mm] W_{0,k}\esmat$ is a basis matrix for $\ker \bsmat A-\l_{0,k}\,I && B \\[1mm] C && D \esmat$, then $V_{0,\ell}=\overline{V}_{0,k}$ and $W_{0,\ell}=\overline{W}_{0,k}$ are such that $\bsmat V_{0,\ell} \\[1mm] W_{0,\ell}\esmat$ is a basis matrix for $\ker \bsmat A-\l_{0,\ell}\,I && B \\[1mm] C && D \esmat$.
Hence, the parameters $k_{0,k}$ and $k_{0,\ell}$ have to be chosen to be complex conjugate, so that constructing the matrices $V_{k_{i,j}}$ and $W_{k_{i,j}}$ as in Theorem  \ref{thepar1B}, the corresponding feedback matrix $F$ is real as shown for example in the proof of \cite[Proposition 1]{Moore-76}.
} 
\end{remark}

\subsubsection{Problem \ref{pro01}C}
%\\
\begin{lemma}
	\label{lem:1C_comp}
	Let $\nu_0=n-\nu_1-\nu_2-\ldots-\nu_p$.
	Problem \ref{pro01}B is solvable if and only if there exist
	\beann
	&& v_{0,0,1},\ldots,v_{0,0,\nu_{0,0}} \in E_{g,0} \\
	&& v_{0,i,j} \in E_{g,i} \\%\qquad  \\
	&& v_{0,i+1,j}=\overline{v}_{0,i,j} \in E_{g,i+1}=\overline{E}_{g,i}  \qquad i\in \{1,3,\ldots,2\,c-1\},\,j=\{1,\ldots,\nu_{0,i}\} \\
	&& v_{i,j} \in \gsR_i \qquad i\in\{1,\ldots,p\}, \;\; \,j\in\{1,\ldots,\nu_i\} 
	\eeann
	which are all linearly independent and such that $v_{0,0,1},\ldots,v_{0,0,\nu_{0,0}}$  are real.
\end{lemma}
\proof
The proof follows immediately form the one of Lemma \ref{lemmac21} by considering that now we can select arbitrary vectors from  $ \gsR_i $ since the model $ \l_{i,j} $ are not assigned.
\endproof

\begin{theorem}
	\label{thm:1C_comp}
	Problem \ref{pro01}C is solvable if and only if there exist $\nu_{0,0},\nu_{0,1},\nu_{0,3},\ldots,\nu_{0,2c-1}\in \mathbb{N}$ such that $\nu_{0,0}+2\nu_{0,1}+2\nu_{0,3}+\ldots+2\nu_{0,2c-1}=\nu_{0}$ and
	\bea
	\dim \left(\sum_{i\in Q_0}\gE_{g,0}+\sum_{i\in Q}\gE_{g,i}+\sum_{i\in Q^\prime}\overline{\gE}_{g,i}+\sum_{i\in P} \gsR_i{}^c\right) \ge 
	\sum_{i\in Q_0} \nu_{0,i}+\sum_{i\in Q} \nu_{0,i}+\sum_{i\in Q^\prime} \nu_{0,i}
	+\sum_{i\in P}\nu_i\qquad 
	\eea
	for all 
	\begin{itemize}
		\item $Q_0\in 2^{J_0}$ with $J_0=\{0\}$, 
		\item $Q,Q^\prime\in 2^J$ with $J= \{1,3,\ldots,2c-1\}$, 
		\item $P\in 2^I$, where $I=\{1,2,\ldots, p\}$,
	\end{itemize}
	and where $\gsR_i{}^c=\spanR_{\scriptscriptstyle \complex} \gsR_i$.
\end{theorem}
\proof
The result follows naturally form the proof of Theorem \ref{thm:1C_comp} by noting that the set $ \gsR_i{}^c $ always comprises pairs of complex conjugate elements.
\endproof

\subsection{Problem \ref{pro02}}
This section is devoted to the solution of Problems \ref{pro02}B-\ref{pro02}C in the presence of complex-conjugate zeros. The necessary subspaces will be defined along the section.  
\subsubsection{Problem \ref{pro02}B}
\begin{lemma}
	\label{lem:2B_complex}
	Let $\nu_0=n-\nu_1-\nu_2-\ldots-\nu_p$.
	Problem \ref{pro02}B is solvable if and only if there exist
	\beann
	&& v_{0,0,1},\ldots,v_{0,0,\nu_{0,0}} \in E_{g,0} \\
	&& v_{0,i,j} \in E_{g,i} \\
	&& v_{0,i+1,j}=\overline{v}_{0,i,j} \in E_{g,i+1}=\overline{E}_{g,i}  \qquad i\in \{1,3,\ldots,2\,c-1\},\,j=\{1,\ldots,\nu_{0,i}\}\\
	&& v_{i,j} \in \widehat{R}_i(\lambda_{i,j}) \qquad i\in\{1,\ldots,p\}, \;\; \,j\in\{1,\ldots,\nu_i\} 
	\eeann
	which are all linearly independent and such that $v_{0,0,1},\ldots,v_{0,0,\nu_{0,0}}$  are real and either $ v_{i,j}=\bar v_{i,k} $ if $ \l_{i,j}\in \gZ_{g,\complex}  $ or $ v_{i,j} $ is real if $ \l_{i,j}\in \real $.
\end{lemma}

\begin{theorem}
	\label{thm:2B_complex}
	Problem \ref{pro02}B is solvable if and only if there exist $\nu_{0,0},\nu_{0,1},\nu_{0,3},\ldots,\nu_{0,2c-1}\in \mathbb{N}$ such that $\nu_{0,0}+2\nu_{0,1}+2\nu_{0,3}+\ldots+2\nu_{0,2c-1}=\nu_{0}$ and
	\bea
	&&\dim \left(\sum_{i\in Q_0}\gE_{g,0}+\sum_{i\in Q}\gE_{g,i}+\sum_{i\in Q^\prime}\overline{\gE}_{g,i}+\sum_{(i,j)\in P_\real} \gR_i^c(\lambda_{i,j})+\sum_{(i,j)\in P_\complex} \gR_i^c(\lambda_{i,j})\sum_{(i,j)\in P_\complex^\prime} \overline{\gR}_i^c(\lambda_{i,j})\right) \nonumber\\ 
	&&\ge \sum_{i\in Q_0} \nu_{0,i}+\sum_{i\in Q} \nu_{0,i}+\sum_{i\in Q^\prime} \nu_{0,i}
	+\operatorname{card} P_\real+\operatorname{card} P_\complex +\operatorname{card} P_\complex^\prime 
	\eea
	for all 
	\begin{itemize}
		\item $Q_0\in 2^{J_0}$ with $J_0=\{0\}$, 
		\item $Q,Q^\prime\in 2^J$ with $J= \{1,3,\ldots,2c-1\}$, 
		\item $P_\complex,P_\complex^\prime\in 2^{I_\complex}$, where $I_\complex=\{(i,j)\;|\;i,j\in\mathbb{N},i\leq p,j\leq\nu_i,\l_{i,j}\in\gZ_{g,\complex},\mathfrak{Im}\l_{i,j}<0 \}$,
		\item $P_\real\in 2^{I_\real}$, where $I_\real=\{(i,j)\;|\;i,j\in\mathbb{N},i\leq p,j\leq\nu_i,\l_{i,j}\in\real\}$,
	\end{itemize}
	and where $\gR_i^c(\lambda_{i,j})=\spanR_{\scriptscriptstyle \complex} \widehat{R}_i(\lambda_{i,j})$.
\end{theorem}

\subsubsection{Problem \ref{pro02}C} 
In order to address Problem \ref{pro02}C, we consider the generalization of the set $ L_i $ to the complex case. We define
\beann
L_i & \defi & \left\{v\in \complex^n\,\Big|\,\exists \,\lambda\in \complex_g,\;\,\exists \,w\in \complex^m,\;\;\exists \;\delta \in \real\setminus \{0\}:\bmat{cc} A-\lambda\,I & B \\ C& D \emat \bmat{c} v \\ w \emat= \bmat{c} 0 \\ \delta\,e_i \emat\right\}.
\eeann
It is immediate to note that Lemma \ref{lem:Li_real} generalizes to the complex case yielding  
\begin{lemma}
	\label{lem:Li_complex}
	For all $i \in \{1,\ldots,p\}$ we have
	\[
	\gL_i=\gR^{\star c}_i+\sum_{\lambda \in {\real_g\cap}\gZ} \spanR_{\scriptscriptstyle \complex}\widehat{R}_i(\lambda)+\sum_{\lambda \in \gZ_{g,\complex}} \spanR_{\scriptscriptstyle \complex}\widehat{R}_i(\lambda).
	\]
\end{lemma}
\proof
The result can be proven using exactly the same procedure employed  in the proof of Lemma \ref{lem:Li_real}. 
\endproof
Following the same approach used in the definition of the set $ E_g $, we decompose the set $ L_i $ into smaller subsets in order to apply Theorems \ref{kimura}-\ref{kimura1}. We can conveniently represent the set $ L_i $ as
\[
L_i=L_{i,0} \cup \bigcup_{\lambda \in \gZ_{g,\complex}} \widehat{R}_i(\lambda)
\]
where $L_{i,0} = \bigcup_{\l\in \complex_g\setminus\gZ_{g,\complex}} \widehat{R}_i(\lambda) $. If there are $c$ pairs of complex conjugate invariant zeros in  $\gZ_{g,\complex}$, we may write 
\[
\bigcup_{\lambda \in \gZ_{g,\complex}} \widehat{R}_i(\lambda)=L_{i,1} \cup L_{i,2}\cup \ldots\cup L_{i,2\,c},
\] 
where the $L_{i,j}$ are conformably indexed, i.e., where for all odd $j \in \{1,\ldots,2\,c-1\}$ we have $L_{i,j}=\overline{L}_{i,j+1}$.\\
Since $\nu_i$ closed-loop eigenvectors are chosen from each $L_i$, in the complex case there must exist $\nu_{i,0},\nu_{i,1},\ldots,\nu_{i,2\,c}$, with $\nu_{i,j}=\nu_{i,j+1}$ for each odd $j$, such that $\nu_i=\sum_{j=0}^{2\,c} \nu_{i,j}$. \begin{lemma}
	\label{lem:2C_complex}
	Let $\nu_0=n-\nu_1-\nu_2-\ldots-\nu_p$.
	Problem \ref{pro02}C is solvable if and only if there exist
	\beann
	&& v_{0,0,1},\ldots,v_{0,0,\nu_{0,0}} \in E_{g,0} \\
	&& v_{0,j,k} \in E_{g,j}\\
	&& v_{0,j+1,k}=\overline{v}_{0,j,k} \in E_{g,j+1}=\overline{E}_{g,j}  \qquad j\in \{1,3,\ldots,2\,c-1\} ,\,k\in\{1,\ldots,\nu_{0,j}\}\\ 
	&& v_{i,0,k},\ldots,v_{i,0,\nu_{i,0}} \in L_{i,0} \\
	&& v_{i,j,k} \in L_{i,j}\\
	&& v_{i,j+1,k}=\overline{v}_{i,j,k} \in L_{i,j+1}=\overline{L}_{i,j} \qquad i\in\{1,\ldots,p\},\, j\in\{1,3,\ldots,2\,c-1\},\,k\in\{1,\ldots,\nu_{i,j}\} \\
	\eeann
	which are all linearly independent and such that $v_{i,0,1},\ldots,v_{i,0,\nu_{i,0}}$  are real.
\end{lemma}

\begin{theorem}
	\label{thm:2C_complex}
	Problem \ref{pro02}C is solvable if and only if there exist $\nu_{i,0},\nu_{i,1},\nu_{i,3},\ldots,\nu_{i,2c-1}\in \mathbb{N}$ such that $\nu_{i,0}+2\nu_{i,1}+2\nu_{i,3}+\ldots+2\nu_{i,2c-1}=\nu_{i}$ and
	\bea
	&&\dim \left(\sum_{j\in Q_0}\gE_{g,0}+\sum_{j\in Q}\gE_{g,i}+\sum_{j\in Q^\prime}\overline{\gE}_{g,i}+\sum_{i\in P_0} \gL_{i,0}+\sum_{(i,j)\in P} \gL_{i,j}+\sum_{(i,j)\in P^\prime} \overline{\gL}_{i,j}\right) \nonumber\\ 
	&&\ge \sum_{j\in Q_0} \nu_{0,j}+\sum_{j\in Q} \nu_{0,j}+\sum_{j\in Q^\prime} +\nu_{0,j}+\sum_{i\in P_0} \nu_{i,0}+\sum_{(i,j)\in P} \nu_{i,j}+\sum_{(i,j)\in P^\prime} \nu_{i,j} 
	\eea
	for all 
	\begin{itemize}
		\item $Q_0\in 2^{J_0}$ with $J_0=\{0\}$, 
		\item $Q,Q^\prime\in 2^J$ with $J= \{1,3,\ldots,2c-1\}$, 
		\item $P_0\in 2^{J_0}$ with $J_0=\{1,\ldots,p\}$, 
		\item $P,P^\prime\in 2^J$ with $J= \{(1,1),(1,3),\ldots,(1,2c-1),(2,1),\ldots,(p,2c-1)\}$,
	\end{itemize}
	where $ \gL_{i,j}\defi\spanR_{\scriptscriptstyle \complex}L_{i,j} $.
\end{theorem}

\subsection{Problem \ref{pro03}}
Finally, in this section, we address Problems \ref{pro03}B-\ref{pro03}C. Again, the necessary subspaces will be generalized to the complex case along the section.

\subsubsection{Problem \ref{pro03}B}%\
%\\
\begin{lemma}
	\label{lem:3B_complex}
	Let $\nu_0=n-\nu_1-\nu_2-\ldots-\nu_p$.
	Problem \ref{pro03}B is solvable if and only if there exist $ {\nu}_i \leq \bar{\nu}_i,\;i\in\{1,
	\ldots,p\}$ and $\nu_0 = n- \nu_1-\ldots- \nu_p \geq \bar \nu_0$ and
	\beann
	&& v_{0,0,1},\ldots,v_{0,0,\nu_{0,0}} \in E_{g,0} \\
	&& v_{0,i,j} \in E_{g,i} \\
	&& v_{0,i+1,j}=\overline{v}_{0,i,j} \in E_{g,i+1}=\overline{E}_{g,i}  \qquad i\in \{1,3,\ldots,2\,c-1\},\,j=\{1,\ldots,\nu_{0,i}\}\\
	&& v_{i,j} \in \widehat{R}_i(\lambda_{i,j}) \qquad i\in\{1,\ldots,p\}, \;\; \,j\in\{1,\ldots,\nu_i\} 
	\eeann
	which are all linearly independent and such that $v_{0,0,1},\ldots,v_{0,0,\nu_{0,0}}$  are real and either $ v_{i,j}=\bar v_{i,k} $ if $ \l_{i,j}\in \gZ_{g,\complex}  $ or $ v_{i,j} $ is real if $ \l_{i,j}\in \real $.
\end{lemma}

\begin{theorem}
	\label{thm:3B_complex}
	Problem \ref{pro02}B is solvable if and only if there exist $ {\nu}_i \leq \bar{\nu}_i,\;i\in\{1,
	\ldots,p\}$, $\nu_0 = n- \nu_1-\ldots- \nu_p \geq \bar \nu_0$, and $\nu_{0,0},\nu_{0,1},\nu_{0,3}\ldots,\nu_{0,2c-1}\in \mathbb{N}$ such that $\nu_{0,0}+2\nu_{0,1}+2\nu_{0,3}+\ldots+2\nu_{0,2c-1}=\nu_{0}$ and
	\bea
	&&\dim \left(\sum_{i\in Q_0}\gE_{g,0}+\sum_{i\in Q}\gE_{g,i}+\sum_{i\in Q^\prime}\overline{\gE}_{g,i}+\sum_{(i,j)\in P_\real} \gR_i^c(\lambda_{i,j})+\sum_{(i,j)\in P_\complex} \gR_i^c(\lambda_{i,j})\sum_{(i,j)\in P_\complex^\prime} \overline{\gR}_i^c(\lambda_{i,j})\right) \nonumber\\ 
	&& \ge \sum_{i\in Q_0} \nu_{0,i}+\sum_{i\in Q} \nu_{0,i}+\sum_{i\in Q^\prime} \nu_{0,i}
	+\operatorname{card} P_\real+\operatorname{card} P_\complex +\operatorname{card} P_\complex^\prime
	\eea
	for all 
	\begin{itemize}
		\item $Q_0\in 2^{J_0}$ with $J_0=\{0\}$, 
		\item $Q,Q^\prime\in 2^J$ with $J= \{1,3,\ldots,2c-1\}$, 
		\item $P_\complex,P_\complex^\prime\in 2^{I_\complex}$, where $I_\complex=\{(i,j)\;|\;i,j\in\mathbb{N},i\leq p,j\leq\nu_i,\l_{i,j}\in\gZ_{g,\complex},\mathfrak{Im}\l_{i,j}<0 \}$,
		\item $P_\real\in 2^{I_\real}$, where $I_\real=\{(i,j)\;|\;i,j\in\mathbb{N},i\leq p,j\leq\nu_i,\l_{i,j}\in\real\}$,
	\end{itemize}
	and where $\gR_i^c(\lambda_{i,j})=\spanR_{\scriptscriptstyle \complex} \gR_i(\lambda_{i,j})$.
\end{theorem}

\subsubsection{Problem \ref{pro03}C} In order to address the last problem we need to generalize the definition of the set $ T_i $ to the complex case
\[
T_i \defi \left\{v\in \complex^n\,\Big|\,\exists \,\lambda\in \complex_g,\;\;\exists \,w\in \complex^m\,:
\bmat{cc} A-\lambda\,I & B \\ C_{(i)} & D_{(i)}  \emat \bmat{c} v \\ w \emat= 0\right\}.
\]
Again, following the procedure previously employed to decompose $ E_g $ and $ L_i $, we can get $ T_i=T_{i,0}\cup\bigcup_{\lambda \in \gZ_{g,\complex}} \gR_i(\lambda)$ where
\[
\bigcup_{\lambda \in \gZ_{g,\complex}} \gR_i(\lambda)=T_{i,1} \cup T_{i,2}\cup \ldots\cup T_{i,2\,c},
\] 
and the $T_{i,j}$ are conformably indexed, i.e., where for all odd $j \in \{1,\ldots,2\,c-1\}$ we have $T_{i,j}=\overline{T}_{i,j+1}$.\\
We can now state the following lemma. It is worth stressing that the resolvability result can be provided in terms of the problem data $ \bar \nu_i $ because, in view of the right invertibility, we have that  $ E_{g,j}\subseteq T_{i,j},\;\forall i\in\{1,\ldots,p\} $. 
\begin{lemma}
\label{lem:3C_complex}
Let $\bar \nu_0=n-\bar \nu_1-\bar \nu_2-\ldots-\bar \nu_p$.
Problem \ref{pro03}C is solvable if and only if there exist $\bar \nu_{i,0},\bar\nu_{i,1},\bar\nu_{i,3}\ldots,\nu_{i,2c-1}\in \mathbb{N}$ such that $\bar\nu_{i,0}+2\bar\nu_{i,1}+2\bar\nu_{i,3}+\ldots+2\bar\nu_{i,2c-1}=\bar\nu_{i}$ 
\beann
&& v_{0,0,1},\ldots,v_{0,0,\nu_{0,0}} \in E_{g,0} \\
&& v_{0,j,k} \in E_{g,j}\\
&& v_{0,j+1,k}=\overline{v}_{0,j,k} \in E_{g,j+1}=\overline{E}_{g,j}  \qquad j\in \{1,3,\ldots,2\,c-1\} ,\,k\in\{1,\ldots,\nu_{0,j}\} \\
&& v_{i,0,k},\ldots,v_{i,0,\nu_{i,0}} \in T_{i,0} \\
&& v_{i,j,k} \in T_{i,j}\\
&& v_{i,j+1,k}=\overline{v}_{i,j,k} \in T_{i,j+1}=\overline{T}_{i,j} \qquad i\in\{1,\ldots,p\},\, j\in\{1,3,\ldots,2\,c-1\},\,k\in\{1,\ldots,\nu_{i,j}\} \\
\eeann
which are all linearly independent and such that $v_{i,0,1},\ldots,v_{i,0,\nu_{i,0}}$  are real.
\end{lemma}

\begin{theorem}
	\label{thm:3C_complex}
	Problem \ref{pro03}C is solvable if and only if there exist $\bar\nu_{i,0},\nu_{i,1},\bar\nu_{i,3},\ldots,\bar\nu_{i,2c-1}\in \mathbb{N}$ such that $\bar\nu_{i,0}+2\bar\nu_{i,1}+2\bar\nu_{i,3}+\ldots+2\bar\nu_{i,2c-1}=\bar\nu_{i}$ and
	\bea
	&&\dim \left(\sum_{j\in Q_0}\gE_{g,0}+\sum_{j\in Q}\gE_{g,i}+\sum_{j\in Q^\prime}\overline{\gE}_{g,i}+\sum_{i\in P_0} \gT_{i,0}+\sum_{(i,j)\in P} \gT_{i,j}+\sum_{(i,j)\in P^\prime} \overline{\gT}_{i,j}\right) \nonumber\\ 
	&&\ge\sum_{j\in Q_0} \bar\nu_{0,j}+\sum_{j\in Q} \bar\nu_{0,j}+\sum_{j\in Q^\prime} +\bar\nu_{0,j}+\sum_{i\in P_0} \bar\nu_{i,0}+\sum_{(i,j)\in P} \bar\nu_{i,j}+\sum_{(i,j)\in P^\prime} \bar\nu_{i,j} 
	\eea
	for all 
	\begin{itemize}
		\item $Q_0\in 2^{J_0}$ with $J_0=\{0\}$, 
		\item $Q,Q^\prime\in 2^J$ with $J= \{1,3,\ldots,2c-1\}$, 
		\item $P_0\in 2^{J_0}$ with $J_0=\{1,\ldots,p\}$, 
		\item $P,P^\prime\in 2^J$ with $J= \{(1,1),(1,3),\ldots,(1,2c-1),(2,1),\ldots,(p,2c-1)\}$,
	\end{itemize}
	where $ \gT_{i,j}\defi\spanR_{\scriptscriptstyle \complex}T_{i,j} $.
\end{theorem}

\begin{remark}
	\label{pizzardonecomplessi}
	{\em
		In this section, for the sake of simplicity, only the case of possibly complex minimum-phase invariant zeros has been considered. The same machinery can easily be employed to tackle the case where some freely assignable closed-loop eigenvalues are selected to be complex (in complex conjugate pairs). Addressing the general case where some pairs of assignable eigenvalues are chosen to be complex  conjugate involves a full characterization of the order of the indexing of the assigned eigenvalues as already done in the indexing of the invariant zeros. This minor extension does not lead to an augmentation of the set of solvable problems; indeed, if a problem is solvable by assigning complex conjugate eigenvalues which are not invariant zeros, it is always solvable by assigning real closed-loop eigenvalues. This is clearly not the case for the minimum-phase invariant zeros, which cannot be selected; this is the reason why this case has been considered in this section.
		% in the definition of the sets $ 
	}
\end{remark}

\subsection{Necessary conditions}%\
%\\
An important consideration is related to the necessary solvability conditions in the  presence of complex conjugate closed-loop modes. Computing the necessary and sufficient conditions provided in this section could result in cumbersome calculations. Hence, the user may prefer to have algorithmically less burdensome necessary condition to check before considering going through the necessary ad sufficient ones.
We show here that the conditions provided in Sections \ref{sec:problem1_real}-\ref{sec:problem3_real} in this case result to be exactly the necessary condition we were looking for. for the sake of brevity, we only address Problem \ref{pro01}B. All the other cases can be treated using the same machinery.
\begin{theorem}
	Let $\nu_0=n-\nu_1-\ldots-\nu_p$. 
	If Problem \ref{pro01}B is solvable then 
	\be
	\label{cond1B_compII_1}
	\dim \left(\gV^\star_g+ \sum_{(i,j)\in P} \gR_{i}(\lambda_{i,j})\right) \ge \operatorname{card}\,P+\nu_0,
	\ee
	and
	\be
	\label{cond1B_compII_2}
	\dim \left(\sum_{(i,j)\in P} \gR_{i}(\lambda_{i,j})\right) \ge \operatorname{card}\,P,
	\ee
	for all $P$ in the power set $2^I$ where $I=\{(1,1),\ldots,(1,\nu_1),\ldots,(p,1),\ldots,(p,\nu_p)\}$.
\end{theorem}
\proof
If the problem is solvable, then \eqref{cond1B_comp} holds true. We first note that for each pair of complex conjugate subspaces $ \gE_{g,i},\,\overline{\gE}_{g,i} $, with $ i\in J =\{1,3\ldots,2c-1\}$, we can find a pair of complex conjugate basis matrices $ A_{g,i} $ and $ \overline{A}_{g,i} $ such that\footnote{Given a real or complex matrix $M$, we denote by $\spanR_{\scriptscriptstyle \complex}\{M\}$ the span of the columns of $M$ over the field $\complex$.}
\[ 
\spanR_{\scriptscriptstyle \complex} E_{g,i}+\spanR_{\scriptscriptstyle \complex}\overline{E}_{g,i}=\gE_{g,i}+\overline{\gE}_{g,i}=\spanR_{\scriptscriptstyle \complex} \{A_{g,i}\}+\spanR_{\scriptscriptstyle \complex}\{\overline{A}_{g,i}\}=
\spanR_{\scriptscriptstyle \complex}\{\bmat{cc} A_{g,i} & \overline{A}_{g,i}\emat\}.
\] 
Since $ A_{g,i} $ and $ \overline{A}_{g,i} $ are complex conjugate, it is always possible to find a complex invertible matrix $ T $ such that  
$ \tilde A_{g,i} = \bmat{cc} A_{g,i} & \overline{A}_{g,i}\emat T $ is real and $ \spanR_{\scriptscriptstyle \complex}\{\bmat{cc} A_{g,i} & \overline{A}_{g,i}\emat\}= \spanR_{\scriptscriptstyle \complex} \{\tilde A_{g,i}\}$. Defining the set $ \tilde E_{g,i} \subset \real^n $ as the set that comprises all the columns of $ \tilde A_{g,i} $, there holds     
\[
\spanR_{\scriptscriptstyle \complex} E_{g,i}+\spanR_{\scriptscriptstyle \complex}\overline{E}_{g,i}=\gE_{g,i}+\overline{\gE}_{g,i}=\spanR_{\scriptscriptstyle \complex}\tilde E_{g,i}.
\]
Moreover, for every pair of complex conjugate sets $ E_{g,i},\,\overline{E}_{g,i} $, with $ i\in J $, if $ \dim \bigl( \spanR_{\scriptscriptstyle \complex} E_{g,i}+\spanR_{\scriptscriptstyle \complex}\overline{E}_{g,i}\bigr) \geq 2n $ for some $ n\in \mathbb{N} $, then $ \dim \bigl(\spanR_{\scriptscriptstyle \complex}E_{g,i}\bigr) =\dim \bigl(\spanR_{\scriptscriptstyle \complex}\overline{E}_{g,i}\bigr)\geq n  $.\\ 
%Now, We need to prove the equivalence between the hypotheses of then theorem to be proven and Lemma \ref{leFE}. In view of the previous considerations, 

Now, \eqref{cond1B_comp} can be  rewritten as  
\[
\dim \left(\sum_{i\in Q_0}\spanR_{\scriptscriptstyle \complex}E_{g,0}+\sum_{i\in Q}\spanR_{\scriptscriptstyle \complex}\tilde E_{g,i} +\sum_{(i,j)\in P} \spanR_{\scriptscriptstyle  \complex}\widehat R_i(\lambda_{i,j})\right) \ge 
\sum_{i\in Q_0} \nu_{0,i}+\sum_{i\in Q} 2\nu_{0,i}+\operatorname{card} P.
\]
The previous equation can be conveniently rewritten as 
\[
\dim \left(\spanR_{\scriptscriptstyle \complex}\Big(\bigcup_{i\in Q_0}E_{g,0}\cup \bigcup_{i\in Q}\tilde E_{g,i} \cup  \bigcup_{(i,j)\in P} \widehat R_i(\lambda_{i,j})\Big)\right) \ge 
\sum_{i\in Q_0} \nu_{0,i}+\sum_{i\in Q} 2\nu_{0,i}+\operatorname{card} P.
\]
Since $ \lambda_{i,j} \in \real $ and %\cancel{$ E_{0,0} \subseteq \real $}
in view of Lemma \ref{lemss}, all the sets appearing in the left hand-side of the latter are real, hence, $ \dim \bigl(\spanR_{\scriptscriptstyle \complex}\{\cdot \}\bigr) =  \dim \bigl(\spanR_{\scriptscriptstyle \real}\{\cdot \} \bigr)$ and we can rewrite
\[
\dim \left(\spanR_{\scriptscriptstyle \real}\Big(\bigcup_{i\in Q_0} E_{g,0}\cup\bigcup_{i\in Q}\tilde E_{g,i}\Big) + \sum_{(i,j)\in P} \spanR_{\scriptscriptstyle  \real} \widehat R_i(\lambda_{i,j})\right) \ge 
\sum_{i\in Q_0} \nu_{0,i}+\sum_{i\in Q} 2\nu_{0,i}+\operatorname{card} P.
\]     
When $ Q_0=\{0\} $ and $ Q=\{1,3\ldots,2c-1\} $ the previous condition is easily seen to be equivalent to \eqref{cond1B_compII_1}, whereas when
$ Q_0 $ and $ Q $ are empty sets, the equivalence with \eqref{cond1B_compII_2} is proven.
\endproof

Similar necessary conditions can be obtained for the other problems considered in this paper, following the same ideas.

\section*{Concluding remarks}
In this paper, we provided necessary and sufficient constructive conditions for the solution of the general eigenstructure assignment problem, which is shown to be equivalent to a tracking problem in which 
a certain number of closed-loop modes appear in each output component. 

This problem is not just important {\em per se}, but also because in the past twenty years it appeared as the prototype of a variety of non-interacting and fault detection problems, for which a set of necessary and sufficient conditions could only be achieved {\em a posteriori} by checking the rank of the matrix of closed-loop eigenvectors. 

Nine problems have been identified in this paper, whose formulation depends on whether the eigenvalues to be assigned coincide or not with invariant zeros of the system, on the fact that we may want to assign only the number, but not the specific numerical value, of the closed-loop modes, and also on whether we want this assignment to take place only within the unobservable, or also in the observable part of the closed-loop spectrum.

The solvability conditions of these problems have been obtained by merging the key ideas of combinatorics with those of geometric control theory. The method for determining the decoupling filter matrix is also outlined.
The new framework developed in this paper has yielded a satisfactory answer to control/estimation problems for 
which, so far, the use alone of standard geometric techniques has not been successful. We expect the same framework to provide important insight into problems that are still open in control theory, such as the input-output (row-by-row) decoupling problem.

\end{document}